\documentclass[a4paper,11pt]{article}
\usepackage[utf8]{inputenc}
\usepackage[english]{babel}
\usepackage{amsthm, amsfonts,amsmath,amscd,amssymb,epsfig,stmaryrd,psfrag, mathtools}

\usepackage{enumitem}
\usepackage{enumerate,bbm,array}
\usepackage{tikz} %For the color command
\usepackage{yhmath} 

\usepackage[active]{srcltx}
\usepackage[T1]{fontenc}
\usepackage{graphicx}
\mathtoolsset{showonlyrefs}

\usepackage{rotating}

\usepackage{xcolor}
\usepackage{framed}

\usepackage{pifont} %flowers

\usepackage{stmaryrd}
\usepackage{nicefrac}

\usepackage{comment}
\usepackage{accents}

\usepackage{longtable,geometry}

\usepackage{hyperref}
\hypersetup{
	colorlinks=true, %set true if you want colored links
	linktoc=all,     %set to all if you want both sections and subsections linked
	linkcolor=blue  %choose some color if you want links to stand out
}

\usepackage{centernot}

\title{Rotational invariance in critical planar lattice models}
\author{
Hugo Duminil-Copin\thanks{Institut des Hautes \'Etudes Scientifiques and Universit\'e Paris-Saclay}\ \thanks{Universit\'e de Gen\`eve}, 
Karol Kajetan Kozlowski\addtocounter{footnote}{0}\thanks{ENS Lyon}, 
Dmitry Krachun\addtocounter{footnote}{0}\thanks{Princeton University}, \\
Ioan Manolescu\addtocounter{footnote}{0}\thanks{University of Fribourg}, 
Mendes Oulamara\addtocounter{footnote}{-4}\footnotemark}

\theoremstyle{plain} \newtheorem{theorem}{Theorem}[section]
\theoremstyle{plain} \newtheorem{corollary}[theorem]{Corollary}
\theoremstyle{plain} \newtheorem{lemma}[theorem]{Lemma}
\theoremstyle{plain} \newtheorem{proposition}[theorem]{Proposition}
\theoremstyle{plain} \newtheorem{claim}[theorem]{Claim}
\theoremstyle{remark}\newtheorem{remark}[theorem]{Remark}
\theoremstyle{definition} \newtheorem{definition}[theorem]{Definition}

\usepackage{mdframed}
\mdfdefinestyle{MyFrame}{%
    linecolor=black,
    outerlinewidth=3pt,
    innertopmargin=\baselineskip,
    innerbottommargin=\baselineskip,
    innerrightmargin=10pt,
    innerleftmargin=10pt,
    backgroundcolor=gray!10!white}

\newcommand{\calC}{\mathcal{C}}

\newcommand{\calE}{\mathcal{E}}
\newcommand{\calF}{\mathcal{F}}

\newcommand{\calH}{\mathcal{H}}

\newcommand{\calL}{\mathcal{L}}

\newcommand{\calQ}{\mathcal{Q}}
\newcommand{\calR}{\mathcal{R}}

\newcommand{\calW}{\mathcal{W}}

\newcommand{\bfS}{\mathbf{S}}
\newcommand{\bfT}{\mathbf{T}}

\newcommand{\bbC}{\mathbb{C}}

\newcommand{\bbE}{\mathbb{E}}

\newcommand{\bbG}{\mathbb{G}}

\newcommand{\bbL}{\mathbb{L}}

\newcommand{\bbN}{\mathbb{N}}

\newcommand{\bbP}{\mathbb{P}}

\newcommand{\bbR}{\mathbb{R}}

\newcommand{\bbT}{\mathbb{T}}

\newcommand{\bbV}{\mathbb{V}}

\newcommand{\bbZ}{\mathbb{Z}}
\newcommand{\sfn}{{\sf n}}

\newcommand{\sfC}{{\sf C}}
\newcommand{\sfD}{{\sf D}}
\newcommand{\sfN}{{\sf N}}

\newcommand{\ep}{\varepsilon}

\newcommand{\rk}[1]{\bgroup\color{red}%
	\par\medskip\hrule\smallskip%
	\noindent\textbf{#1}%
	\par\smallskip\hrule\medskip\egroup}

\newcommand{\Top}{{\rm T}}
\newcommand{\Bottom}{{\rm B}}
\newcommand{\Left}{{\rm L}}
\newcommand{\Right}{{\rm R}}
\newcommand{\eps}{\epsilon}

\newcommand\xlra{\xleftrightarrow}
\renewcommand{\comment}[1]{}

\newcommand{\ind}{\mathbbm{1}}

\begin{document}
\maketitle

\begin{abstract}
We prove that the large-scale properties of a number of two-dimensional lattice models are rotationally invariant. More precisely, we prove that the random-cluster model on the square lattice with cluster-weight~$1\le q\le 4$ exhibits rotational invariance at large scales. This covers the case of Bernoulli percolation on the square lattice as an important example. We deduce that the correlations of the critical Potts models with~$q\in\{2,3,4\}$ colours are rotationally invariant at large scales. 

Our result is instrumental in proving the convergence of the six-vertex model to the Gaussian Free Field in a separate paper.
\end{abstract}

\tableofcontents

\section{Introduction}

\subsection{Motivation}

Physical systems undergoing a continuous phase transition have been the focus of much attention in the past seventy years, both on the physical and the mathematical sides. Since Onsager's revolutionary solution of the 2D Ising model, mathematicians and physicists have tried to understand the delicate features of the critical phase of these systems. In the sixties, the arrival of the renormalization group (RG) formalism (see \cite{Fis98} for a historical
exposition) led to a generic (non-rigorous) deep understanding of continuous phase transitions. The RG formalism suggests that ``coarse-graining'' renormalization transformations correspond to appropriately changing the scale and the parameters of the model under study. The large-scale limit of the critical regime then arises as the fixed
point of the renormalization transformations. 

A striking consequence of the RG formalism is that the assumption that the critical fixed point is unique leads one to the prediction that  the scaling limit at the critical point must satisfy translation, rotation
and scale invariance, which allows one to deduce some information about correlations.  
In \cite{PolyakovArgumentAboutCFTInvarianceCriticalCorrelators}, Polyakov outlined a set of arguments pointing towards  a much stronger invariance of statistical physics models at
criticality: since the scaling limit field theory is a local field, it should be invariant
under any map which is locally a composition of translation, rotation and homothety,
which leads to postulate full conformal invariance. In \cite{BPZ84a,BPZ84b}, Belavin, Polyakov and Zamolodchikov went even further by considering massless field theories that enjoy full conformal invariance from the outset, which allowed them to derive explicit expressions for their correlation functions, hence giving birth to conformal field theories.
Once conformal invariance is proved, a whole world of new techniques becomes available thanks to Conformal Field Theory and the Schramm-Loewner Evolution \cite{Law05}, and it is therefore a problem of fundamental importance to prove conformal invariance of the scaling limits of lattice models.

Proving conformal invariance is quite difficult  for most lattice models. The examples of models for which such a statement has been obtained can be counted on the fingers of one's hand: 
site Bernoulli percolation on the triangular lattice \cite{Smi01,KhrSmi21,CamNew07,CamNew06} (respectively for Cardy's formula for the first two, SLE(6) convergence, and CLE(6) convergence), 
Ising and FK-Ising models \cite{Smi10,CheSmi12,HonSmi13,CheHonIzy15,CheDumHon12a,BenHon16} (respectively for the fermionic observables in FK-Ising, in Ising, the energy and the spin fields, SLE convergence, and CLE convergence), uniform spanning trees \cite{LawSchWer04}, dimers \cite{Ken00}, level lines of the discrete GFF \cite{SchShe09}. 

In all the cases mentioned, the proof relied, in one form or another, on discrete holomorphic observables satisfying some discrete version of conformally covariant boundary value problems. 
Mathematicians were therefore able to prove conformal invariance directly, bypassing the road suggested by physicists consisting in first proving scaling and rotation invariance (translation invariance is obvious), and then deducing from it conformal invariance. Unfortunately, the current mathematical strategy is very dependent on discrete properties of the system, which explains why we are currently limited to very few instances of proofs of conformal invariance.

In this paper, we perform one step towards the strategy inspired by field theory 
and prove rotational invariance of the large-scale properties of a family of planar models at their critical point. 
Our strategy is quite general and applies to a number of integrable planar systems. 
We treat here the case of the random-cluster model (also called Fortuin-Kasteleyn percolation),
and deduce properties for the associated Potts model using the known mapping between the two. 

While a previous (unpublished) version of this paper used an eigenvalue computation specifically related to the random-cluster model, the present argument only relies on the Yang-Baxter equation (or star-triangle transformation) satisfied by the random-cluster model in the isoradial setting. As such, it indicates that {\em planar models satisfying the Yang-Baxter equation are rotationally invariant}. In particular, we expect the argument to extend to a large variety of other settings such as the Ashkin-Teller model and certain loop models. 

The result of this paper is an important pillar of the proof in \cite{DumKozLamMan26} of the convergence of the full-plane six-vertex model to the Gaussian Free Field (GFF) in the range of parameters corresponding to the critical random-cluster model with $1 \leq q \leq 4$. 
This proof of convergence thus follows the heuristic of  \cite{PolyakovArgumentAboutCFTInvarianceCriticalCorrelators}, with the present paper providing an essential input about the emerging symmetries of the scaling limit.

% We believe that the reasoning has also applications for the Askhin-Teller model and certain loop models. 
% The proof will proceed by focusing on the random-cluster model and then extending its rotational invariance to other planar lattice models using known mapping between the models. 

\subsection{Definitions: random-cluster model and large-scale topology}\label{sec:1.2}

As mentioned in the previous section, the model of central interest in this paper is the random-cluster model, 
 introduced by Fortuin and Kasteleyn around 1970 \cite{ForKas72,For70}, which we now define. 
For background, we direct the reader to the monograph \cite{Gri06} and to the lecture notes \cite{Dum17a} for an exposition of more recent results.

Consider the square lattice~$(\mathbb Z^2,\mathbb E)$, that is the graph with vertex-set~$\mathbb Z^2=\{(n,m):n,m\in\mathbb Z\}$ and edges between nearest neighbours. In a slight abuse of notation, we write~$\mathbb Z^2$ for the graph itself. 
Consider a finite subgraph~$G$ of the square lattice with vertex-set~$V$ and edge-set~$E$. For instance, think of~$G=\Lambda_n$ as being the subgraph of~$\mathbb Z^2$ spanned by the vertex-set~$\{-n,\dots, n\}^2$ (we will use the notation~$\Lambda_n$ throughout the paper). A percolation configuration~$\omega$ on~$G$ is an element of~$\{0,1\}^{E}$. An edge~$e$ is {\em open} (in~$\omega$) if~$\omega_e=1$, otherwise it is {\em closed}. A configuration $\omega$ can be seen as a subgraph of~$G$ with vertex-set~$V$ and edge-set~$\{e\in E:\omega_e=1\}$. When speaking of connections in~$\omega$, we view~$\omega$ as a graph. A {\em cluster} is a connected component of~$\omega$. 

%The {\em boundary conditions}~$\xi$ on~$G$ are given by a partition of the set~$\partial G$ of vertices in~$V$ incident to at most three edges in~$E$. We say that two vertices of~$G$ are {\em wired together} if they belong to the same element of the partition~$\xi$.
\begin{definition} The random-cluster measure on~$G$ with edge-weight~$p\in[0,1]$, cluster-weight~$q>0$, and free boundary conditions   is given by
\begin{align}\label{eq:RCM_def1}
	\phi_{G,p,q}^0[\omega]:=\frac{p^{|\omega|}(1-p)^{|E|-|\omega|}q^{k(\omega)}}{Z^0_{\mathrm{RCM}}(G,p,q)} ,
\end{align}
where~$|\omega|:=\sum_{e\in E}\omega_e$ is the number of open edges,~$k(\omega)$ is the number of connected components of the graph, and~$Z^0_{\rm RCM}(G,p,q)$ is a normalising constant called the {\em partition function} chosen in such a way that~$\phi_{G,p,q}^0$ is a probability measure. 
\end{definition}

For~$q\ge1$, the family of measures~$\phi_{G,p,q}^0$ converges weakly as~$G$ tends to the whole square lattice to an infinite-volume measure~$\phi_{p,q}^0$ on~$\{0,1\}^{\mathbb E}$. The random-cluster model undergoes a phase transition \cite{BefDum12,DumRaoTas17} at a critical parameter 
\[
p_c=p_c(q)=\frac{\sqrt q}{1+\sqrt q}
\]
in the sense that the~$\phi_{p,q}^0$-probability that there exists an infinite cluster is 0 if~$p<p_c(q)$, and is 1 if~$p>p_c(q)$. 

It was also proved in \cite{DumGagHarManTas20,DumSidTas16}  that the phase transition is continuous (i.e.~that the probability that 0 is connected to infinity tends to 0 as~$p\searrow p_c$) if and only if~$q\le 4$ (see also \cite{GlaLam25} and \cite{RaySpi19} for alternative proofs of continuity when $q \leq 4$ and discontinuity for $q >4$, respectively). In the whole paper we restrict our attention to the range~$q\in[1,4]$. For this reason,
\begin{center}
{\em fix~$q\in[1,4]$ and~$p=p_c(q)$ and drop them from notation.}
\end{center}

We will be interested in measuring how close the large-scale properties of two random percolation configurations really are. In order to do that, we introduce a rescaling of the lattice and define the random-cluster model on subgraphs of~$\delta\mathbb Z^2$ with~$\delta>0$. To highlight on which lattice we are working, we will consistently use the subscript~$\delta$ to refer to a percolation configuration on a subgraph of the lattice~$\delta\mathbb  Z^2$, and write~$\omega_\delta$ for such a configuration. When~$\Omega$ is a simply connected domain of~$\mathbb R^2$, write~$\Omega_\delta$ for the intersection of~$\Omega$ with~$\delta \mathbb  Z^2$.

In \cite{CamNew06}, Camia and Newman introduced a convenient way of measuring the geometry of large clusters in a percolation configuration  in the plane. Let~$\mathfrak C=\mathfrak C(\Omega)$ be the collection of sets~$\calF=\calF_0\sqcup\calF_1$  of two locally finite families~$\calF_0$ and~$\calF_1$ of non-self-crossing loops in some simply connected domain~$\Omega$ that do not intersect each other (even between loops in~$\calF_0$ and~$\calF_1$). Define the metric on~$\mathfrak C$, 
\[
d_{\rm CN}(\calF,\calF')\le \ep \ \Longleftrightarrow\ \Big(\begin{array}{c} \forall i\in\{0,1\},\forall \gamma\in \calF_i\text{ with }\gamma\subset B(0,1/\ep),\exists \gamma'\in \calF
_i',d(\gamma,\gamma')\le \ep\\
\text{ and similarly when exchanging~$\calF'$ and~$\calF$}\end{array}\Big)
\]
where, for two loops~$\gamma_1$ and~$\gamma_2$, we set 
\begin{align}\label{eq:loop_dist}
	d(\gamma_1,\gamma_2):=\inf\sup_{t\in\mathbb S^1}|\gamma_1(t)-\gamma_2(t)|,
\end{align}
with the infimum running over all continuous one-to-one parametrizations of the loops~$\gamma_1$ and~$\gamma_2$ by~$\mathbb S^1$.

The distance~$d_{\rm CN}$ induces a distance on the space of probability measures on~$\mathfrak C$, which we also denote by~$d_{\rm CN}$:
for~$\phi$ and~$\phi'$ two measures on~$\mathfrak C$, 
%let~$d_{\rm CN}(\phi,\phi')$ be the minimal~$\epsilon \geq 0$ such that there exists a coupling~$\bbP$ of~$\phi$ and~$\phi'$ with
\begin{align}\label{eq:dis_law}
d_{\rm CN}(\phi,\phi') := \inf\big\{\ep \text{ so that } \exists \text{ coupling~$\bbP$ of~$\phi$ and~$\phi'$ with }	\mathbb P[d_{\rm CN}(\omega,\omega')> \ep] <\ep\big\}.
%\smallskip
\end{align}

Another way of encoding the geometry of large clusters was proposed by Schramm and Smirnov in \cite{SchSmi11}. In order to define it formally, let a {\em quad}~$Q$ be the image of a  homeomorphism from~$[0,1]^2$ to~$\mathbb C$, and let~$a,b,c,d$ be the images of the corners of~$[0,1]^2$. A {\em crossing} of~$Q$ is a continuous path in~$Q$ going from~$(ab)$ to~$(cd)$. 
 Let~$\mathcal Q$ be the set of quads, endowed with the distance between quads given by
 \[
 d_{\mathcal Q}(Q,Q'):=d(\partial Q,\partial Q')+|a-a'|+|b-b'|+|c-c'|+|d-d'|.
 \]
 Call $S\subset\mathcal Q$ {\em hereditary} if whenever~$Q\in S$, every quad~$Q'$ such that any crossing of~$Q$ contains a crossing of~$Q'$ also belongs to~$S$. Let~$\mathfrak H=\mathfrak H(\Omega)$ be the set of closed hereditary subsets of~$\mathcal Q$. Endow~$\mathfrak H$ with the smallest topology generated by the sets of the type~$\{S  \in \mathfrak H: Q\in S\}_{Q\in\mathcal Q}$ and~$\{S  \in \mathfrak H: S\cap U=\emptyset\}_{U\text{ open set in~$\mathcal Q$}}$. The set~$\mathfrak H$ with this topology is metrizable, and we denote the metric (whose definition is implicit) by~$d_{\rm SS}(\cdot,\cdot)$.
 As for~$d_{\rm CN}$, the distance~$d_{\rm SS}$ extends to probability measures on~$\mathfrak H$ via a formula similar to~\eqref{eq:dis_law}. 
\medbreak

A configuration~$\omega$ can be identified with the (automatically hereditary) set~$S\in \mathfrak H$ containing all the quads that are crossed by an open path in~$\omega$ (seen as a continuous path in the plane). Similarly,~$\omega$ can be seen as an element of~$\mathfrak C$ by considering the loop representation of the model obtained as follows (see Section~\ref{sec:elementary} for details): to each~$\omega$ is associated a {\em dual configuration}~$\omega^*$ on the dual graph, as well as a {\em loop configuration}~$\overline\omega$ on the medial graph, corresponding basically to the boundaries between the primal and dual clusters. Then, we say that a loop is in~$\calF_1$ if it is the exterior boundary of a primal cluster, and in~$\calF
_0$ if it is the exterior boundary of a dual cluster. 
Whether~$\omega$ is seen as an element of~$\mathfrak H$ or~$\mathfrak C$ will depend on the context (it will always be clear which identification is used, if any).

%For future reference, we introduce the distances~$\mathbf d_\#$ for~$\#={\bf CC}$ and~${\bf SS}$ (we will also use the same for the homotopy distance~$d_{\bf H}$ introduced below) between two measures~$\mu,\nu$ defined as follows:
%\begin{align}
%\mathbf d_\#(\mu,\nu)< \ep \ \Longleftrightarrow\ \exists \text{ a coupling }\mathbb P\text{ of }\omega\sim\mu\text{ and }\omega'\sim\nu\text{ such that }\mathbb P[d_\#(\omega,\omega')>\ep]<\ep.
%\end{align}

\subsection{Rotational invariance of the random-cluster model}

The main theorem of our paper is the following.

\begin{theorem}[Rotation invariance of critical random-cluster model]\label{thm:rotation_invariance_infinite_vol}
	Fix~$q\in[1,4]$. There exist constants $c,C > 0$ such that the following holds. 
	For~an angle $\alpha \in [0,2\pi]$ write~$\phi_{e^{i\alpha}\delta \bbZ^2}$
	for the critical random-cluster measure on the rescaled lattice~$\delta \bbZ^2$, rotated by~$\alpha$. 
	Then, for any~$\alpha \in [0,2\pi]$ and~$\delta >0$,
   \begin{align*}
	{d}_{\rm CN}(\phi_{\delta \bbZ^2},\phi_{e^{i\alpha}\delta \bbZ^2}) \leq C \delta^c.
    \end{align*}
    Furthermore, we also have that 
    \begin{align*}
	{d}_{\rm SS}(\phi_{\delta \bbZ^2},\phi_{e^{i\alpha}\delta \bbZ^2}) \xrightarrow[\delta \to 0]{} 0.
    \end{align*}
\end{theorem}
In light of the above, we say that~$\phi_{ \bbZ^2}$ is {\em asymptotically rotationally invariant}. It implies that any scaling limit is rotationally invariant.

A similar result holds in simply connected domains~$\Omega$ with a~$C^1$-smooth boundary, meaning that~$\partial\Omega$ can be parametrized by a~$C^1$-function whose differential does not vanish at any point.

\begin{corollary}[Rotation invariance in domains]\label{cor:rotation_invariance_domain}
	Fix~$q\in[1,4]$ and a simply connected domain~$\Omega$ with a~$C^1$-smooth boundary. 
	For~$\alpha \in [0,2\pi]$ write~$\phi_{\Omega_\delta}^0$ for the critical random-cluster model on a discretisation of $\Omega \cap   \delta \bbZ^2$ with free boundary conditions,
	and $\phi_{(e^{{\rm i}\alpha} \Omega)_\delta}^0 \circ e^{{\rm i}\alpha}$ for the measure obtained in a similar way but using the rotated lattice $e^{i\alpha}\delta \bbZ^2$.
	Then,    \begin{align*}
	{d}_{\rm CN}\big(\phi_{\Omega_\delta}^0,\phi_{(e^{{\rm i}\alpha} \Omega)_\delta}^0 \circ e^{{\rm i}\alpha}\big)&\xrightarrow[\delta \to 0]{}0,\\
	{d}_{\rm SS}\big(\phi_{\Omega_\delta}^0,\phi_{(e^{{\rm i}\alpha} \Omega)_\delta}^0 \circ e^{{\rm i}\alpha}\big)&\xrightarrow[\delta \to 0]{}0. 
    \end{align*}
\end{corollary}

Such a $C^1$ condition on domains may be relaxed to cover any Jordan domain, yet we avoid such considerations to focus on the most interesting aspects of the problem at hand (which are already encompassed in the present framework).
For a similar reason, we state the result in a non-quantitative form, as a quantitative version would require a substantially more technical proof.

Theorem~\ref{thm:rotation_invariance_infinite_vol} has a number of applications for the random-cluster model. First, the definition of the Schramm-Smirnov topology implies, in particular, that crossing probabilities are invariant under rotation in the following sense. For a quad~$Q$, let~$\{\omega\in\mathcal C(Q)\}$ be the event that~$Q$ is crossed in the percolation configuration~$\omega$.

\begin{corollary}[Rotation invariance of crossing probabilities]\label{cor:crossing}
    Fix~$q\in[1,4]$, a simply connected domain~$\Omega$ with a~$C^1$-smooth boundary
    and~$Q$ a quad contained in~$\Omega$.
    Then, for each~$\alpha \in [0,\pi]$, 
    \begin{align}\label{eq:crossing}
    	\phi_{(e^{{\rm i}\alpha}\Omega)_\delta}^0[\mathcal C(e^{{\rm i}\alpha} Q)]-\phi_{\Omega_\delta}^0[\mathcal C(Q)] \xrightarrow[\delta \to 0]{} 0.
    \end{align}
    Furthermore, for~$1 \leq q < 4$, the condition~$Q \subset \Omega$  may be replaced 
    with the more natural condition that~$Q \subset \overline\Omega$. 
\end{corollary}

The convergence in Corollary \ref{cor:crossing} may be proved to be uniform in~$\alpha \in [0,\pi]$. 
The second part of the corollary is restricted to $q < 4$ due to the particular way free boundary conditions interact with the existence of open paths touching the boundary when $q = 4$. 
For instance, it is expected that, when $\Omega = Q$ and $q = 4$, both quantities in~\eqref{eq:crossing} converge to zero as~$\delta$ tends to 0. %\im{TO DO: add reference to $CLE_4$ paper when ready}. 

We  turn to ``pointwise correlations''. For points~$x_1,\dots,x_n$ and a partition~$\mathcal P$ of~the set $\{x_1,\dots,x_n\}$, let~$\mathcal E(\mathcal P,x_1,\dots,x_n)$ be the event that~$x_i$ and~$x_j$ are connected if and only if they belong to the same element of~$\mathcal P$. The following corollary will be useful when studying spin-spin correlations in the Potts model. 

\begin{corollary}[Rotation invariance of connectivity correlations]\label{cor:connectivity_correlations}
	Fix~$q\in[1,4]$ and a simply connected domain~$\Omega$ with a~$C^1$-smooth boundary. 
	For every~$\ep>0$ and~$n$, there exists~$\delta_0=\delta_0(q,n,\ep,\Omega)>0$ 
	such that for every~$\alpha\in[0,2\pi]$ and~$\delta\le \delta_0$, every~$x_1,\dots,x_n\in\Omega_\delta$ 
	at a distance at least~$\ep$ from each other and from the boundary of~$\Omega$, and every partition~$\mathcal P$ of~$\{x_1,\dots,x_n\}$, 
	\[
		|\phi^0_{(e^{{\rm i}\alpha}\Omega)_\delta}[\mathcal E(\mathcal P,e^{{\rm i}\alpha} x_1,\dots,e^{{\rm i}\alpha} x_n)]-\phi^0_{\Omega_\delta}[\mathcal E(\mathcal P,x_1,\dots,x_n)]|
		\le \ep\, \phi^0_{\Omega_\delta}[\mathcal E(\mathcal P,x_1,\dots,x_n)],
	\]
	where we use, in a slight abuse of notation,~$e^{{\rm i}\alpha} x_i$ to denote a vertex~$x$ of~$(e^{{\rm i}\alpha}\Omega)_\delta$ 
	within a distance~$\delta$ from the image of~$x_i$ under the rotation by the angle~$\alpha$.
\end{corollary}

\begin{remark}
	We may also study the edge-density variables~$\epsilon_e^\Omega:=\omega_e-\phi_\Omega^0[\omega_e]$ and attempt to prove some rotation invariance for these variables. 
	Obtaining this result requires combining the present result with coupling techniques developed in \cite{DumMan20}. For the sake of brevity, we do not consider these observables here.
	% \im{TO DO  add reference}
\end{remark}

\subsection{Universality  on isoradial rectangular graphs}\label{sec:results_universality}

The previous rotational invariance result (Theorem~\ref{thm:rotation_invariance_infinite_vol}) is accompanied by a universality result for certain isoradial lattices; see Theorem~\ref{thm:universalCNSS} below. In order to state it, we describe now an inhomogeneous random-cluster model on some distorted embedding of the square lattice~$\bbZ^2$. 
These notions may appear strange at first, but will be shown to fit in the more general framework of isoradial graphs (see Section~\ref{sec:2}). Fix $q \in [1,4]$.

\begin{figure}
\begin{center}
\includegraphics[width=0.45\textwidth,page=1]{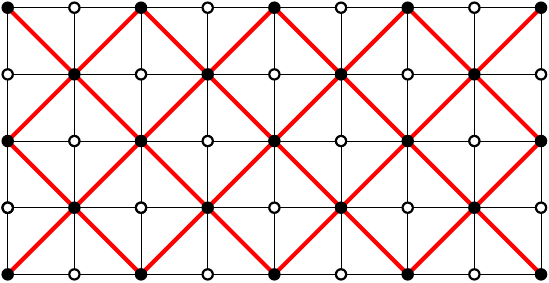}\hspace{.05\textwidth}
\includegraphics[width=0.45\textwidth,page=2]{figures/isoradial-lattice-pure.pdf}
\caption{In red: the lattice $\bbL(\pi/2)$ on the left, and a lattice $\bbL(\alpha)$ on the right. 
The dual vertices are hollow; the black graph is the diamond graph, introduced in Section~\ref{sec:isoradial_graphs}. }
\label{fig:isoradial-lattice-pure}
\end{center}
\end{figure}

For~$\alpha \in (0,\pi)$, let~$\bbL(\alpha)$ be the embedding of~$\bbZ^2$ in which horizontal edges have length~$2 \cos(\alpha/2)$ 
and vertical edges have length~$2\sin (\alpha/2)$, rotated by an angle~$\alpha/2$; see Figure~\ref{fig:isoradial-lattice-pure}. 
Consider the random-cluster model on~$\bbL(\alpha)$ with cluster-weight~$q$ and different edge-parameters~$p_{\rm hor}$ and~$p_{\rm vert}$ for the ``horizontal'' and ``vertical'' edges (that is, the edges of lengths~$2 \cos(\alpha/2)$ and~$2\sin (\alpha/2)$) respectively, given by
\begin{align}
	\frac{p_{\rm hor}}{1-p_{\rm hor}} = q \frac{1-p_{\rm vert}}{p_{\rm vert}}  = \sqrt{q} \frac{\sin( r \alpha)}{\sin( r (\pi-\alpha))}
	\qquad \text{ with } r := \tfrac{1}{\pi} \cos^{-1} \left( \tfrac{\sqrt{q}}{2} \right).
\end{align}
Write~$\phi_{\bbL(\alpha)}$ for the unique infinite-volume measure on~$\bbL(\alpha)$ with the parameters above (the uniqueness of the measure will be discussed below).
We call~$\bbL(\alpha)$ an {\em isoradial rectangular lattice}, and~$\phi_{\bbL(\alpha)}$ its associated random-cluster model.

In the particular case of~$\alpha = \pi/2$,~$\phi_{\bbL(\pi/2)}$ is the critical random-cluster on~$\bbZ^2$, rotated by~$\pi/4$ and dilated by a factor $\sqrt 2$. 
It was proved in \cite{DumLiMan18} that~$\phi_{\bbL(\alpha)}$ shares the qualitative features of critical random-cluster on~$\bbZ^2$, most importantly an RSW property (see Theorem~\ref{thm:RSW_iso}) and the uniqueness of the infinite-volume measure. 

The models~$\phi_{\bbL(\alpha)}$ for different values of~$\alpha$ may be related via a series of star-triangle transformations (see Section~\ref{sec:star-triangle}). Features that are stable under these transformations may then be transferred from~$\phi_{\bbL(\pi/2)}$ to all other~$\phi_{\bbL(\alpha)}$. 
This strategy was used in \cite{GriMan14,DumLiMan18} for RSW estimates, and will be used here to prove the following universality result.

\begin{theorem}[Asymptotic universality for isoradial rectangular graphs]\label{thm:universalCNSS}
	Fix~$q\in[1,4]$. There exist constants $c,C > 0$ such that,  for every~$\alpha \in (0,\pi)$ and $\delta>0$,
	\begin{align}\label{eq:universalCNSS}
		d_{\rm CN}(\phi_{\delta\bbL(\alpha)}, \phi_{\delta\bbL(\frac\pi2)}) < C \delta^c.
    \end{align}
      Furthermore, we also have that 
    \begin{align*}
	d_{\rm SS}(\phi_{\delta\bbL(\alpha)}, \phi_{\delta\bbL(\frac\pi2)}) \xrightarrow[\delta \to 0]{} 0.
    \end{align*}

\end{theorem}

The statement above is a remarkable instance of universality. It states that the isoradial embedding is the exact embedding that compensates the inhomogeneity between the vertical and horizontal edge-parameters of the random-cluster measures $\phi_{\bbL(\alpha)}$, 
and produces a model that is expected to have a conformally-invariant scaling limit. 
The previously mentioned results \cite{GriMan14,DumLiMan18} do not depend on the embedding, 
and use isoradiality only as a guarantee that the star-triangle transformation applies. 
As such, Theorem~\ref{thm:universalCNSS} is more subtle than previous uses of this type of strategy. 

The universality of Theorem~\ref{thm:universalCNSS} is generalised to all bi-periodic isoradial graphs in~\cite{HanMan25}. 

\begin{remark}\label{rem:univ_implies_rot}
Theorem~\ref{thm:universalCNSS} readily implies Theorem~\ref{thm:rotation_invariance_infinite_vol}.
Indeed, for any~$\alpha \in (0,\pi)$,~$e^{{i}\alpha/2}\mathbb R$ is an axis of symmetry for~$\bbL(\alpha)$, and therefore for~$\phi_{\bbL(\alpha)}$.
Then,~\eqref{eq:universalCNSS} implies that~$\phi_{\bbL(\pi/2)}$ is  asymptotically invariant under the reflection with respect to~$e^{{i}\alpha/2}\mathbb R$. 
The composition of the reflections with respect to the horizontal axis and to~$e^{{i}\alpha/2}\mathbb R$ produces the rotation by~an angle $\alpha$, 
hence the asymptotic invariance of~$\phi_{\bbL(\pi/2)}$ under this rotation. 
\end{remark}

In light of the remark above, it would be tempting to assume that we will first prove Theorem~\ref{thm:universalCNSS}, then deduce Theorem~\ref{thm:rotation_invariance_infinite_vol}. 
This is not the case, as we explain below.

Theorem~\ref{thm:universalCNSS} is proved by gradually transforming $\bbL(\frac\pi2)$ into $\bbL(\alpha)$, and keeping track of the configuration throughout. This produces a Markov chain of configurations $\omega_0,\dots, \omega_T$ on various intermediate graphs, with 
$\omega_0 \sim \phi_{\delta\bbL(\frac\pi2)}$ and  $\omega_T \sim \phi_{\delta \bbL(\alpha)}$. Because every passage from $\omega_t$ to $\omega_{t+1}$ is the composition of star-triangle transformations (which preserve connections) the macroscopic clusters of $\omega_0$ are preserved throughout the process, but may progressively change shape. 
Two essential ingredients are needed for Theorem~\ref{thm:universalCNSS}:
\begin{itemize}
\item[(i)] clusters change shape in a "stationary" way throughout the process $(\omega_t)_{0\leq t\leq T}$, implying that the change between the clusters of $\omega_0$ and of $\omega_T$ is affine (up to polynomially diminishing errors) and
\item[(ii)] the expected move of the clusters (hereafter called {\em drift}) is $0$, which implies that the affine function mentioned above is  the identity. 
\end{itemize}

To better separate these ingredients, %, and ultimately to obtain a simpler proof of Theorem~\ref{thm:universalCNSS}, 
we start by proving the following weaker version of Theorem~\ref{thm:universalCNSS}, which encapsulates point (i) above. 

\begin{theorem}[Universality up to linear deformation]\label{thm:linear}
	For every~$q\in[1,4]$ and~$\alpha,\beta\in (0,\pi)$, 
	there exist constants~$c,C > 0$ and an invertible linear map~$M_{\beta,\alpha}:\bbR^2 \to \bbR^2$ such that 
	\begin{align}\label{eq:linear}
		d_{\rm CN}\big[ \phi_{\delta\bbL(\beta)}, \phi_{\delta\bbL(\alpha)}\circ M_{\beta,\alpha}\big] \leq C\,\delta^c 
		\quad \text{ for all~$\delta > 0$}.
	\end{align}
\end{theorem}

The statement above should be understood as~$\omega \sim \phi_{\delta\bbL(\beta)}$ and~$\omega' \sim  \phi_{\delta\bbL(\alpha)}$ may be coupled so that the loop representations of~$\omega$ and~$M_{\beta,\alpha}^{-1}(\omega')$ are close (or equivalently such that~$M_{\beta,\alpha}(\omega)$ and~$\omega'$ are close) for~$d_{\rm CN}$ with high probability. 

Remarkably, Theorem~\ref{thm:linear} only uses qualitative features of the critical random-cluster and the star-triangle transformation.
Ignoring the dependence of $c,C$ on $\alpha$ and $\beta$, 
Theorem~\ref{thm:universalCNSS} is equivalent to Theorem~\ref{thm:linear} with the additional input that~$M_{\pi/2,\alpha} = {\rm id}$.
Deducing that $M_{\pi/2,\alpha} = {\rm id}$ requires some form of exact integrability that encodes the relationship between the model and the isoradial embedding.

A direct way to prove that~$M_{\pi/2,\alpha} = {\rm id}$ is via explicit computations of certain partition functions of the six-vertex model based on the Bethe ansatz (for instance relying on \cite{DumKozKraManTik20}). These are then translated into probabilities of events for the random-cluster model on lattices $\bbL(\alpha)$, 
which are then used to prove that the drift of the transformations above is null. 
This approach is technically challenging and uses an advanced form of exact integrability that is ultimately not necessary. It was used in an early version of the present paper and may be found in \cite{MendesThesis}. 

We favour another strategy relying on the symmetries of the different lattices $\mathbb L(\alpha)$. This argument avoids the use of any form of exact integrability other than the star-triangle transformation. 
More precisely, once Theorem~\ref{thm:linear} is proved, we show the asymptotic rotational invariance of~$\phi_{\bbL(\pi/2)}$, i.e.~Theorem~\ref{thm:rotation_invariance_infinite_vol}, by an argument similar to the one described in Remark~\ref{rem:univ_implies_rot}; see Section~\ref{sec:deducing_rot_inv}.
Then, in Section~\ref{sec:drift0}, we use  Theorem~\ref{thm:rotation_invariance_infinite_vol} to prove that~$M_{\pi/2,\alpha} = {\rm id}$ for all~$\alpha \in (0,\pi)$, thus deducing Theorem~\ref{thm:universalCNSS}. 
That~$M_{\pi/2,\alpha} = {\rm id}$ follows from a subtle interplay between Theorem~\ref{thm:linear}, an additional symmetry of~$\phi_{\bbL(\pi/2)}$, and the way that the  star-triangle transformation acts on isoradial graphs. 
In addition to ultimately being shorter and more self-contained than the direct approach, it illustrates that the star-triangle transformation alone implies that the isoradial embedding is the ``correct'' embedding to ensure universality.

\paragraph{Uniformity in angles.}

In Theorem~\ref{thm:universalCNSS}, the constants $c$ and $C$ may be chosen uniformly in $\alpha$. 
This requires some extra work.
Indeed, Theorem~\ref{thm:linear} is not uniform in the angles $\alpha$, $\beta$, so even if we assume $M_{\frac{\pi}2,\alpha} = {\rm id}$,
Theorem~\ref{thm:linear} only implies a non-uniform version of Theorem~\ref{thm:universalCNSS}. In Section~\ref{sec:unif_angles}, we explain how, if we assume $M_{\frac{\pi}2,\alpha} = {\rm id}$ for all $0<\alpha<\pi$, 
the proof of Theorem~\ref{thm:linear} may be adapted so that the constants are uniform in $\alpha$. 
In Section~\ref{sec:drift0}, we will first show that $M_{\frac{\pi}2,\alpha} = {\rm id}$ for all $0<\alpha<\pi$
by using (the non-uniform version of) Theorem~\ref{thm:linear}, then invoke its uniform version to obtain Theorem~\ref{thm:universalCNSS}. 
We advise the reader to ignore this issue in a first reading. 

 As a consequence, $\alpha$ may be taken to tend to $0$, yielding results for the so-called quantum (or continuum) random-cluster model    (see \cite[Sec.~9.2]{Gri10} for a description). 
    This is a model of cuts placed on the vertical axes of $\bbZ \times \bbR$, and bridges between these axes.
    Its critical version is obtained as a limit of $\phi_{\bbL(\alpha)}$ as $\alpha$ tends to $0$ 
    (see \cite[Sec.~1.3 and 5]{DumLiMan18} for a full description and properties of the critical quantum random-cluster model). 
    
    Write $\phi_{\calQ, 2\delta(\bbZ \times \bbR)}$ for the critical quantum random-cluster model on the rescaled set of vertical axes $2\delta(\bbZ \times \bbR)$. The factor $2$ in the spacing between axes is chosen so that the model is exactly the limit of $\phi_{\delta\bbL(\alpha)}$ as $\alpha$ tends to $0$. 
    
\begin{corollary}\label{cor:quantum}
 	For any $q \in [1,4]$, there exist constants $c,C$ such that, for any $\delta >0$, 
    \begin{align*}
    		d_{\rm CN}(\phi_{\calQ, 2\delta(\bbZ \times \bbR)}, \phi_{\delta\bbL(\tfrac\pi2)}) < C \delta^c.    	\end{align*}
		    Furthermore, we also have that 
    \begin{align*}
	d_{\rm SS}(\phi_{\calQ, 2\delta(\bbZ \times \bbR)}, \phi_{\delta\bbL(\tfrac\pi2)})  \longrightarrow 0.
    \end{align*}
\end{corollary}

The above follows directly from Theorem~\ref{thm:universalCNSS} by taking $\alpha \to 0$. We will not provide further details here. 

  \subsection{Applications to the Potts model}

In this section, we explain certain consequences of Theorem~\ref{thm:rotation_invariance_infinite_vol} for the Potts model. 
The model is defined as follows. Let~$\bbT_q$ be the simplex in~$\bbR^{q-1}$ containing~$(1,0,\dots,0)$ such that for any~$a,b\in\bbT_q$,
$$a\cdot b:=\begin{cases}\ \ 1&\text{ if~$a=b$,}\\
\ -\frac1{q-1}&\text{ otherwise}\end{cases}$$
(above and below~$\cdot$ denotes the scalar product). Attribute a {\em spin} variable~$\sigma_x\in\bbT_q$ to each vertex~$x\in V$. A
 {\em spin configuration}~$\sigma=(\sigma_x:x\in V)\in\bbT_q^{V}$ is given by the collection of all the spins.
  Introduce the Hamiltonian of~$\sigma$ defined by
$$H_G(\sigma):=-\sum_{xy\in E}\,\sigma_x\cdot\sigma_y.$$
The {\em Gibbs measure on~$G$ at inverse temperature~$\beta\ge0$} is defined by the formula, for every~$f:\bbT_q^{V}\rightarrow \bbR$,
\begin{align}\label{eq:Gibbs}\mu_{G,\beta,q}[f]:=\frac1{Z_{\rm Potts}(G,\beta,q)}\sum_{\sigma\in \bbT_q^{V}}f(\sigma)\exp[-\beta H_{G}(\sigma)].\end{align}
Similarly to the random-cluster model, the Potts model exhibits a phase transition at inverse temperature~$\beta_c(q):=\tfrac{q-1}q\log(1+\sqrt q)$, which separates a phase where correlations decay exponentially fast from a phase where they do not decay. 
When~$q\in\{2,3,4\}$, the  phase transition is continuous, as predicted by Baxter (see e.g.~the book \cite{Bax89}) and proved in \cite{DumSidTas16}. 
We will work with fixed~$q \in \{2,3,4\}$ and~$\beta=\beta_c$ and therefore drop them from the subscript in the measure. 

The following corollary, stating the rotational invariance of the spin field, is an immediate application (via the Edwards-Sokal coupling) of the corresponding one for the random-cluster model.

\begin{corollary}[Rotation invariance of spin-spin correlations]\label{cor:Potts_correlations}
Fix~$q\in\{2,3,4\}$ and a simply connected domain~$\Omega$ with a~$C^1$-smooth boundary. For every~$n$ and~$\ep > 0$, there exists~$\delta_0=\delta_0(q,n,\ep, \Omega)>0$ such that for every~$\alpha\in[0,2\pi]$ and~$\delta\le \delta_0$, every~$\tau_1,\dots,\tau_n\in \mathbb T_q$, and every~$x_1,\dots,x_n\in\Omega_\delta$ at a distance at least~$\ep$ from each other and from the boundary of~$\Omega$,
\[
|\mu_{\Omega_\delta}[\sigma_{x_i}=\tau_i,1\le i\le n]-\mu_{(e^{{\rm i}\alpha}\Omega)_\delta}[\sigma_{e^{{\rm i}\alpha} x_i}=\tau_i,1\le i\le n]|\le \ep\, \mu_{\Omega_\delta}[\sigma_{x_i}=\tau_i,1\le i\le n],%\,\mu_{\Omega_\delta}[\sigma_{x_i}=\tau_i,1\le i\le n],
\]
where we use, in a slight abuse of notation,~$e^{{\rm i}\alpha} x_i$ to denote a vertex~$x$  of~$(e^{{\rm i}\alpha}\Omega)_\delta$ within a distance~$\delta$ from the image of~$ x_i$ under the rotation by the angle~$\alpha$.
\end{corollary}

%As before one may check that both terms on the left are of the order of~$f(\delta)^n$. The previous corollary implies that by rescaling the correlations by~$f(\delta)^{-n}$, and then taking convergent sub-sequential limits, one necessarily ends up with rotationally invariant correlations.

\begin{remark} 
	Deducing the rotation invariance of energy~$n$-point correlations 
	(i.e.~the correlations of the random variables~$\epsilon_e^\Omega:=\sigma_x\cdot\sigma_y-	\mu_\Omega[\sigma_x\cdot\sigma_y]$ for~$e=xy$ an edge of~$G$) requires proving the rotation invariance of the energy correlations of the random-cluster model. For this reason, we do not discuss this result here.\end{remark}

These results are known for the Ising model (i.e.~the~$q=2$ Potts model). In fact, in this case the existence of the scaling limit and its conformal invariance were obtained in \cite{CheHonIzy15} for the spin field, and  in \cite{HonSmi13} for the energy field.
  
\subsection*{Structure of the paper}
Section~\ref{sec:2} contains background on isoradial graphs, the associated random-cluster model and the star-triangle transformation. 
It introduces some of the notation used in the following sections, as well as the results of \cite{GriMan14, DumLiMan18}, which are precursors to those presented here. 
Theorem~\ref{thm:linear} is proved in Section~\ref{sec:universality_isoradial_lin}. 
This is the most technical part of the paper; it is self-contained and may be skipped in a first reading. 
Theorems~\ref{thm:rotation_invariance_infinite_vol} and \ref{thm:universalCNSS} are proved in Sections~\ref{sec:deducing_rot_inv} and~\ref{sec:drift0} respectively. %In Section~\ref{sec:BA} we describe the direct approach to proving that  $M_{\frac{\pi}2,\alpha} = {\rm id}$ for all $\alpha$. This is not necessary for the results of the paper and is only presented for illustration. 
Section~\ref{sec:consequences} contains the derivation of the various consequences of Theorem~\ref{thm:rotation_invariance_infinite_vol}, 
namely, Corollaries~\ref{cor:rotation_invariance_domain}, \ref{cor:crossing}, \ref{cor:connectivity_correlations}, and~\ref{cor:Potts_correlations}.

\section{Preliminaries}\label{sec:2}

This section contains a brief introduction to isoradial graphs and the random-cluster model associated to them.

\subsection{Isoradial graphs}\label{sec:isoradial_graphs}

Isoradial graphs were introduced by Duffin in~\cite{Duf68} in the context of discrete complex analysis, 
and later appeared in the physics literature in the work of Baxter~\cite{Bax78} under the name~$Z$-invariant graphs.
They have been studied extensively, in particular in the context of statistical mechanics; we refer to~\cite{CheSmi12,KS-quad-graphs, Mer01, GriMan14,DumLiMan18} for literature on the subject.

A rhombic tiling~$\bbG^\diamond$  is a tiling of the plane by rhombi of edge-length~$1$. Any such graph is bipartite, and we may divide its vertices in two sets of non-adjacent vertices~$\bbV_{\bullet}$ and~$\bbV_{\circ}$. 
The \emph{isoradial graph}~$\bbG$ associated with~$\bbG^\diamond$ is the graph with vertex set~$\bbV_{\bullet}$ and edge-set given by the diagonals of the faces of~$\bbG^\diamond$ between vertices of~$\bbV_{\bullet}$. If the roles of~$\bbV_{\bullet}$ and~$\bbV_{\circ}$ are exchanged, we obtain the dual of~$\bbG$, which is also isoradial. 
The rhombic tiling~$\bbG^\diamond$ is called the {\em diamond graph} of~$\bbG$.

A {\em train-track} (or simply track) of~$\bbG$ is a bi-infinite sequence of adjacent faces~$(r_i)_{i\in \bbZ}$ of the diamond graph~$\bbG^\diamond$, 
with the edges shared by each~$r_i$ and~$r_{i+1}$ being parallel. The angle formed by any such edge with the horizontal axis
is called the {\em transverse angle} of the track.

\begin{figure}[t]
  \centering
  \includegraphics[width=0.70\textwidth]{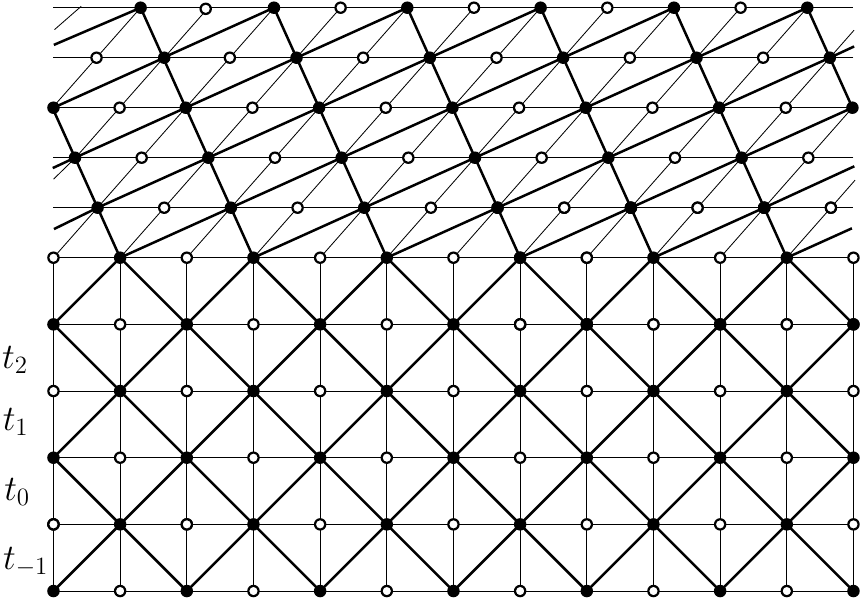}
  \caption{An example of a graph~$\bbL(\pmb\alpha)$, where~$\alpha_i$ is equal to~$\tfrac\pi2$ for~$i\le 3$, and some angle~$\alpha$ above. The diamond graph is drawn in light black lines; the solid and hollow dots are the vertices of~$V_\bullet$ and~$V_\circ$, respectively. The actual isoradial graph~$\bbG$ is drawn in thicker black lines. % Both the lower and upper parts are portions of rectangular lattices; the former is simply a rotated square lattice.
  }
 \label{fig:isoradial lattice unaltered}
\end{figure}

Isoradial graphs considered in this paper are of a very special type, see Figure~\ref{fig:isoradial lattice unaltered}.
They will all be isoradial embeddings of the square lattice in which all rhombi of~$\bbG^\diamond$ have bottom and top edges that are horizontal. 
A consequence of this assumption is that the diamond graph contains {\em horizontal tracks}~$t_i$ with transverse angles~$\alpha_i \in (0,\pi)$ and {\em vertical tracks}~$s_j$, all of which have transverse angle~$0$. Each track of one category intersects all tracks of the other category but no track of the same category. 
Write~$t_{i-1}^+ = t_i^-$ for the set of vertices between~$t_{i-1}$ and~$t_i$.

%This angle, denoted~$\alpha_i$ of the track~$t_i$ is the angle made between the horizontal line and the non-horizontal segments of the diamonds (see the figure).  
For a sequence of track angles~$\pmb\alpha=(\alpha_i)_{i\in \mathbb Z}\in(0,\pi)^\bbZ$, denote by~$\mathbb  L(\pmb\alpha)$ the graph whose horizontal tracks have transverse angles~$\alpha_i$, in increasing vertical order.
When~$\alpha_i=\alpha$ for every~$i$, simply write~$\bbL(\alpha) = \bbL(\pmb\alpha)$.
Note that~$\bbL(\alpha)$ has identical rectangular faces and is rotated such that~$e^{{i}\alpha/2}\mathbb R$ acts as an axis of symmetry.  
In particular,~$\bbL(\tfrac\pi2)$ is simply a rescaled and rotated (by an angle of~$\pi/4$) version of~$\mathbb Z^2$.
These are indeed the isoradial rectangular lattices described in Section~\ref{sec:results_universality}.

When considering isoradial graphs~$\bbG = (V,E)$, we keep the notation~$\Lambda_n = [-n,n]^2$ and identify it with the subgraph spanned by the vertices of~$V$ contained in~$\Lambda_n$. We write~$\partial \Lambda_n$ for the set of vertices~$v \in V \cap \Lambda_n$ that have at least one neighbour in~$V \cap \Lambda_n^c$. 
Let us conclude this section by mentioning that we will (almost) always consider~$\bbG=\bbL(\pmb\alpha)$.

\subsection{Random-cluster model on isoradial graphs}

For a graph~$G=(V,E)$ included in an isoradial graph~$\bbG=(\bbV,\bbE)$ with vertex-set~$V$ and edge-set~$E$, a {\em boundary condition}~$\xi$ on~$G$ is given by a partition of the set~$\partial G$ of vertices in~$V$ incident to at least one vertex in~$\bbV\setminus V$. We say that two vertices of~$G$ are {\em wired together} if they belong to the same element of the partition~$\xi$. Recall that a {\em cluster} is a connected component of~$\omega$. 

\begin{figure}[t]
  \centering
  \includegraphics[width=0.3\textwidth]{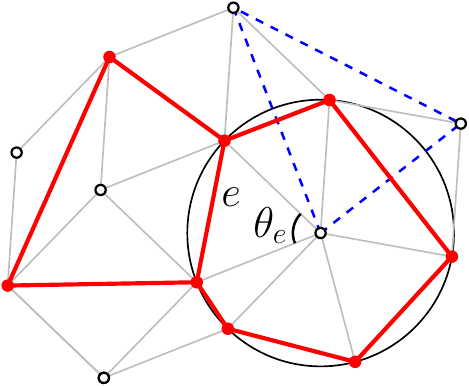}
  \caption{The edge~$e$ and its subtended angle~$\theta_e$; the red edges are those of~$G$ and the grey ones are those of the diamond graph. The dual graph (blue edges, hollow vertices) is also isoradial.}
  \label{fig:subtended_angle}
\end{figure}

In the paper, we will always work with the random-cluster model on an isoradial graph with specific weights, called {\em isoradial weights}, associated with this graph as follows. 
If~$e$ is an edge of~$G$ and~$\theta_e$ is the angle of the rhombus of~$G^\diamond$ containing~$e$ that is not bisected by~$e$ (see Figure~\ref{fig:subtended_angle}), we set
\begin{align}\label{eq:isoraial_p_e}
	p_e:=\begin{cases}\displaystyle \frac{\sqrt{q}\sin( r (\pi - \theta_e))}{\sin( r \theta_e)+\sqrt{q}\sin( r (\pi - \theta_e))} 
    &\text{if } q<4, \\[5mm]
\qquad\ \ \displaystyle  \frac{2\pi - 2\theta_e}{2\pi-\theta_e} &  \text{if } q = 4,\\[5mm]
\displaystyle \frac{\sqrt{q}\sinh( r (\pi - \theta_e))}{\sinh( r \theta_e)+\sqrt{q}\sinh( r (\pi - \theta_e))} 
    &\text{if }q>4,\end{cases}
\end{align}
where~$r := \tfrac{1}{\pi} \cos^{-1} \left( \tfrac{\sqrt{q}}{2} \right)$ for~$q\le 4$ and $r := \tfrac{1}{\pi} \cosh^{-1} \left( \tfrac{\sqrt{q}}{2} \right)$ for~$q>4$. The case~$q>4$ is not relevant for this presentation, but we give the formula to emphasise that the weights exist for all~$q \geq 1$ and change nature at~$q = 4$.
Notice that when~$\theta_e = \pi/2$ we find 
\begin{align}
	p_e:=\tfrac{\sqrt q}{1 + \sqrt q}.
\end{align}

\begin{definition} The random-cluster measure with  isoradial edge-weights and cluster-weight~$q>0$ on a finite graph~$G$  with boundary conditions~$\xi$ is given by
\begin{align}\label{eq:RCM_def2}
	\phi_{G,q}^\xi[\omega]:=\frac{q^{k(\omega^\xi)}}{Z^\xi_{\mathrm{RCM}}(G,q)} \prod_{e\in E}p_e^{\omega_e}(1-p_e)^{1-\omega_e},
\end{align}
where~$k(\omega^\xi)$ is the number of connected components of the graph~$\omega^\xi$  which is obtained from~$\omega$ by identifying wired vertices together, and~$Z^\xi_{\mathrm{RCM}}(G,q)$ is a normalising constant called the {\em partition function} chosen in such a way that~$\phi_{G,q}^\xi$ is a probability measure. 
\end{definition}

Two specific families of boundary conditions will be of special interest to us. On the one hand, the {\em free} boundary conditions, denoted 0, correspond to no wirings between boundary vertices. On the other hand, the {\em wired} boundary conditions, denoted 1, correspond to all boundary vertices being wired together. 

We will also consider the random-cluster model on infinite isoradial graphs~$\bbG$ with free boundary conditions obtained by taking the limit of the measures with free boundary conditions on larger and larger finite graphs~$G$ tending to~$\bbG$. Set~$\phi_{\bbG,q}$ for the measure in infinite volume, which, as shown in \cite{DumLiMan18}, is unique for~$1\le q\le 4$. 

The choice of the isoradial parameters is such that the model is critical. This result was obtained in the case of the square lattice in \cite{BefDum12} and for isoradial graphs in \cite{DumLiMan18} (see also the anterior paper \cite{BefDumSmi12} for the case~$q>4$).

\begin{center}
{\em Henceforth, fix $q \in [1,4]$. For isoradial graphs $\bbG$, $\phi_\bbG$  denotes an infinite-volume measure on $G$ with parameter $q$ and edge-weights given by~\eqref{eq:isoraial_p_e}. } 
\end{center}
For all practical purposes, one may consider the infinite-volume measure to be unique.

\subsection{Elementary properties of the random-cluster model}\label{sec:elementary}

We will use the following standard properties of the random-cluster model. They can be found in \cite{Gri06}, and we only recall them briefly below. 

\bigbreak\noindent
{\em Monotonic properties.} Fix~$G$ as above. 
An event~$A$ is called {\em increasing} if for any~$\omega\le\omega'$ (for the partial ordering on~$\{0,1\}^E$ given by~$\omega\le\omega'$ if~$\omega_e\le \omega'_e$ for every~$e\in E$),~$\omega\in A$  implies that~$\omega'\in A$.
Fix~$q\ge1$ and boundary conditions~$\xi'\ge\xi$, where~$\xi'\ge\xi$ means that any wired vertices in~$\xi$ are also wired in~$\xi'$. 
Then, for any increasing events~$A$ and~$B$,
\begin{align*}\tag{FKG}\label{eq:FKG} 
	\phi_{G}^\xi[A\cap B]&\ge \phi_{G}^\xi[A]\phi_{G}^{\xi}[B],\\
	\tag{CBC}\label{eq:CBC} 
	\phi_{G}^{\xi'}[A]&\ge \phi_{G}^\xi[A].
\end{align*}

The inequalities above will respectively be referred to as the {\em FKG inequality} and the {\em comparison between boundary conditions}. 

\bigbreak\noindent{\em Spatial Markov property.} For any configuration~$\omega'\in\{0,1\}^E$ and any~$F\subset E$,
\begin{align}\label{eq:SMP} \tag{SMP}
	\phi_{G}^\xi[\cdot_{|F}\,|\,\omega_e=\omega'_e,\forall e\notin F]= \phi_{H}^{\xi'}[\cdot],
\end{align}
where~$H$ denotes the graph induced by the edge-set~$F$, and~$\xi'$ are the boundary conditions on~$H$ defined as follows: 
$x$ and~$y$ on~$\partial H$ are wired if they are connected in~$(\omega'_{|E\setminus F})^\xi$. 

\begin{figure}[htb]
  \centering
  \includegraphics[width=0.5\textwidth]{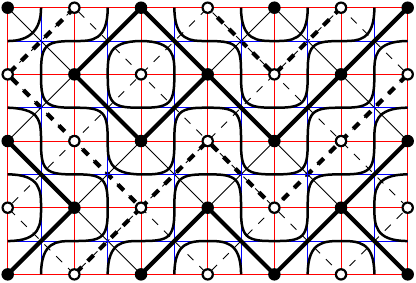}
  \caption{We depict in black, dotted black, red and blue respectively the primal, dual, diamond and medial lattices. The primal configuration~$\omega$ is in bold and the dual one~$\omega^*$ in dashed bold. Finally, the loop configuration~$\overline\omega$ is in black. }
  \label{fig:configurations}
\end{figure}

A direct consequence of the spatial Markov property is the {\em finite-energy property} guaranteeing that, conditioned on the states of all the other edges in a graph, the probability that an edge $e$ is open is between~$p_e/(p_e+q(1-p_e))$ and~$p_e$.
\bigbreak\noindent
{\em Dual model.} Define the dual graph~$G^*=(V^*,E^*)$ of~$G$ as follows: 
place dual sites at the centres of the faces of~$G$ (the external face, when considering a graph in the plane, must be counted as a face of the graph), 
and for every edge~$e\in E$, place a dual edge between the two dual sites corresponding to faces bordering~$e$ (see Figures~\ref{fig:subtended_angle} and~\ref{fig:configurations}). When the graph is isoradial, we make the following choice for the position of dual vertices in~$V^*$: the vertex~$v^*$ corresponding to a face of~$G$ is placed at the center of the corresponding circumcircle. The dual of an isoradial graph is by construction an isoradial graph.

Given a subgraph configuration~$\omega$, construct a configuration~$\omega^*$ on~$G^*$ by
declaring any edge of the dual graph to be open (resp.\ closed) if the
corresponding edge of the primal lattice is closed (resp.\ open) for the
initial configuration. The new configuration is called the \emph{dual
  configuration} of~$\omega$.
The dual model on the dual graph given by the dual configurations then
corresponds to a random-cluster measure 
with isoradial weights and dual boundary conditions. We will not discuss the details of how dual boundary conditions are defined (we refer to \cite{Gri06} for details in the general case and to \cite{DumLiMan18} for the isoradial setting) and we simply observe that the dual of the free boundary conditions are the wired ones, and vice versa. 
\bigbreak
\noindent
{\em Loop model.} The loop representation of a configuration on~$G$ is supported on the {\em medial graph} of~$G$ defined as follows (see Figure~\ref{fig:configurations}). 
Recall that, for~$\bbG$ an isoradial graph,~$\bbG^\diamond$  is the associated diamond graph, with each face corresponding to a pair of mutually dual edges $e$, $e^*$ of $\bbG$ and $\bbG^*$, respectively. 
The medial graph of $\bbG$ is the dual $(\bbG^\diamond)^*$ of its diamond graph. Each of its vertices thus corresponds to a pair of mutually dual edges.

Let~$G^\diamond$ be the subgraph of~$\bbG^\diamond$ spanned by the edges of~$\bbG^\diamond$ adjacent to a face corresponding to a vertex of~$G$.
Consider a configuration~$\omega$ on~$G$ and recall its dual configuration~$\omega^*$. 
Draw self-avoiding paths on the medial graph as follows: a path arriving at a vertex of the medial graph 
always  turns so as not to cross the open edges of~$\omega$ or~$\omega^*$. 
The loop configuration~$\overline\omega$ thus defined is formed of  disjoint loops and paths between boundary vertices.
Together these form a partition of the edges of~$(G^\diamond)^*$.

\subsection{Box crossing (RSW) property and its consequences}\label{sec:RSW}

As in the case of the square lattice, it was proved in~\cite{DumLiMan18} that, for~$q \in [1,4]$, there exists a unique infinite-volume measure on a large variety of isoradial graphs and that it has similar properties to the critical random-cluster model on~$\bbZ^2$. 
Most importantly for us, the RSW estimates of Theorem~\ref{thm:RSW_iso} also hold for the isoradial rectangular lattice which we will use below. 

When discussing isoradial graphs, it is convenient to think of edges as closed segments in the plane~$\bbR^2 = \bbC$.
Then, (open) paths are piecewise linear paths running along (open) edges; they are not required to end at vertices of the graph. 

\begin{theorem}[RSW on isoradial rectangular lattices]\label{thm:RSW_iso}
	For any~$1\le q\le 4$ and~$\rho,\ep>0$, there exists~$c=c(\rho,\ep)>0$ 
	such that for any~$\alpha,\beta\in(0,\pi)$ and any $\pmb\alpha=(\alpha_i:i\in\bbZ) \in \{\alpha,\beta\}^\bbZ$, 
	any~$n\geq1$ and any event~$A$ depending on the edges at a distance at least~$\eps n$ from the rectangle~$R:=[0,\rho n]\times[0,n]$,
    \begin{align}\label{eq:RSW_iso}    	
    	c\le \phi_{ \bbL(\pmb\alpha)}\big[\{0\} \times [0,n] \xlra{\, R\,\,} \{\rho n\} \times [0,n]\,\big|\,A\big]\le 1-c.\tag{RSW}
    \end{align}
\end{theorem}

In the above $\phi_{ \bbL(\pmb\alpha)}$ denotes any infinite-volume random-cluster  measure on $\bbL(\pmb\alpha)$.
The statement is actually about finite volume measure with arbitrary boundary conditions; it is stated using the infinite-volume setting to avoid introducing additional notation.  

It is in fact a direct consequence of~\eqref{eq:RSW_iso} that, for all $1 \leq q \leq 4$ and  an isoradial rectangular lattice $\bbL(\pmb\alpha)$ as in the statement of Theorem~\ref{thm:RSW_iso}, there exists a {\em unique} infinite-volume random-cluster  measure $\phi_{\bbL(\pmb\alpha)}$.
Furthermore, under this measure, neither the primal nor the dual configurations contain infinite clusters. 

\begin{remark}\label{rem:RSW_uniform_angles}
Notice that the bounds in~\eqref{eq:RSW_iso} are uniform in the angles of the lattice $\bbL(\pmb\alpha)$;
in particular, we have uniform RSW estimates for $\phi_{ \bbL(\alpha)}$ as $\alpha$ tends to $0$. 
This fact requires additional arguments beyond RSW estimates which only hold for angles $\alpha$ bounded away from $0$ and $\pi$.
A full proof may be found in~\cite{DumLiMan18}.\end{remark}

\begin{remark}
We limit Theorem~\ref{thm:RSW_iso} to lattices $\bbL(\pmb\alpha)$ with at most two different angles for the horizontal tracks as this statement follows directly from~\cite{DumLiMan18} and suffices for our purposes. An RSW estimate holding uniformly for all rectangular lattices $\bbL(\pmb\alpha)$ with no restriction on the angle sequence may also be derived from~\cite{DumLiMan18} but requires additional work. 
\end{remark}

An important consequence of~\eqref{eq:RSW_iso} we will be using repeatedly  is the mixing property.

\begin{proposition}[Mixing property]\label{prop:mixing}
There exist~$C_{\rm mix},c_{\rm mix}\in(0,\infty)$ such that for any~$\alpha,\beta\in(0,\pi)$ and any $\pmb\alpha=(\alpha_i:i\in\bbZ) \in \{\alpha,\beta\}^\bbZ$, 
every~$r\le R/2$, every event~$A$ depending on edges in~$\Lambda_r$, and every event~$B$ depending on edges outside~$\Lambda_R$, we have that 
\begin{align}
	\big|\phi_{\bbL(\pmb\alpha)}[A\cap B]-\phi_{\bbL(\pmb\alpha)}[A]\phi_{\bbL(\pmb\alpha)}[B]\big|\le C_{\rm  mix}(r/R)^{c_{\rm mix}}\phi_{\bbL(\pmb\alpha)}[A]\phi_{\bbL(\pmb\alpha)}[B].\label{eq:mixing}
\end{align}
\end{proposition}

\begin{proof}
	The argument is identical to the one on the square lattice, see e.g.~\cite[Proposition~2.9]{DumMan20}.
\end{proof}

\subsection{Incipient infinite cluster in the half-plane}\label{sec:IIC}

Fix two angles~$\alpha, \beta \in (0,\pi)$. 
Write~$\bbL_{\rm mix}$ for the lattice~$\bbL(\pmb\alpha)$ with~$\pmb\alpha = (\alpha_i)_{i\in\bbZ}$ with~$\alpha_i = \alpha$ for~$i$ even and~$\alpha_i = \beta$ for~$i$ odd.

\begin{figure}
\begin{center}
\includegraphics[width=0.65\textwidth]{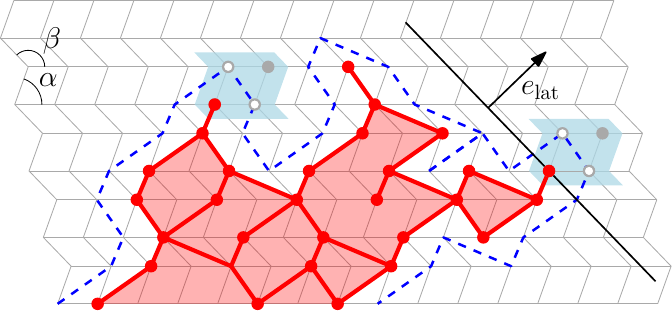}
\caption{A red cluster~$\sfC$ with the cells containing~${\rm Top}(\sfC)$ and~${\rm Right}(\sfC)$ marked in blue. 
Note that there are three topmost cells visited by~$\sfC$, the marked one is the leftmost.}
\label{fig:IIC}
\end{center}
\end{figure}

Recall the following notation: the horizontal tracks of~$\bbL_{\rm mix}$ are denoted~$(t_k)_{k\in \bbZ}$ and the vertical ones~$(s_k)_{k\in \bbZ}$.
Consider the origin of~$\bbR^2$ to be the vertex between~$t_0$ and~$t_1$ and between~$s_0$ and~$s_1$; we will assume it to be a primal vertex.
The cell~$(i,j)$ is the set of primal and dual vertices contained between~$s_{2i-1}$ and~$s_{2i+1}$ and between~$t_{2j-1}$ and~$t_{2j+1}$. 
We associate to the cell~$(i,j)$ its lower left lattice point (due to our definition, this is a primal point); see the blue parts in Figure~\ref{fig:IIC}.

Define the unit vectors~$e_{\rm vert} = (0,1)$ and~$e_{\rm lat} = (\sin\beta,-\cos \beta)$    of~$\bbR^2$. We choose to work here in the basis~$(e_{\rm vert},e_{\rm lat})$ as the two directions are well adapted to the lattice~$\bbL_\beta$. This is ultimately visible in Proposition~\ref{prop:drift_RT}.

For a finite cluster~$\sfC$ of a configuration~$\omega$ on~$\bbL_{\rm mix}$, 
consider the leftmost cell of maximal vertical coordinate intersected by~$\sfC$.
Write~${\rm Top}(\sfC)$ for its associated lattice point. 
Write~${\rm Bottom}(\sfC)$ for the lattice point of the rightmost cell of~$\bbL_{\rm mix}$ of minimal vertical coordinate intersected by~$\sfC$.
Consider the lattice points~${\rm Left}(\sfC)$ and~${\rm Right}(\sfC)$ corresponding to the cells intersected by~$\sfC$ of minimal, respectively maximal, scalar product with~$e_{\rm lat}$. 
When multiple choices are possible, let~${\rm Left}(\sfC)$ be the bottommost such point and~${\rm Right}(\sfC)$ be the topmost. 
Finally, write 
\begin{align}\label{eq:TBLR}
    {\rm T}(\sfC) &= \langle {\rm Top}(\sfC), e_{\rm vert} \rangle;&
    {\rm B}(\sfC) &= \langle {\rm Bottom}(\sfC), e_{\rm vert} \rangle;\nonumber\\
    {\rm L}(\sfC) &= \langle {\rm Left}(\sfC), e_{\rm lat} \rangle;&
    {\rm R}(\sfC) &= \langle {\rm Right}(\sfC), e_{\rm lat} \rangle.
\end{align}
We call these the {\em extremal coordinates} of~$\sfC$ with respect to~$(e_{\rm vert},e_{\rm lat})$.

Write~$x_{n}$ for the primal vertex of~$\bbL_{\rm mix}$ between~$t_{2n}$ and~$t_{2n+1}$ and between~$s_0$ and~$s_1$. 
Similarly, write~$y_n$ for the primal vertex between~$t_{0}$ and~$t_{1}$ and between~$s_{2n}$ and~$s_{2n+1}$. 
Write~$\sfC_{x}$ for the cluster of the vertex~$x$. 
The IIC measure is given by the following proposition. 

\begin{proposition}(IIC measure)\label{prop:IIC_def}
	The following limits exist and are called the half-plane IIC measures in the lower, upper, right and left half-planes, respectively: 
	\begin{align*}
	\phi^{\rm IIC,T}_{\bbL_{\rm mix}}&:= \lim_{n \to - \infty}	\phi_{\bbL_{\rm mix}} [\cdot\,|\, {\rm Top}(\sfC_{x_n}) = 0 ],\\
	\phi^{\rm IIC,B}_{\bbL_{\rm mix}}&:= \lim_{n \to \infty}	\phi_{\bbL_{\rm mix}} [\cdot\,|\, {\rm Bottom}(\sfC_{x_n}) = 0 ],
	\\
		\phi^{\rm IIC,L}_{\bbL_{\rm mix}}&:= \lim_{n \to \infty}	\phi_{\bbL_{\rm mix}} [\cdot\,|\, {\rm Left}(\sfC_{y_n}) = 0 ],
	\\ 
	\phi^{\rm IIC,R}_{\bbL_{\rm mix}}&:= \lim_{n \to -\infty}	\phi_{\bbL_{\rm mix}} [\cdot\,|\, {\rm Right}(\sfC_{y_n}) = 0 ].
	\end{align*}
\end{proposition}

The~\eqref{eq:RSW_iso} property applies to $\bbL_{\rm mix}$ and ultimately allows us to define the half-plane IIC measures. 
%This construction is tedious but relatively standard, so we will only state the results here, and direct the reader to~\cite{Kes86,MendesThesis} for details. 

The measures above refer to the {\em local} environment around the extrema (in the sense of~\eqref{eq:TBLR}) of large clusters. 
Indeed, it may be shown that extrema of a typical large cluster are distributed according to~$\phi^{\rm IIC,*}_{\bbL_{\rm mix}}$, even when the cluster is conditioned on (reasonable) large-scale features. 
This is a manifestation of a wider principle of mixing between the local and large-scale features of the model and is a consequence of the~\eqref{eq:RSW_iso} property.

\begin{lemma}\label{lem:IIC}
	There exist constants~$C,c_{\rm IIC}>0$ such that the following holds for all~$r \leq R/2$. 
	For any configuration~$\omega_0$ on~$\Lambda_R^c$, 
	any union~$\calC$ of clusters of~$\omega_0$ with some~$x \in \calC$,
	and any event~$A$ depending on~$\Lambda_r$, 
	\begin{align*}
        \big|\,	\phi^{\rm IIC,T}_{\bbL_{\rm mix}}[A]-\phi_{\bbL_{\rm mix}}[A\,|\, \omega = \omega_0\text{ on } \Lambda_R^c,\,\sfC_x \cap \Lambda_R^c = \calC,\, {\rm Top}(\sfC_x)=0]\,\big|\,&\le C(r/R)^{c_{\rm IIC}},
    \end{align*}
    as long as the conditioning is non-degenerate, meaning that there exists at least one configuration in $\Lambda_R$ such that $\sfC_x \cap \Lambda_R^c = \calC$ and ${\rm Top}(\sfC_x)=0$. 
    Moreover, the same holds for all other directions~${\rm Bottom}(\cdot)$,~${\rm Left}(\cdot)$ and~${\rm Right}(\cdot)$, with the corresponding IIC measure. 
\end{lemma}

In the above, the constants~$C,c_{\rm IIC}>0$ only depend on the RSW property, and therefore are uniform over the choice of lattice~$\bbL_{\rm mix}$. As mentioned above, Lemma~\ref{lem:IIC} is a consequence of~\eqref{eq:RSW_iso} 
in the spirit of previous IIC constructions~\cite{Kes86,Jar03,Jar03b,damron2011outlets,garban2013pivotal,BasSap17,MendesThesis}.

We will not give a full proof here, but will mention some differences with previously published versions. 
First, most previous versions refer only to Bernoulli percolation. In light of~\eqref{eq:RSW_iso}, the proofs adapt readily; see also \cite{GasManMoh26} for a generic arm separation statement which applies to FK-percolation and is useful when adjusting the arguments. 
Secondly, the IIC here is in the half-plane, in that the three arms imposed by the conditioning are contained in a half-plane. 
Thirdly, the IIC is not rooted at a single point, but rather at one of the two primal vertices in the cell of $0$. 
This poses no additional difficulties, but it should be kept in mind that the choice between the two possible values of the extremum is part of the IIC measure. 

Lastly, and most importantly from a technical point of view, the trace~$\calC$ of~$\sfC_x$ on~$\Lambda_R^c$ is an arbitrary set, which will generically contain a divergent number of components. More traditionally, the conditioning is on~$\omega$ outside of a well-separated flower domain (or some other convenient domain), which essentially boils down to~$\calC$ containing a single component; see \cite{garban2013pivotal, MendesThesis} and the notion of simple flower domain in \cite[Exercise 5.7]{Man25}. 
For cleaner applications, we choose here this stronger form of Lemma~\ref{lem:IIC}. 
To reduce it to its more traditional versions, one may use the following claim. 

Write~${\rm Only}(\calC)$ for the event that no point of ~$\Lambda_R^c \setminus \calC$ is connected to~$\Lambda_{R/2}$. 
Let~${\rm Merge}(\calC) $ be the event that all components of~$\calC$ are connected to each other in~$\Lambda_{R/4}^c$. 

\begin{claim}\label{claim:traceIIC}
	There exists a constant~$c > 0$ such that, for any~$R$ large enough and~$\omega_0$,~any $\calC$ and~$x \in \calC$ as in Lemma~\ref{lem:IIC},
	\begin{align}\label{eq:traceIIC}
	\phi_{\bbL_{\rm mix}}\big[ {\rm Only}(\calC) \cap {\rm Merge}(\calC)  \,\big|\, \omega = \omega_0 \text{ on } \Lambda_R^c,\,\sfC_x \cap \Lambda_R^c = \calC,\, {\rm Top}(\sfC_x)=0\big] \geq c.
	\end{align}
\end{claim}

As this is the only element of the proof of Lemma~\ref{lem:IIC} that is not contained in existing literature, we sketch its proof below. 

\begin{proof}[Proof sketch of Claim~\ref{claim:traceIIC}]
	Start by proving a lower bound on the probability of~${\rm Only}(\calC)$ under the measure in~\eqref{eq:traceIIC}.
	Conditionally on any realisation of the clusters of~$\calC$ in the whole of~$\Lambda_R$, the measure in their complement is dominated by the unconditioned measure 
	$\phi_{\bbL_{\rm mix}}[. \,|\, \omega = \omega_0\text{ on } \Lambda_R^c]$.  
	As such,~\eqref{eq:RSW_iso} applies to the dual model, and one may construct with positive probability a dual ``barrier'' that prevents any point of~$\Lambda_R^c \setminus \calC$  from being connected to~$\Lambda_{R/2}$. 
	Thus
	\begin{align}\label{eq:traceIIC0}
	\phi_{\bbL_{\rm mix}}\big[ {\rm Only}(\calC) \,\big|\, \omega = \omega_0 \text{ on } \Lambda_R^c,\,\sfC_x \cap \Lambda_R^c = \calC,\, {\rm Top}(\sfC_x)=0 \big] \geq c_0, 
	\end{align}
	for some universal~$c_0>0$.

	Next, we prove that 
	\begin{align}\label{eq:traceIIC1}
	\phi_{\bbL_{\rm mix}}\big[ {\rm Merge}(\calC)  \,\big|\, \omega = \omega_0 \text{ on } \Lambda_R^c,\,\sfC_x \cap \Lambda_R^c = \calC,\, {\rm Top}(\sfC_x)=0,\, {\rm Only}(\calC) \big] \geq c_1,
	\end{align}
	for some universal~$c_1 > 0$.
	Explore all primal clusters of~$\Lambda_R^c \setminus \calC$. Due to~${\rm Only}(\calC)$, this exploration reveals no edges of~$\Lambda_{R/2}$. 
	Then, explore the top boundary~$\Gamma$ of~$\sfC_x$ in~$\Lambda_R$, and everything above it. 
	Now, in the unexplored region, the measure is conditioned by an increasing event: that all components of~$\calC$ be connected together. 
	By~\eqref{eq:FKG}, we may apply~\eqref{eq:RSW_iso} for the primal model in this part of the space. 
	Thus, with positive probability there exist open paths in~$\Lambda_{R/2}\setminus \Lambda_{R/4}$ separating~$\Lambda_{R/4}$ from~$\Lambda_{R/2}^c$ in the unexplored region. 
	When such paths exist, due to the conditioning,~$ {\rm Merge}(\calC) $ necessarily occurs. 
	This concludes the proof of~\eqref{eq:traceIIC1}. Combining it with~\eqref{eq:traceIIC0}, we obtain~\eqref{eq:traceIIC}. 
\end{proof}

We close by observing that Proposition~\ref{prop:IIC_def} follows readily from Lemma~\ref{lem:IIC}.

\subsection{Star-triangle and the track-exchange transformations}\label{sec:star-triangle}

In this section, we present the track-exchange transformation. In order to do so, we first introduce the star-triangle transformation and then define the track-exchange transformation as the result of a sequence of star-triangle transformations.

\paragraph{Star-triangle transformation} The \emph{star-triangle transformation}, also known as the \emph{Yang-Baxter relation}, was first discovered by Kennelly in 1899 in the context of electrical networks~\cite{Kennelly}.
Then, it became a key relation in different models of statistical mechanics~\cite{Bax89,Onsager} indicating the integrability of the system.
We do not plan to do a full review on this transformation (see for instance \cite{DumLiMan18} for more details) and focus directly on the context of the random-cluster model on isoradial graphs with isoradial edge-weights.

\begin{figure}[htb]
  \centering
  \includegraphics[width=0.6\textwidth]{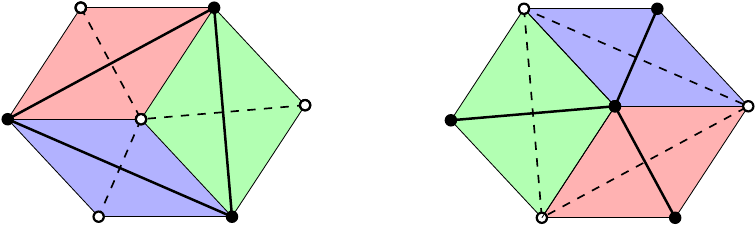}
  \caption{\label{fig:diamond}The three diamonds together with the drawing, on the left, of the triangle (in which case the dual graph in dashed is a star) and, on the right, of the star (in which case the dual graph is a triangle).}
\end{figure}

First of all, note that for any triangle~$ABC$ contained in an isoradial graph, there exists a unique choice of point~$O$ (namely the orthocenter) such that, if the triangle~$ABC$ is replaced by the star~$ABCO$, 
    the resulting graph is also isoradial. Conversely, for every star~$ABCO$ in an isoradial graph, the graph obtained by removing this star and putting the triangle~$ABC$ is isoradial. This process of changing the graph is called the {\em star-triangle transformation}. Note that triangles and stars of isoradial graphs correspond to hexagons formed of three rhombi in the diamond graph. 
Thus, when three such rhombi are encountered in a diamond graph, 
they may be permuted as in Figure~\ref{fig:diamond} using a star-triangle transformation.

%We now introduce the {\em cube flip} (or star-triangle) operator~$F^\cube_k:\Omega_{\bbL^k}\times [0,1]\rightarrow \Omega_{\bbL^{k+1}}$ on configurations. 

The star-triangle transformation was first used to prove that the laws of connections between vertices of a graph~$G$ with a triangle~$ABC$ and the graph~$G'$ obtained from~$G$ with the star~$ABCO$ instead of~$ABC$ are the same, except for the additional vertex~$O$ in~$G'$. In other words,  the star-triangle transformation produces a coupling between the random-cluster models on~$G$ and that on~$G'$ that preserves connections. This coupling interpretation appears in a number of works, see for instance \cite{DumLiMan18}, as well as \cite{GriMan13,GriMan13a,GriMan14} for Bernoulli percolation. 
That the star-triangle transformation preserves the partition functions of models on isoradial graphs with isoradial weights goes further back, for instance~\cite{Bax78, Bax89, BaxEnt78, Gri06, Mer01} among others.

\begin{figure}[htb]
  \centering
  \includegraphics[width=0.6\textwidth]{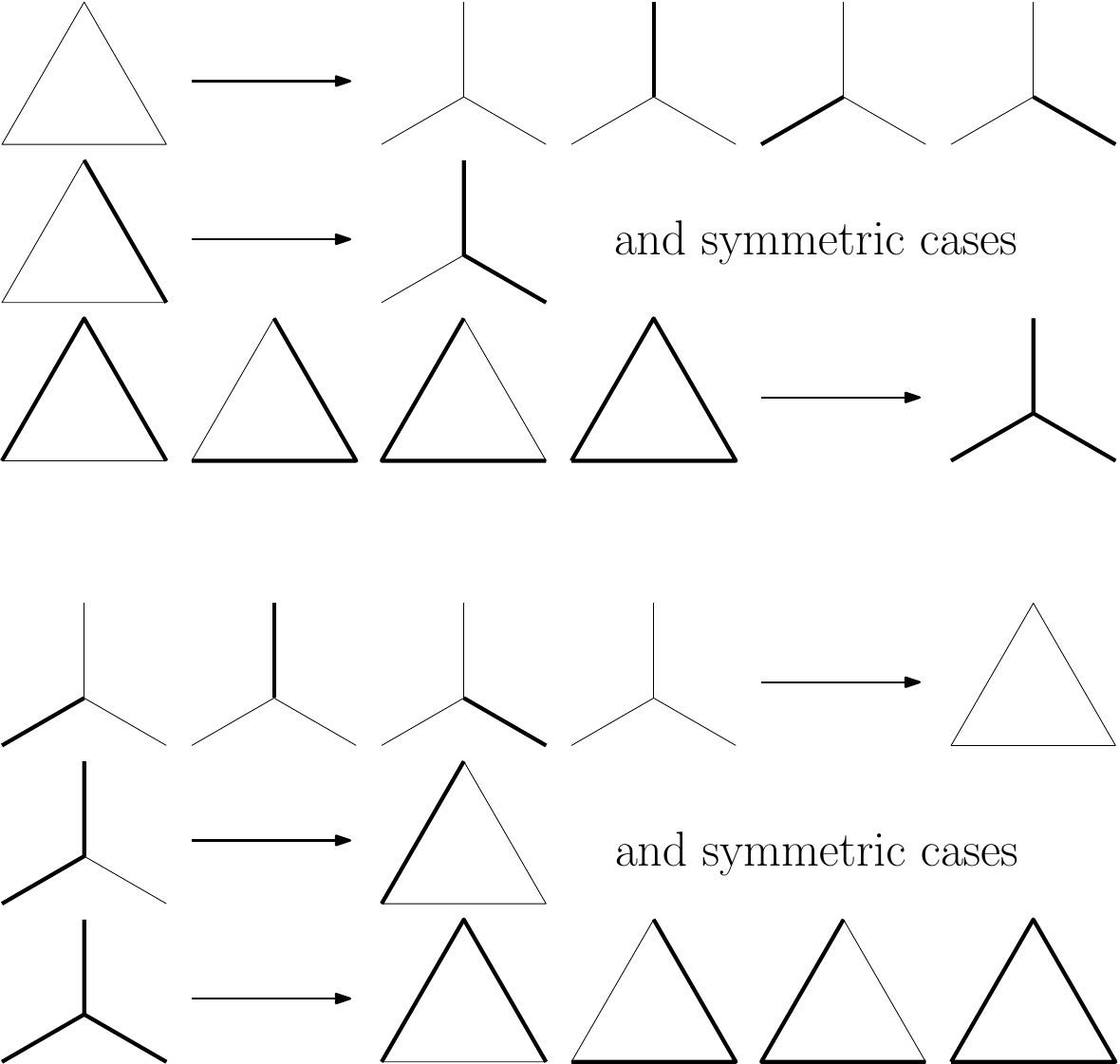}
  \caption{A picture of the possible transformations in the star-triangle coupling (the probabilities in the case of multiple outcomes are described in the definition). We also picture the reverse map.}
  \label{fig:simple_transformation_coupling}
\end{figure}

\begin{definition}[Star-triangle coupling \cite{DumLiMan18}] Consider a graph~$G$ containing a triangle~$ABC$ and let~$G'$ be the graph with the star~$ABCO$ instead. Introduce the {\em star-triangle coupling} between~$\omega\sim\phi_{G}^\xi$ and~$\omega'\sim\phi_{G'}^\xi$ defined as follows (see Figure~\ref{fig:simple_transformation_coupling}):
\begin{itemize}[noitemsep]
\item For every edge~$e$ which does not belong to~$ABCO$,~$\omega'_e=\omega_e$,
\item If two or three of the  edges of~$ABC$ are open in~$\omega$, then all the edges in~$ABCO$ are open in~$\omega'$,
\item If exactly one of the edges of~$ABC$ is open in~$\omega$, say~$BC$, then the edges~$BO$ and~$OC$ are open in~$\omega'$, and the third edge of the star is closed in~$\omega'$,
\item If no edge of~$ABC$ is open in~$\omega$, then~$\omega'_{OABC}$ has 
\begin{itemize}[noitemsep]
\item no  open edge with probability equal to~$\tfrac{1-p_{OA}}{p_{OA}}\tfrac{1-p_{OB}}{p_{OB}}\tfrac{1-p_{OC}}{p_{OC}}$, 
\item the edge~$OA$ is open and the other two closed with probability~$q\tfrac{1-p_{OB}}{p_{OB}}\tfrac{1-p_{OC}}{p_{OC}}$, 
\item similarly with cyclic permutations for~$B$ and~$C$.
\end{itemize}
\end{itemize}
\end{definition}

Let us make a few observations concerning the coupling.
First, note that the transformation uses extra randomness in one case and that it is not a deterministic matching of the different configurations.
Second, the coupling preserves the connectivity between the vertices, except at the vertex~$O$. 
Third, in the coupling, given~$\omega'$, the edges of~$ABC$ in~$\omega$ are sampled as follows:
\begin{itemize}[noitemsep]
\item If there is one or zero edge of~$ABCO$ that is open in~$\omega'$, then none of the edges in~$ABC$ is open in~$\omega$,
\item If exactly two of the edges in~$ABCO$ are open in~$\omega'$, say~$AO$ and~$BO$, then the edge~$AB$ is the only edge of~$ABC$ that is open in~$\omega$,
\item If all the edges of~$ABCO$ are open in~$\omega'$, then 
\begin{itemize}[noitemsep]
\item all the edges of~$ABC$ are open in~$\omega$ with probability~$\tfrac1q\tfrac{p_{AB}}{1-p_{AB}}\tfrac{p_{BC}}{1-p_{BC}}\tfrac{p_{CA}}{1-p_{CA}}$, 
\item~$AB$ and~$BC$  are open and~$CA$ is closed with probability equal to~$\tfrac1q\tfrac{p_{AB}}{1-p_{AB}}\tfrac{p_{BC}}{1-p_{BC}}$, 
\item similarly with cyclic permutations.
\end{itemize}\end{itemize}

\begin{figure}[htb]
  \centering
  \includegraphics[width=0.70\textwidth]{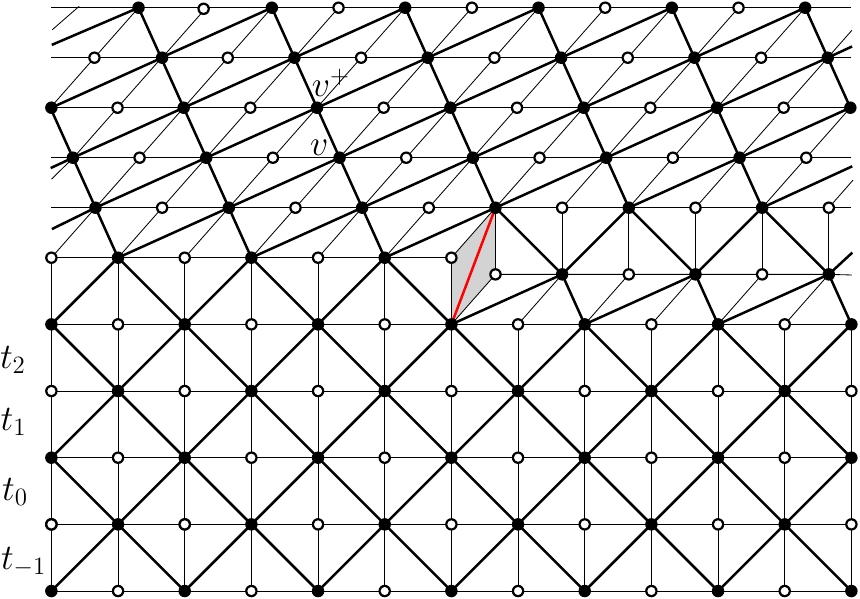}
  \caption{An example of a graph~$\bbL_k$ for~$k=4$. What happens between tracks~$t_2$ and~$t_5$ is a mixture of the isoradial lattice~$\bbL$ with angles~$\pi/2$  for~$i\le 3$ and~$\alpha$ for~$i\ge 4$, and~$\bbL'$ is obtained by exchanging the tracks~$4$ and~$5$. The only diamond that does not belong to~$\bbL$ or~$\bbL'$ is in grey. }
 \label{fig:isoradial lattice}
\end{figure}

\paragraph{Track-exchange operator} Recall the definition of $t_i$ from Section~\ref{sec:isoradial_graphs}. The previous star-triangle operator gives rise to a track-exchange operator defined as follows. For~$\bbL=\bbL(\pmb\alpha)$ and~$i\in\bbZ$, let~$\bbL'=\bbL(\pmb\alpha')$ be the lattice obtained by exchanging the tracks~$t_i$ and~$t_{i-1}$ (i.e.~by exchanging~$\alpha_i$ and~$\alpha_{i-1}$ in the sequence~$\pmb\alpha$). Index the vertices of~$t_{i-1}^-$ from left to right by~$(x_k:k\in\bbZ)$ and assume that~$\alpha_{i-1}>\alpha_i$. Also, let~$\bbL_k$ be the isoradial graph, see Figure~\ref{fig:isoradial lattice}, obtained by 
\begin{itemize}[noitemsep]
\item taking the same diamonds as~$\bbL$ (or equivalently~$\bbL'$) on~$t_j$ with~$j\notin\{i-1,i\}$;
\item taking the same diamonds as~$\bbL$ on the part of~$t_{i-1}$ and~$t_i$ on the left of~$x_k$;
\item taking the same diamonds as~$\bbL'$ on the part of~$t_{i-1}$ and~$t_i$ on the right of~$x_{k}$;
\item adding a diamond above~$x_k$ to complete the gap.
\end{itemize}
Note that the properties above determine all the diamonds in~$\bbL_k$, and that there is only one diamond in~$\bbL_k$ which does not belong to either~$\bbL$ or~$\bbL'$. Denote this diamond by~$\lozenge$.
We now  define an operator sending configurations on~$\bbL$ to configurations on~$\bbL'$, that gives a formal meaning to the intuitive idea of inserting~$\lozenge$ at the position~$+\infty$ and using the star-triangle transformation to exchange the tracks by moving~$\lozenge$ step by step to~$-\infty$.

Let~$\omega$ be some configuration on~$\bbL$ and define for every~$k\in\bbZ$ the configuration~$\tilde{\omega}_k$ on~$\bbL_{k}$ coinciding with~$\omega$ on the diamonds common to~$\bbL_k$ and~$\bbL$ (i.e.~outside~$t_{i-1},t_{i}$ and on the left of~$x_k$), and defined arbitrarily otherwise. Denote~$\tilde{\omega}_k^k := \tilde{\omega}_k$ and for every~$j< k$, define inductively~$\tilde{\omega}_k^j$ to be the result of the star-triangle transformation mapping a configuration on~$\bbL_{j+1}$ to a configuration on~$\bbL_j$, applied to~$\tilde{\omega}_k^{j+1}$. Define~$\omega_k := \lim_{j\to -\infty} \tilde{\omega}_k^j$, which is a configuration on~$\bbL'$. 

\begin{figure}
\begin{center}
\includegraphics[width = .5\textwidth]{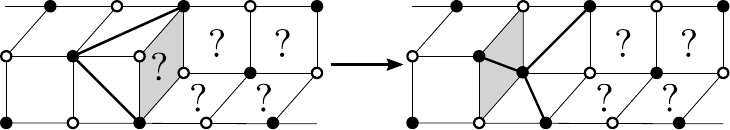}
\caption{An example where the configuration on the two diamonds left of the grey
diamond~$\lozenge$ determines the configuration on~$\lozenge$ after the
star-triangle transformation, regardless of the configuration on or
right of~$\lozenge$.
}
\label{fig:track_exchange}
\end{center}
\end{figure}

Now, observe the important fact that if we have three integers~$k,k' \ge j$ such that 
$\tilde{\omega}_{k}^{j}$ and~$ \tilde{\omega}_{k'}^{j}$ coincide on~$\lozenge$, then the (local) outcome of the star-triangle transformation from~$\tilde{\omega}_{k}^{j}$ and~$\tilde{\omega}_{k'}^{j}$ will be the same (as long as it uses the same external randomness). More generally, applying all the subsequent steps we see that~$\omega_{k}$ and~$\omega_{k'}$ coincide on the part of~$t_{i-1}\cup t_i$ that is to the left of~$x_j$. Finally, notice that some configurations on the two diamonds left of~$\lozenge$ in~$\bbL_k$ fix deterministically the state of~$\lozenge$ in~$\bbL_{k-1}$ after a single star-triangle transformation (see e.g.~Figure~\ref{fig:track_exchange}). Denote by~$F_k$ this event. If~$F_k$ occurs for~$\omega$, then for all~$k', k'' > k$, it also does by definition for~$\omega_{k'}$ and~$ \omega_{k''}$, and therefore~$\omega_{k'}$ and~$\omega_{k''}$ coincide left of~$x_k$. This leads to the following definition.

%Let the {\em cube flip} (or star-triangle) operator~$F^\cube_k:\Omega_{\bbL^k}\times [0,1]\rightarrow \Omega_{\bbL^{k+1}}$ be defined in Figure \ref{fig:def-cube-flip}, see also Section~\ref{sec:star-triangle}.

\begin{definition}[Track exchange by star-triangle transformation]\label{def:track exchange}
If~$\alpha_{i-1}>\alpha_i$, and~$\omega$ is a percolation configuration on~$\bbL$ such that~$\omega\in F_k$ occurs for an infinite number of indices~$k>0$, define the track-exchange operator~$\mathbf T_i$ by~$\mathbf T_i(\omega)=\lim_{k\rightarrow \infty}\omega_k$, where~$\omega_k$ is defined as in the previous paragraph.
\end{definition}

We will only work with measures (random cluster measures, IIC measures) that verify some finite energy property so that~$F_k$ occurs for an infinite number of~$ k>0$ almost surely. Hence the operator~$\mathbf T_i$ is well defined on almost all configurations~$\omega$.

Furthermore, the configuration $\mathbf T_i(\omega)$ for edges adjacent to $x_j$ is determined by $\omega \cap \Lambda_r(x_j)$ (and the external randomness used in the transformations affecting this part of the space), as soon as some $F_k$ occurs for $j < k \leq j + r/2$. By the finite energy property, the probability of this event approaches $1$ exponentially in $r$. We therefore say that $\mathbf T_i$ is {\em local}. 

If~$\alpha_i>\alpha_{i-1}$, we construct~$\mathbf T_i$ similarly by inverting the left and the right, and~$-\infty$ and~$+\infty$. 
When~$\alpha_i = \alpha_{i-1}$, define~$\mathbf T_i$ as the identity map. 

It should be noted that the mixing property of the random-cluster model implies that the random-cluster measure on~$\bbL$ is the limit of the random-cluster measures on~$\bbL_k$. Therefore, 
\begin{align}\label{eq:T_law}
	\text{if~$\omega$ is distributed according to~$\phi_\bbL$, then~$\mathbf T_i(\omega)$ has law~$\phi_{\bbL'}$.}
\end{align}
Let us also insist on the fact that~$\mathbf T_i$ is not a deterministic map, as at each step where a star-triangle operator is used, there is extra randomness in the outcome of the transformation. 

We finish this section with an important proposition.

\begin{proposition}\label{prop:same_law}
	If~$\pmb\alpha$ and~$\pmb\beta$ satisfy~$\alpha_i=\beta_i$ for~$a\le i\le b$, 
	the law of~$\omega$ restricted to the strip between~$t_a^-$ and~$t_b^+$ and
	of the connections between the points of $t_b^+$ in the half-space above $t_b^+$
	and between the points of $t_a^-$ in the half-space below $t_a^-$
	%the homotopy 	classes of loops in~$\omega$ with respect to points in this strip 
	is the same in~$\phi_{\bbL(\pmb\alpha)}$ and~$\phi_{\bbL(\pmb\beta)}$. 
\end{proposition}

%\im{Removed homotopy from here as it was not introduced and not really needed for this statement

\begin{proof}
As a sequence of star-triangle transformations, the track-exchange operator preserves the connection properties of the vertices that are not on the tracks which are exchanged. From this, one may deduce that for every event~$A$ involving only edges inside the strip and the connections mentioned above 
\[
\phi_{\bbL(\pmb\alpha)}[A]=\lim_{R\rightarrow\infty}\phi_{\bbL(\pmb\alpha(R))}[A]=\lim_{R\rightarrow\infty}\phi_{\bbL(\pmb\beta(R))}[A]=\phi_{\bbL(\pmb\beta)}[A],
\]
where 
\[
\pmb\alpha(R):=\begin{cases}\alpha_i&\text{ if }|i|\le R,\\
\beta_{i-R+b}&\text{ if }i>R,\\
\beta_{R-i+a}&\text{ if }i<-R,\end{cases}\quad\text{and}\quad\pmb\beta(R):=\begin{cases}\beta_i&\text{ if }|i|\le R,\\
\alpha_{i-R+b}&\text{ if }R<i<2R-b,\\
\alpha_{R-i+a}&\text{ if }-2R+a<i<-R,\\
\beta_i&\text{otherwise}.\end{cases}
\] (In the first  and last equalities, we use the measurability and the uniqueness of the infinite-volume measure, and in the second equality, the track-exchange operator.)
\end{proof}

\subsection{Arm exponents}	
Arm exponents generally describe the speed of decay of certain connection probabilities in the critical random-cluster model. 
The RSW property allows one to obtain bounds for such exponents (but does not generally prove their existence), and compute some of them exactly. We state here the bounds necessary for our arguments. 

\newcommand{\onearmbound}{\xi_1}

\begin{proposition}\label{prop:arms_bounds_iso}
	There exist constants~$C,c,\onearmbound > 0$ (with $\xi_1 \leq 1$) such that, 
	for any~$\alpha,\beta\in(0,\pi)$, any $\pmb\alpha=(\alpha_i:i\in\bbZ) \in \{\alpha,\beta\}^\bbZ$,
	and any~$r < R$,
	\begin{align}
	\phi_{\bbL(\pmb\alpha)} [\Lambda_r \longleftrightarrow \Lambda_R^c] &\leq C (r/R)^{\onearmbound}, \label{eq:one_arm_bound_iso}\\
	\tfrac1C (r/R)^{2}\leq \phi_{\bbL(\pmb\alpha)} [\exists \sfC \text{ with }{\rm Top}(\sfC) \in \Lambda_r \text{ and }\sfC \cap \Lambda_R^c \neq \emptyset] &\leq C (r/R)^{2}. \label{eq:3hp_arm}
	\end{align}
	The second line is also valid for~${\rm Bottom}$ instead of~${\rm Top}$.
	Moreover,
	\begin{align}
	\tfrac1C (r/R)^{2}\leq \phi_{\bbL_{\rm mix}} [\exists \sfC \text{ with }{\rm Left}(\sfC) \in \Lambda_r \text{ and }\sfC \cap \Lambda_R^c \neq \emptyset] &\leq C (r/R)^{2}, \label{eq:3hp_arm2} \\
	\phi_{\bbL(\pmb\alpha)} [\exists \sfC \text{ with }{\rm Left}(\sfC) \in \Lambda_r \text{ and }\sfC \cap \Lambda_R^c \neq \emptyset] &\leq C (r/R)^{1+c}.\label{eq:3hp_arm3}
	\end{align}
	 The same holds with~${\rm Right}$ instead of~${\rm Left}$. 
	% \im{Formally, we should replace here $c$ by $\xi_1$, or change the definition \eqref{eq:Rchoice} of $R$. We can also just leave it as is.}
\end{proposition}

It is remarkable that the exponent of the half-plane three-arm probability in~\eqref{eq:3hp_arm} and~\eqref{eq:3hp_arm2} may be computed exactly; 
we say that the three-arm half-plane exponent is {\em universal}, in that it does not depend on~$q \in [1,4]$. 
Notice that~\eqref{eq:3hp_arm2} only holds for~$\bbL_{\rm mix}$, while for generic lattices we only have the weaker bound~\eqref{eq:3hp_arm3} at this stage. 
This is because the proofs of~\eqref{eq:3hp_arm} and~\eqref{eq:3hp_arm2} use the invariance of the lattice under translations orthogonal to the direction of the extremum. 
All lattices~$\bbL(\pmb\alpha)$ are invariant under horizontal translations, which allows us to deduce the exponent for~${\rm Top}$ and~${\rm Bottom}$, but not for~${\rm Left}$ or~${\rm Right}$. The lattices~$\bbL_{\rm mix}$ are invariant under two non-collinear translations, which allow one to extend~\eqref{eq:3hp_arm} to the lateral extrema. 

Theorem~\ref{thm:universalCNSS} ultimately implies that the exponent for~\eqref{eq:3hp_arm3} is also equal to two, but it is not possible to prove this at this stage. We mention the sub-optimal bound~\eqref{eq:3hp_arm3}  here, as it is necessary for the proof of Theorem~\ref{thm:linear}.

The events above are called the one-arm and the half-plane three-arm events, respectively, from a distance~$r$ to a distance~$R$. 
The direction of the half-plane depends on which extremum of the cluster we consider. The exponents in~\eqref{eq:3hp_arm} and~\eqref{eq:3hp_arm2} are the equivalent, on $\mathbb L(\pi/2)$, of the so-called three-arm event in a half space, previously derived in \cite{KesSidZha98,Nol08,BefDum13,CheDumHon16,DumManTas20}. 	

The bounds of Proposition~\ref{prop:arms_bounds_iso} are relatively standard consequences of \eqref{eq:RSW_iso}, so we will only sketch their proofs below. 
The one-arm bound \eqref{eq:one_arm_bound_iso} is obtained exactly as for the square lattice, and we do not give additional details.

\begin{figure}
\begin{center}
\includegraphics[height = 3.5cm]{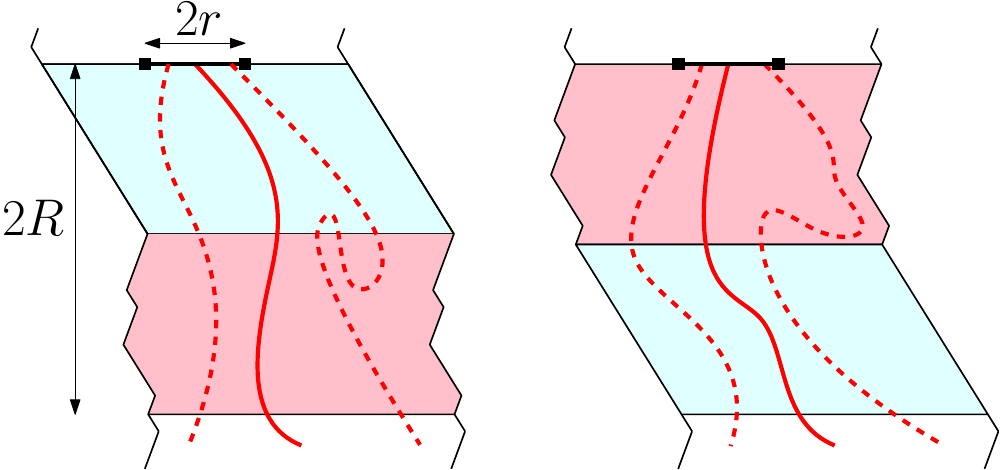}\hspace{2cm}
\includegraphics[height = 3.5cm]{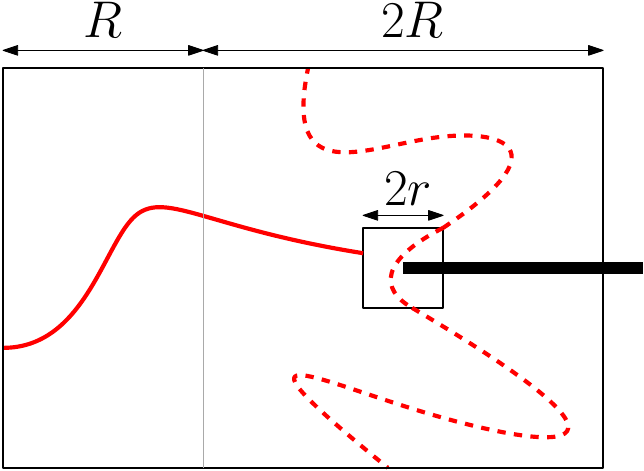}
\caption{{\em Left:} The lattices $\bbL$ and $\bbL'$ with the event $\overline{\rm Arm}_{\rm T}(r,2R+2)$ occurring in both. 
The blue block is formed of tracks of angle $\alpha_0$, the pink one has angles $\alpha_{-1},\dots, \alpha_{-M}$. 
One may exchange these blocks using a series of track-exchanges which do not affect $\overline{\rm Arm}_{\rm T}(r,2R+2)$. 
{\em Right:} The event ${\rm Arm}'_{\rm slit}(r,R)$: the two dual paths and the primal path between them are not allowed to intersect the slit on the right. Observe that at most one of  $R/r$ translates of ${\rm Arm}'_{\rm slit}(r,R)$ may occur in any one configuration. }
\label{fig:3arms}
\end{center}
\end{figure}

\begin{proof}[Proof sketch of~\eqref{eq:3hp_arm} and~\eqref{eq:3hp_arm2}]
	The derivation of~\eqref{eq:3hp_arm} and~\eqref{eq:3hp_arm2} closely follows the previous proofs for the square lattice, requiring only minor modifications to account for the absence of certain symmetries.

	We start by considering the case of~$\bbL_{\rm mix}$ and proving the upper bounds (these are the ones we use in this paper).  
	Write~$\bbT_N$ for the torus formed of~$N\times N$ cells of~$\bbL_{\rm mix}$. 
	Let~$\phi_{\bbT_N}$ be the random-cluster  measure on the torus with edge-weights given by~\eqref{eq:isoraial_p_e}.
	Fix constants~$R = \eps N$ and~$r < R$, where~$\eps > 0$ is some fixed quantity. All constants below depend on~$\eps$, 
	but do not depend on~$r, R$ or~$N$.

We call {\em almost-extremum} a point~$x \in \sfC$ with~$\sfC$ being a cluster of diameter at least~$R$ and
		\begin{align}\label{eq:almost_extremum}
		\langle x, e_{*}\rangle  \geq \max\{ \langle y, e_{*}\rangle : y \in \sfC\} - 2 
		\qquad  \text{ for some~$e_* \in \{\pm e_{\rm vert}, \pm e_{\rm lat}\}$}.
		\end{align}
		 The number of almost-extrema has an exponential tail, uniformly in~$N$. Indeed, using~\eqref{eq:RSW_iso}, one may show (the proof is classical and omitted) that, uniformly in~$N$, the number of clusters in~$\bbT_N$ with diameter at least~$R$ has an exponential tail. Now, assuming that a configuration~$\omega$ has ``many'' almost-extrema compared to large clusters, 
		one may produce configurations~$\omega'$ by modifying~$\omega$ 
		so that each large cluster of~$\omega'$ has at most~$K$ almost-extrema (where~$K$ is some fixed number, say~$K = 100$). 
		This can be done by only modifying $\omega$ in the vicinity of  a single almost-extremum for each large cluster. 
		In particular the probabilities of $\omega$ and $\omega'$ have uniformly bounded ratio. 
		
		Then, on the one hand, the many almost-extrema of $\omega$ provide many ways to produce such configurations~$\omega'$ from~$\omega$. 
		On the other hand, since the modifications are local near the almost-extrema of $\omega'$, of which there are at most $K$ per cluster, 
		there are relatively few ways to reconstruct~$\omega$ from~$\omega'$. As a consequence, the probability to have many almost-extrema but few large clusters is small (such arguments are sometimes referred to as applications of a ``multi-valued map'' principle). 
		
		The argument above implies that the expected number of almost-extrema is uniformly bounded in $N$. Now that this is done, the periodicity of the torus shows that the probability that~$\Lambda_r$ contains an almost-extremum is of order $(r/R)^2$.
	A simple use of mixing enables one to transfer this estimate to the full-plane measure~$\phi_{\bbL_{\rm mix}}$.		
	This proves the upper bounds in~\eqref{eq:3hp_arm} and~\eqref{eq:3hp_arm2} for lattices~$\bbL_{\rm mix}$.
	
	For the lower bound, it suffices to observe that, by \eqref{eq:RSW_iso}, the torus $\bbT_N$ contains at least one cluster of diameter $R$ with positive probability. 
	This produces at least one extremum in $\bbT_N$ for each direction. This translates to the lower bounds in~\eqref{eq:3hp_arm} and~\eqref{eq:3hp_arm2} by the same periodicity and mixing arguments as above. 
	\medskip 
	
	Let us now deduce~\eqref{eq:3hp_arm} for a general lattice~$\bbL(\pmb\alpha)$, where~$\pmb\alpha = (\alpha_j)_{j\in \bbZ}\in\{\alpha,\beta\}^\mathbb Z$. Note that the invariance under (even) vertical translations used in the previous argument is not available anymore. We therefore rely on track-exchanges. Again, we focus first on the upper bound, as this is most useful. 
	
	For some fixed~$R$, write $N  = \lceil \frac{R}{\sin\alpha_0}\rceil $. 
	Also, let $M$ be the lowest number such that $t_{-M}$ is below $\bbR \times \{-R\}$ in $\bbL(\pmb\alpha)$. 
	Define the lattices $\bbL$ and $\bbL'$ of the form $\bbL(.)$ as follows. 
	
	The horizontal tracks $t_j$ of $\bbL$ have angles $\alpha_j$ from the sequence $\pmb\alpha$ for $j\geq 0$. 
	For $j <0$, the track $t_{-N + j}$ has angle $\alpha_j$. All other tracks, that is $t_{-1}, \dots, t_{-N}$, have transverse angle $\alpha_0$. 
	Note that, on the horizontal strip~$\mathbb R\times[-R,0]$, $\bbL$ is identical to~$\bbL(\alpha_0)$.
	
	Similarly $\bbL'$ is the lattice in which $t_j$ has the angle $\alpha_j$ for $j \geq -M$, angle $\alpha_{j +N}$ for $j <- M- N$ and 
	angle $\alpha_0$ for the remaining tracks. It is identical to $\bbL(\pmb\alpha)$ on~$\mathbb R\times[-R,0]$.
	See Figure~\ref{fig:3arms}, left diagram, for an illustration.

	Write \begin{align}
			{\rm Arm}_{\rm T}(r,R) &= \{\exists \sfC \text{ with }{\rm Top}(\sfC) \in \Lambda_r \text{ and }\sfC \cap \Lambda_R^c \neq \emptyset\}
%		\qquad \text{ and }\\
%		\overline{\rm Arm}_{\rm T}(r,2R) &= \{\exists \sfC \text{ with }{\rm Top}(\sfC) \in \Lambda_r \text{ and }\sfC \cap \bbR \times (-\infty,-2R] \neq \emptyset\}.
	\end{align}
	and $\overline{\rm Arm}_{\rm T}(r,R)$ for the event that there exist three paths 
	$\gamma_1,\gamma_2, \gamma_3$ with~$\gamma_1,\gamma_3 \in \omega^*$ and~$\gamma_2 \in \omega$,
	all contained in the lower half plane, connecting $[-r,r] \times \{0\}$ to $\bbR \times  \{-R\}$
	and ordered from left to right. 
	Note that $\overline{\rm Arm}_{\rm T}(r,R)$ only depends on the configuration in the strip~$\mathbb R\times[-R,0]$.
		
	Using~\eqref{eq:RSW_iso} and the arm-separation technique (see for instance \cite{GasManMoh26}),
	one may show that there exists a universal constant~$c >0$ such that
		\begin{align*}
			\phi_{\bbL(\alpha_0)} [{\rm Arm}_{\rm T}(r,R)] 
			\geq c\phi_{\bbL(\alpha_0)} [\overline{\rm Arm}_{\rm T}(r,R)] 
			= c\phi_{\bbL} [\overline{\rm Arm}_{\rm T}(r,R)] 
			\geq c \phi_{\bbL} [\overline{\rm Arm}_{\rm T}(r,2R+2)].
		\end{align*} 	
	The equality in the above is due to Proposition~\ref{prop:same_law}. Similarly
		\begin{align*}
		\phi_{\bbL'} [\overline{\rm Arm}_{\rm T}(r,2R+2)] 
		\geq c\phi_{\bbL'} [\overline{\rm Arm}_{\rm T}(r,R)] 
		= c\phi_{\bbL(\pmb\alpha)} [\overline{\rm Arm}_{\rm T}(r,R)] 
		\geq c^2\phi_{\bbL(\pmb\alpha)} [{\rm Arm}_{\rm T}(r,R)].
		\end{align*} 	

	Finally, there exists a sequence of track exchanges  acting between $t_{-1},\dots, t_{-M-N}$ that transforms $\bbL$ into $\bbL'$. 
	Since $\overline{\rm Arm}_{\rm T}(r,2R+2)$ is expressed in terms of connections between points above $t_{-1}$ and those below $t_{-M-N}$, 
	and since these track exchanges do not alter the existence of such connections, we find 
	\begin{align}
		 \phi_{\bbL} \big[\overline{\rm Arm}_{\rm T}(r,2R+2)\big] = 	\phi_{\bbL'} \big[\overline{\rm Arm}_{\rm T}(r,2R+2)\big].
	\end{align}	
		
	Combining the three previous displayed equations, we find
		\begin{align}
			c^3 \phi_{\bbL(\pmb\alpha)} [{\rm Arm}_{\rm T}(r,R)] \leq  \phi_{\bbL(\alpha_0)} [{\rm Arm}_{\rm T}(r,R)].
		\end{align} 
	Applying the upper bound~\eqref{eq:3hp_arm} to the quantity on the right (which fits in the context of the lattices $\bbL_{\rm mix}$ treated above), 
	we deduce a similar upper bound for  $\phi_{\bbL(\pmb\alpha)}$. The same strategy applies for the lower bound. 
\end{proof}

\begin{proof}[Proof sketch of~\eqref{eq:3hp_arm3}]
	Fix a lattice~$\bbL = \bbL(\pmb\alpha)$ for some sequence~$\pmb\alpha = (\alpha_j)_{j\in \bbZ}\in\{\alpha,\beta\}^\mathbb Z$. 
	
	For~$r < R$, let~${\rm Arm}_{\rm slit}(r,R)$ be the event that there exist three connections~$\gamma_1$,~$\gamma_2$ and~$\gamma_3$ in~$\bbL \setminus (\bbR_+\times \{0\})$
	between~$\Lambda_r$ and~$\partial \Lambda_R$, in counter-clockwise order, with~$\gamma_1,\gamma_3 \in \omega^*$ and~$\gamma_2 \in \omega$. Let~${\rm Arm}'_{\rm slit}(r,R)$ be the event that, in $\bbL \setminus (\bbR_+\times \{0\})$ there exist
		a dual connection between the top and bottom of $\Lambda_R$ 
		and a primal connection between $\Lambda_r$ and $\{-2R\}\times [-R,R]$ contained in $[-2R,R]\times [-R,R]$.
		See Figure~\ref{fig:3arms}, right diagram, for an illustration. 
		
		Using~\eqref{eq:RSW_iso} and the arm-separation technique (see for instance \cite{GasManMoh26})
		we conclude the existence of some constant~$c > 0$ independent of~$r$ and~$R$ such that 
		\begin{align}
			\phi_{\bbL}[{\rm Arm}'_{\rm slit}(r,R)] \geq c \, \phi_{\bbL}[{\rm Arm}_{\rm slit}(r,R)].
		\end{align}
	Since at most one translate of the event~${\rm Arm}'_{\rm slit}(r,R)$ by~$(kr,0)$ for~$0 \leq k < R/r$ occurs, 
	the invariance under horizontal translation of $\phi_{\bbL}$ implies that 
		\begin{align}
			\phi_{\bbL}[{\rm Arm}'_{\rm slit}(r,R)] \leq r/R.
		\end{align}
		
		Using again~\eqref{eq:RSW_iso} and arm-separation, one can prove that there exist universal constants~$c,C > 0$ such that 
		\begin{align}
		\phi_{\bbL} \big[\exists \sfC \text{ with }{\rm Right}(\sfC) \in \Lambda_r \text{ and }\sfC \cap \Lambda_R^c \neq \emptyset\big] 
		\leq C (r/R)^{c}\phi_{\bbL}[{\rm Arm}_{\rm slit}(r,R)].
		\end{align}
		The two previous displayed equations imply~\eqref{eq:3hp_arm3}.
\end{proof}

% Indeed, imagine that Theorem~\ref{thm:linear} holds, but that~$M_{\beta,\alpha}\neq {\rm id}$. Then the probabilities in~\eqref{eq:3hp_arm} under~$\phi_{\bbL_{t}}$ would have the same exponent as a three-arm in a cone of~$\phi_{\bbL_{0}}$ with opening different from~$\pi$. 

\section{Universality up to linear transformation: proof of Theorem~\ref{thm:linear}}\label{sec:universality_isoradial_lin}

%
%The goal of this section is to prove Theorem~\ref{thm:linear}. 
%To start, let us sketch an approach to proving its more powerful version (Theorem~\ref{thm:universalCNSS}), 
%so as to see exactly why  Theorem~\ref{thm:linear} is indeed more accessible than  Theorem~\ref{thm:universalCNSS}.

The goal of this section is to prove Theorem~\ref{thm:linear}. 
The idea % behind proving Theorem~\ref{thm:linear} 
is to progressively transform~$\bbL(\beta)$ into~$\bbL(\alpha)$ using track exchanges, as in Figure~\ref{fig:transformations}. 
Throughout the transformation, we follow the evolution of the configurations, specifically of their large clusters. 
As track exchanges do not break primal nor dual connections, large clusters survive the transformations, but their boundaries may be progressively altered. 
We would like to argue that these alterations act as i.i.d.\ modifications of the cluster boundaries with a fixed law, whose expected value leads to the linear transformation $M_{\beta,\alpha}$.

There are two difficulties in the above reasoning: 
encoding the cluster boundary in a tractable way and proving that the modifications are i.i.d.\ (up to small errors).
The encoding of clusters will be done through their extremal coordinates in the basis $(e_{\rm vert}, e_{\rm lat} ) $ defined in~\eqref{eq:TBLR};
that following these coordinates suffices to determine the more general shape of clusters is a subtle fact and will be explained in Section~\ref{sec:homotopy}.

That the increments of large clusters' extrema are distributed (almost) according to a fixed law is due to the IIC construction of Section~\ref{sec:IIC}. The issue of independence throughout the process will be solved by resampling the configuration in the vicinity of cluster extrema before every transformation. 
 
The expectations of the increments of $\Top(\sfC)$ and $\Right(\sfC)$ for a large cluster $\sfC$ when applying one step of the transformation are called the {\em vertical} and {\em lateral} drifts (note that the drifts of $\Top(\sfC)$ and $\Right(\sfC)$ are the same as those of $\Bottom(\sfC)$ and $\Left(\sfC)$, respectively) and are denoted by ${\rm Drift}_{\rm vert}(\beta,\alpha)$ and $ {\rm Drift}_{\rm lat}(\beta,\alpha)$. A proper definition is given in~\eqref{eq:drift}.
The linear transformation~$M_{\beta,\alpha}$, is then written in the canonical basis of~$\bbR^2$ as
\begin{align}\label{eq:Mbetaalpha}
M_{\beta,\alpha} 
= \begin{pmatrix}
1 &  \frac{{\rm Drift}_{\rm lat} +  {\rm Drift}_{\rm vert} \cos \beta }{\sin \beta}\cdot \frac{\sin\alpha + \sin\beta}{(\sin \alpha -{\rm Drift}_{\rm vert})\sin\beta}\\[8pt]
0 & 1  + {\rm Drift}_{\rm vert} \cdot \frac{\sin\alpha + \sin\beta}{(\sin \alpha -{\rm Drift}_{\rm vert})\sin\beta}
\end{pmatrix}.
\end{align} 
The apparently complicated form of $M_{\beta,\alpha}$ is due to the particular way in which we encode extrema. See Remark~\ref{rem:M_origin} for the relation between the drift and~$M_{\beta,\alpha}$.

\begin{remark}\label{eq:drift_in_any_direction}
Drifts may be defined as below in any direction, and our proof  would work for any choice of $e_{\rm lat}$ not colinear to~$e_{\rm vert}$.
Furthermore, the drifts in different systems of coordinates are related via linear maps, and will ultimately produce the same matrix $M_{\beta,\alpha}$. 
In particular, if they are proved to be null in one basis (i.e. if $M_{\beta,\alpha} = {\rm id}$), then the drift in any direction is also null.

Our choice of basis is designed to highlight the symmetry in Proposition~\ref{prop:drift_RT}.
\end{remark}

We mention the transformation~$M_{\beta,\alpha}$ here so as to emphasise that it does not depend on the different quantities involved in the proof of Theorem~\ref{thm:linear}. 
Indeed, the drifts are determined by the incipient infinite cluster measures (see Section~\ref{sec:IIC}). 
In particular, we will prove that it is continuous in~$\alpha$ and~$\beta$, which will be useful later on. 

\begin{lemma}\label{lem:M_cont}
	The map~$(\beta,\alpha) \mapsto M_{\beta,\alpha}$ is continuous and 
	$M_{\beta,\alpha}$ is invertible for all~$\alpha,\beta \in (0,\pi)$.
\end{lemma}

The proof of Lemma~\ref{lem:M_cont} is deferred to Section~\ref{sec:lemM_cont}, when all the necessary notions will have been introduced. 
The rest of the section is dedicated to proving Theorem~\ref{thm:linear}. The structure of the proof is described at the end of Section~\ref{sec:5sketch}, after having introduced the relevant notation.

\subsection{Setting of the proof of Theorem~\ref{thm:linear}}\label{sec:5sketch}

We will now fix some notation necessary for the proof of Theorem~\ref{thm:linear} and sketch some of its elements. 
Fix two angles~$\alpha$ and~$\beta$ and~$N \geq 1$. We will be interested in the ``large'' loops contained in the {\em observation window}~$[-N,N]\times [0,N]$ (this is a rectangle in Euclidean geometry, containing a number of tracks that depends on the lattice) for configurations~$\omega$ and~$\omega'$ sampled according to~$\phi_{\bbL(\beta)}$ and~$\phi_{\bbL(\alpha)}$, respectively. 
More specifically, our goal is to construct a coupling between~$\omega$ and~$\omega'$ 
so that the loops of~$\omega$ of a diameter larger than~$\eps N$ contained in~$[-N,N]\times [0,N]$ 
may be paired up with loops of~$\omega'$ which are close for the distance~\eqref{eq:loop_dist} and vice-versa. Here, $\eps > 0$ denotes some arbitrarily small quantity, ultimately chosen as a small negative power of~$N$. 

The exact choice of~$\eps$, as well as the proper definition of diameter of a cluster will be given below. For this section, however, we limit ourselves to a rough description. Hereafter, call clusters of diameter at least~$\eps N$ {\em macroscopic}. 

\paragraph{Lattices and track exchanges.}
Write $K := \lceil N/\sin \beta \rceil$ and assume for simplicity that $K$ is even. Set~$\bbL_0 = \bbL(\pmb\alpha)$ for the sequence~$\pmb\alpha$ given by~$\alpha_i = \beta$ if~$i \leq K$ and~$\alpha_i = \alpha$ otherwise. 

\begin{figure}
\begin{center}
\includegraphics[width = 0.33\textwidth, page = 1]{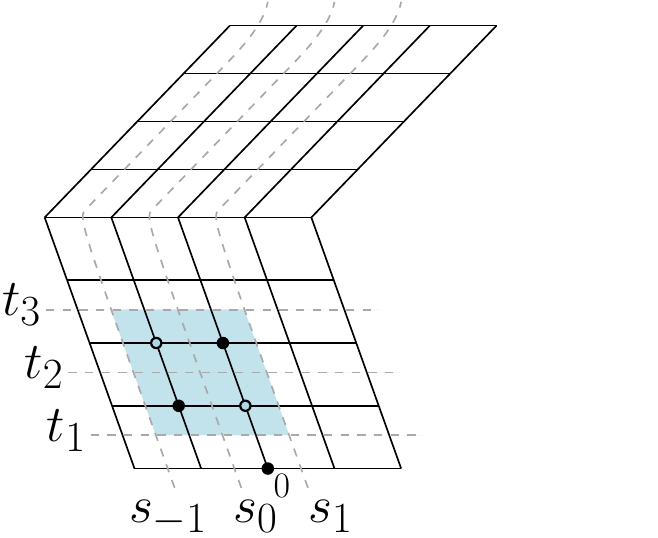}\hspace{-0.12\textwidth}
\includegraphics[width = 0.33\textwidth, page = 3]{figures/transformations.pdf}\hspace{-0.14\textwidth}
\includegraphics[width = 0.33\textwidth, page = 6]{figures/transformations.pdf}\hspace{-0.12\textwidth}
\includegraphics[width = 0.33\textwidth, page = 8]{figures/transformations.pdf}
\caption{From left to right: The initial lattice~$\bbL_0$, with tracks of angle~$\beta$ at the bottom and~$\alpha$ at the top; the blue region is the cell~$(0,1)$ of the lattice. 
After a series of transformations~$\bfS_0,\dots, S_i$ a mixed block appears (red).
After more transformations, a block of angle~$\alpha$ starts appearing at the bottom (blue), which by time~$K+K'$ covers the whole window~$[-N,N]\times [0,N]$.}
\label{fig:transformations}
\end{center}
\end{figure}

Define the lattices~$(\bbL_t)_{t\geq 0}$ inductively as follows. For~$t \geq 0$  set 
\begin{align}\label{eq:5S_transformation}
\bbL_{t+1} = \begin{cases}
	\bfT_1 \circ \bfT_3 \circ \bfT_5 \circ \dots (\bbL_t) \qquad &\text{ if~$t$ odd}\\
	\bfT_2 \circ \bfT_4 \circ \bfT_6 \circ \dots (\bbL_t) \qquad &\text{ if~$t$ even}.
\end{cases}
\end{align}
Write~$\bfS_t$ for the transformation applied when passing from~$\bbL_t$ to~$\bbL_{t+1}$. See Figure~\ref{fig:transformations} for an illustration. 

It may appear as though we apply an infinite number of track exchanges when passing from~$\bbL_{t}$ to~$\bbL_{t+1}$, but this is not actually the case. 
Indeed, the transformations~$\bfT_i$ are non-trivial only when the tracks~$t_{i-1}$ and~$t_{i}$ have different transverse angles, which only happens for a finite number of indices. 
Also note that, due to the parity of $K$,  $\bfS_0$ is trivial and $\bbL_1 = \bbL_0$. 

Generally, each lattice~$\bbL_t$ for~$t \leq K$ is formed of a{\em~$\beta$-block} of transverse angle~$\beta$; 
above which is a {\em mixed block} of alternating tracks of angles~$\alpha$ and~$\beta$, 
above which is an infinite block of tracks of angle~$\alpha$ called the{\em~$\alpha$-block}. 

For~$t > K$, the~$\beta$-block disappears and a block of tracks of angle~$\alpha$ starts forming above~$t_0$. Above this block, there is a mixed block with alternating tracks of angles~$\alpha$ and~$\beta$ (in total there are~$2K$ such tracks), above which we find again a block of tracks of angle~$\alpha$. 
A quick computation shows that the lower boundary of the mixed block is at height~$(K -t)\sin\beta$ for~$t \leq K$ and at height~$(t-K)\sin\alpha$ for~$t > K$. 

Set~$K'$ to be the smallest even integer larger than~$N/\sin \alpha$. Then, in~$\bbL_{K+K'}$, the window~$[-N,N]\times [0,N]$ is fully contained in the bottom~$\alpha$-block of the lattice. 

\paragraph{The configuration chain.}
To the lattices~$(\bbL_t)_{t\geq 0}$, we will associate a Markov chain of configurations~$(\omega_t)_{t\geq0}$ with the property that 
\begin{align*}
	\omega_t \sim \phi_{\bbL_t}\qquad  \text{ for all~$t \geq 0$}. 
\end{align*}
As such, we obtain a {\em coupling} of the laws~$(\phi_{\bbL_t})_{t\geq 0}$. In this coupling, the loops of diameter at least~$\eps N$ of~$\omega_0$ may be paired up with loops of~$\omega_{K+K'}$ that are close in the sense of~\eqref{eq:loop_dist}. 
Proposition~\ref{prop:same_law} states that~$\omega_0 \cap ([-N,N]\times [0,N])$ has the law of a configuration sampled according to~$\phi_{\bbL(\beta)}$, while~$\omega_{K+K'}\cap ([-N,N]\times [0,N])$ has the law~$\phi_{\bbL(\alpha)}$.

The actual Markov chain~$(\omega_t)_{t\geq0}$ is a little complicated to define (see Section~\ref{sec:5MarkovChain}), but the following construction offers a good understanding of it. 
Sample~$\omega_0$ according to~$\phi_{\bbL_0}$. 
Since the track exchanges appearing in~$\bfS_t$ act on disjoint tracks, they may be performed simultaneously to~$\omega_t$
and we write 
\begin{align}\label{eq:omega_t_simple}
	\omega_{t+1} = \bfS_t(\omega_t) \qquad \text{ for all~$t\geq 0$}.
\end{align}

\paragraph{Dynamics of macroscopic clusters.}
Notice that for~$t$ even,~$\bfS_t$ does not affect vertices of~$t_j^+$ with~$j$ even. 
For~$t$ odd, it does not affect the vertices in~$t_j^+$ with~$j$ odd. 
In particular, the connections between any such vertices are preserved when applying~$\bfS_t$. 
As a consequence, any primal or dual cluster of diameter at least~$2$ of~$\omega_t$ has a naturally associated cluster in~$\omega_{t+1}$. 

Clusters may still disappear throughout the process~$(\omega_t)_{t\geq 0}$, but they do so by progressively shrinking to a single point, then disappearing when that point is erased from the lattice. Conversely, new small clusters may appear and potentially grow in size. We will argue below that new clusters do not grow to size~$\eps N$ by time~$K+K'$, nor do initial macroscopic clusters disappear by this time (except on events of vanishing probability). 
This will allow us to follow the evolution of all macroscopic clusters of~$\omega_0$, and pair them with macroscopic clusters of the final configuration~$\omega_{K+K'}$. 

In addition, we would like to argue that the contours of macroscopic clusters are affected throughout the process by i.i.d.\ increments, with a potential drift. Making sense of this statement is challenging as the contour of a cluster is a complicated object. To simplify things, we will only follow the extremal coordinates (see~\eqref{eq:TBLR}) of {\em mesoscopic clusters}. By mesoscopic clusters, we mean clusters of diameter at least~$\eta N$ for some~$0 <\eta \ll \eps$. 
We will eventually prove that controlling the evolution of the extremal coordinates of mesoscopic clusters suffices to control the shape of macroscopic clusters in the sense of~\eqref{eq:loop_dist} (see Section~\ref{sec:homotopy}). 
% We do this via the concept of {\em homotopy information} of a configuration~$\omega_t$, which we define in Section~\ref{sec:5MarkovChain}.

Finally, we would like to argue that each transformation affects the extremal coordinates of mesoscopic clusters in (almost) i.i.d.\ ways. 
As will be apparent below, the increments of the extremal coordinates of the mesoscopic clusters of~$\omega_t$ when applying~$\bfS_t$ depend on the local environments in~$\omega_t$ around these points. 
These local environments will be shown to have a local law independent of the large-scale features of~$\omega_t$, in the spirit of the IIC construction of Section~\ref{sec:IIC}. This will guarantee that the increments are identical in law (up to errors vanishing as~$N$ tends to $\infty$) and that they have expectations ${\rm Drift}_{\rm vert}$ and ${\rm Drift}_{\rm lat}$ given by the IIC measures. 

That the increments are i.i.d.\ in time is more challenging. While we expect some approximate form of independence to hold in the chain~$(\omega_t)_{t\geq0}$ defined above, we are unable to prove it. To circumvent this problem, we will perform a {\em resampling} when passing from~$\omega_t$ to~$\omega_{t+1}$, which will ensure that the large-scale information of~$\omega_t$ is preserved, but the local environments around the extrema of mesoscopic clusters are resampled.

\paragraph{Structure of the proof.}
In Section~\ref{sec:5MarkovChain} we define the notion of mesoscopic cluster and the actual configuration chain~$(\omega_t)_{t\geq0}$. This includes an explicit procedure for the resampling of the configuration around the extrema of mesoscopic clusters.
 
The goal of Sections~\ref{sec:control_error} and~\ref{sec:stability_meso} is to prove the stability of mesoscopic clusters. 
Their conclusion is Proposition~\ref{prop:stability_meso}, which states that, at first order, 
the extremal coordinates of mesoscopic clusters move during the process~$(\omega_t)_{t\geq0}$ as dictated by the drift. 

Section~\ref{sec:homotopy} explains how the knowledge of the extremal coordinates of mesoscopic clusters determines the shape of macroscopic clusters.
The key here is to consider the homotopy classes of contours of macroscopic clusters in the punctured plane obtained by removing the mesoscopic clusters. It will be a consequence of our definition of~$(\omega_t)_{t\geq0}$ and of the stability of mesoscopic clusters that these homotopy classes are preserved throughout the process. 

Finally, in Section~\ref{sec:5linear_proof}, we put together the elements of the previous sections in order to prove Theorem~\ref{thm:linear}. 

Sections~\ref{sec:5equality_drifts} and~\ref{sec:lemM_cont} contain certain additional properties of the matrix~$M_{\beta,\alpha}$, namely the proof of Lemma~\ref{lem:M_cont} 
and a relation between the vertical and lateral drift, which will be crucial when proving Theorem~\ref{thm:universalCNSS}. 

Section~\ref{sec:unif_angles} explains how to obtain the uniformity in the angle $\alpha$ in Theorem~\ref{thm:universalCNSS}.

\subsection{Actual coupling: resampling of extrema}\label{sec:5MarkovChain}
 
In this section, we define the actual Markov chain~$(\omega_t)_{t\geq 0}$. Write~$\bbP$ for the probability measure governing this chain. 
Recall that we assumed that~$K$ and~$K'$ are even. 

Fix some~$0<\eta < 1$; the choice of~$\eta$ will be explained in Section~\ref{sec:5linear_proof}. 
We will use the term ``universal constant'' to mean a positive constant independent of~$N$ and~$\eta$. 

A cluster~$\sfC$ in some configuration~$\omega$ on some lattice~$\bbL_t$ is called {\em$\eta$-mesoscopic} (or simply mesoscopic) 
if 
\begin{align}\label{eq:meso}
	\eta  N \leq {\rm T}(\sfC)-{\rm B}(\sfC) &\leq \eta^{1/2} N,\nonumber\\
	{\rm T}(\sfC) &\leq  N,\nonumber\\
	-N\le {\rm L}(\sfC) \le {\rm R}(\sfC) &\leq  N.
\end{align}
We call~${\rm T}(\sfC)-{\rm B}(\sfC)$ the {\em vertical diameter} of~$\sfC$. 
%A cluster~$\sfC$ on in some configuration~$\omega$ on some lattice~$\bbL_t$ is called~$\eta$-mesoscopic (or simply mesoscopic) 
%if~$\|{\rm Top}(\sfC)-{\rm Bottom}(\sfC)\| \geq \eta N$
%and if both~${\rm Top}(\sfC)$ and~${\rm Bottom}(\sfC)$ are contained in the window~$[-N,N] \times [0,N]$. 
Write~${\rm Meso}(\omega)$ for the set of mesoscopic clusters of~$\omega$ and 
\begin{align}\label{eq:calH_def}
	\calH(\omega) = \big({\rm T}(\sfC), {\rm B}(\sfC), {\rm L}(\sfC), {\rm R}(\sfC)\big)_{\sfC \in {\rm Meso}(\omega)}
\end{align}
for the extremal coordinates of all mesoscopic clusters. 

Sample~$\omega_0 \sim \phi_{\bbL_0}$. For~$t \geq0$ even, assuming~$\omega_0,\dots, \omega_t$ already defined, 
sample a configuration~$\omega_{t+1/2}$ on~$\bbL_{t}$ as follows. 
The goal here is to have~$\omega_{t+1/2} \sim \phi_{\bbL_t}$.

Recall that~$\bbL_{\rm mix}$ denotes the lattice~$\bbL(\pmb\alpha)$ with~$\pmb\alpha = (\alpha_i)_{i\in\bbZ}$, where~$\alpha_i = \alpha$ for~$i$ even and~$\alpha_i = \beta$ for~$i$ odd.
First, sample i.i.d.\ configurations~$(\zeta_j^{{\rm Top}})_{j\in \bbN}$,~$(\zeta_j^{{\rm Bottom}})_{j\in \bbN}$, 
$(\zeta_j^{{\rm Left}})_{j\in \bbN}$ and~$(\zeta_j^{{\rm Right}})_{j\in \bbN}$
 according to the  lower, upper, right and left half-plane IIC measures~$\phi^{\rm IIC,T}_{\bbL_{\rm mix}}$, $\phi^{\rm IIC,B}_{\bbL_{\rm mix}}$, $\phi^{\rm IIC,L}_{\bbL_{\rm mix}}$ and~$\phi^{\rm IIC,R}_{\bbL_{\rm mix}}$, respectively (see Proposition~\ref{prop:IIC_def}). 

%Recall the mesoscopic cluster coordinates~$({\rm Top}(\sfC), {\rm Bottom}(\sfC))_{\sfC \in {\rm Meso}(\omega_t)}$ of~$\omega_t$. 
Fix~
\begin{align}\label{eq:Rchoice}
R := \eta^C N,\qquad \text{ with $C > 0$ a constant such that $(C -1)\onearmbound -2 \geq 1$},
\end{align}
% where~$C > 0$ is a sufficiently large universal constant. We will assume that~$(C -1)\onearmbound -2 \geq 1$ (this choice will become clear later), 
where~$\onearmbound > 0$ is given by~\eqref{eq:one_arm_bound_iso}; the reason for this choice will become clear later. Define the {\em extremum boxes} of~$\omega_t$ to be the boxes~
$$\Lambda_{ R } (x)\text{ for }
x \in \{{\rm Top}(\sfC),\,{\rm Bottom}(\sfC),\,{\rm Left}(\sfC),\,{\rm Right}(\sfC)\,:\, \sfC \in {\rm Meso}(\omega_t)\}.$$ 
We say that an {\em overlap error} occurs at step~$t$ if any of the following occurs: 
\begin{itemize}
	 \item two extremum boxes intersect;
	 \item two distinct clusters of vertical diameter at least~$\eta N/2 - 2R$ intersect the same extremum box; 
	 \item there exists an extremum box that intersects two distinct blocks of the lattice (i.e.,~$\alpha$-block,~$\beta$-block or mixed block).
\end{itemize}
If an overlap error occurs, set~$\omega_{t+1/2} = \omega_t$.

If no overlap error occurs, set~$\omega_{t+1/2} = \omega_t$ everywhere except in the extremum boxes contained in the mixed block. 
We will now describe how to sample~$\omega_{t+1/2}$ in each of these boxes. 
We will do this in such a way that the mesoscopic clusters of~$\omega_t$ and~$\omega_{t+1/2}$ are identical outside of the extremum boxes
and have the same extrema in~$\omega_t$ and~$\omega_{t+1/2}$. 

Enumerate the boxes in an arbitrary way and sample them sequentially. 
Suppose that the configuration~$\omega_{t+1/2}$ in the boxes numbered~$1,\dots, i-1$ has been sampled and let us describe how to sample it in the~$i^{\rm th}$ box. 
Assume that this box is centred at a point~$x = {\rm Top}(\sfC)$ for some~$\sfC \in {\rm Meso}(\omega_t)$; 
the same procedure applies if~$x$ is ${\rm Bottom}(\sfC)$,~${\rm Left}(\sfC)$ or~${\rm Right}(\sfC)$.

To start, let us describe the law~$\calL$ of the configuration in~$\Lambda_{ R } (x)$ under~$\phi_{\bbL_t}$, 
conditionally on the previously sampled edges,
on the fact that the mesoscopic clusters of~$\omega$ and~$\omega_{t}$ are identical outside of the extremum boxes,
and on their extrema being the same in~$\omega$ and~$\omega_{t}$.
It is~$\phi_{\Lambda_R(x)}^\xi$, where~$\xi$ are the  (random\footnote{the boundary conditions are random, as other unsampled boxes may influence the boundary conditions induced by the outside configuration on~$\Lambda_{ R } (x)$.}) boundary conditions given by the configuration outside, 
conditioned on the fact that all pieces of~$\sfC$ in~$\omega_t\setminus \Lambda_{ R } (x)$ are connected in~$\omega$ inside~$\Lambda_{ R } (x)$, 
but not to any other points of~$\Lambda_{ R } (x)^c$,
and that~${\rm Top}(\sfC) = x$.
Indeed, this conditioning ensures that~$\sfC \cap \Lambda_R(x)^c$ remains the restriction to~$\Lambda_R(x)^c$ of a mesoscopic cluster of~$\omega$ 
with same extremal coordinates as in~$\omega_t$. 
Moreover, the resampling in~$\Lambda_{ R } (x)$ may not produce additional mesoscopic clusters, 
since~$\sfC$ is the only cluster of vertical diameter at least~$\eta N/2 - 2R$ intersecting~$\Lambda_{ R } (x)$.

We sample~$\omega_{t+1/2}$ in~$\Lambda_{ R } (x)$ according to~$\calL$, as follows. 
Let~$r := \eta^C R$ for some large fixed constant~$C > 0$. Formally, we need~$C \geq 3/c_{\rm IIC}$, where~$c_{\rm IIC} >0$ is the constant given by Lemma~\ref{lem:IIC}.
Write~$d_{\rm TV}$ for the total variation distance between the restriction of~$\calL$ to~$\Lambda_{ r }(x)$ and the law of~$\zeta_i^{{\rm Top}}$ inside~$\Lambda_{ r }$, translated by~$x$.
%\begin{itemize}
%	\item[(a)] the law of the configuration~$\omega_{t+1/2}$ inside~$\Lambda_{ r }(x)$ 
%	conditionally on the already sampled part, 
%	on the fact that the mesoscopic clusters of~$\omega_t$ and~$\omega_{t+1/2}$ are equal outside of the extremum boxes 
%	and on~$\calH(\omega_{t+1/2}) = \calH(\omega_t)$ and 
%	\item[(b)] the law of~$\zeta_i^{{\rm Top}}$ inside~$\Lambda_{ r }$, translated by~$x$.
%\end{itemize}

Sample a Bernoulli variable~${\rm err}_i$ of parameter~$d_{\rm TV}$
and sample~$\omega_{t+1/2}$ inside~$\Lambda_{ R }(x)$ according to~$\calL$
so that,  
\begin{align}\label{eq:iic_err_couple}
{\rm err}_i = 0 \quad\Longrightarrow \quad \omega_{t+1/2} = \zeta_i^{{\rm Top}} + x \text{ on~$\Lambda_r(x)$},
\end{align}
where~$\zeta_i^{{\rm Top}} + x$ is the translate of~$\zeta_i^{{\rm Top}}$ by~$x$.
This is indeed possible due to the definition of~$d_{\rm TV}$. 
We say that a {\em coupling error} occurs at step~$t$ if there exists~$i$ with~${\rm err}_i = 1$.

This concludes the construction of~$\omega_{t+1/2}$. Notice that~$\omega_{t+1/2}$ has the same law as~$\omega_t$, namely~$\phi_{\bbL_t}$. 
Set
\begin{align*}
	\omega_{t+1} = \bfS_t(\omega_{t+1/2}) \qquad \text{ and  } \qquad \omega_{t+2} = \bfS_{t+1}(\omega_{t+1}).
\end{align*}

We have thus constructed the Markov chain $(\omega_t)_{t\geq 0}$ and have ensured that $\omega_t \sim \phi_{\bbL_t}$ for all $t$. 
Keep in mind the fact that the resampling is done only at even time steps. 

\subsection{Controlling the errors}\label{sec:control_error}

For~$t \geq 0$ even and~$\sfC$ a mesoscopic cluster of~$\omega_t$, 
write~$\tilde\sfC$ for the corresponding cluster in~$\omega_{t+2}$; note that~$\tilde\sfC$ is a well-defined object, but may be non-mesoscopic if it ends up being just below the size required to qualify as mesoscopic. 
For~${\rm Ext} \in \{ {\rm Top}, {\rm Bottom}, {\rm Left}, {\rm Right}\}$, we say~${\rm Ext}(\sfC)$ is  {\em frozen} if it is at a distance at least two from the mixed block, otherwise we call it {\em active}.
That~${\rm Top}(\sfC)$ and~${\rm Bottom}(\sfC)$ are frozen is determined by~${\rm T}(\sfC)$ and~${\rm B}(\sfC)$, respectively, and hence may be read off~$\calH(\omega_t)$. The same is not true for ${\rm Left}(\sfC)$ and ${\rm Right}(\sfC)$. 

Define the increments of the extrema of mesoscopic clusters as
\begin{align*}
	\Delta_t {A}(\sfC)& =  {A}(\tilde \sfC) - {A}(\sfC)\qquad \text{ for~$A \in \{{\rm T,B,L,R}\}$}.
\end{align*}
If~${\rm Top}(\sfC)$ is frozen, then~$\Delta_t {\rm T}(\sfC) = 0$. The same holds for~${\rm Bottom}(\sfC)$. 
For the left and right extrema, this is not the case as, even when~${\rm Left}(\sfC)$ or~${\rm Right}(\sfC)$ are frozen,
a point which is almost extremal in the direction~$e_{\rm lat}$ and which is not frozen may move and become~${\rm Left}(\tilde\sfC)$ or~${\rm Right}(\tilde\sfC)$, respectively. 
However, the increment of all extrema is zero if both~${\rm Top}(\sfC)$ and~${\rm Bottom}(\sfC)$ are frozen and on the same side of the mixed block.

%
%If~${\rm Top}(\sfC)$ is not {\em frozen}, write 
%\begin{align}
%	\Delta_t {\rm Top}(\sfC) =  {\rm Top}(\tilde \sfC) - {\rm Top}(\sfC) \in \bbR^2. % \qquad \text{ and }\qquad
%%	\Delta_t {\rm Bottom}(\sfC) =  {\rm Bottom}(\tilde \sfC) - {\rm Bottom}(\sfC), 
%\end{align}
Define the same quantities for the IIC. 
More precisely, if~$\zeta_i^{{\rm Top}}$ is the IIC configuration used in the sampling of~$\omega_{t+1/2}$ around the top of~$\sfC$, 
and if we write~$\sfC^{\rm IIC}$ for the IIC cluster of~$\zeta_i^{{\rm Top}}$ and~$\tilde\sfC^{\rm IIC}$ for the corresponding IIC cluster of 
$(\bfS_{t+1} \circ \bfS_{t})(\zeta_i^{{\rm Top}})$, set 
\begin{align}
	\Delta_t^{\rm IIC} {\rm T}(\sfC) =  \langle {\rm Top}(\tilde \sfC^{\rm IIC}) - {\rm Top}(\sfC^{\rm IIC}), e_{\rm vert}\rangle.
\end{align}
When applying~$\bfS_{t}$ to~$\bbL_{\rm mix}$, each track of angle~$\beta$ is exchanged with the track of angle~$\alpha$ above it. 
As such, the vertices at the bottom of the  tracks of angle~$\beta$ in~$\bbL_{\rm mix}$ remain unchanged, and we can make sense of the image of any cluster containing at least one such point. 
In particular, there exists a well-defined notion of image of~$\sfC^{\rm IIC}$ after application of~$\bfS_t$. 
The same holds for~$\bfS_{t+1}$. 

Notice that~$\bfS_t(\bbL_{\rm mix})$ is simply a translate of~$\bbL_{\rm mix}$ and that the effect of~$\bfS_{t+1}$ on~$\bfS_t(\bbL_{\rm mix})$ 
is to again exchange each track of angle~$\beta$ with the track above it. 
As such, it would be tempting to think that the effect of the two transformations on the top of the IIC is identical. That is not the case, as the cells of~$\bbL_{\rm mix}$ are distinct from those of~$\bfS_t(\bbL_{\rm mix})$, which affects the definition of~${\rm Top}(\cdot)$. 
It is, however, the case that~$\bbL_{\rm mix}$ and~$(\bfS_{t+1} \circ \bfS_t)(\bbL_{\rm mix})$ are identical (up to translation), including in their partition in cells.

If no IIC configuration was used in the sampling of~${\rm Top}(\sfC)$ (for instance, because the top was frozen or because an overlap error occurred), 
then define~$\Delta_t^{\rm IIC} {\rm T}(\sfC)$ using some arbitrary~$\zeta_i^{{\rm Top}}$, different from the ones used for other clusters. 
Define~$\Delta_t^{\rm IIC} {\rm B}(\sfC)$,~$\Delta_t^{\rm IIC} {\rm L}(\sfC)$ and~$\Delta_t^{\rm IIC} {\rm R}(\sfC)$ analogously. 

If~${\rm Top}(\sfC)$ is active, set 
\begin{align}\label{eq:Delta_err_def}
	\Delta_t^{\rm err} {\rm T}(\sfC) =  \Delta_t {\rm T}(\sfC)  - \Delta_t^{\rm IIC} {\rm T}(\sfC);
\end{align}
if it is frozen, set~$\Delta_t^{\rm err} {\rm T}(\sfC) =0$. Define~$\Delta_t^{\rm err} {\rm B}(\sfC)$,~$\Delta_t^{\rm err} {\rm L}(\sfC)$ and~$\Delta_t^{\rm err} {\rm R}(\sfC)$ in the same way. 

Throughout this process, the variables~$\Delta_t^{\rm IIC} {A}(\cdot)$ are i.i.d., and therefore their sums behave like independent random walks. 
The error part may affect this dynamics, but will be proved to be small. We do so by a crude~$L^1$-bound.

\begin{proposition}\label{prop:error_exists}
	There exist universal constants~$c,C >0$ such that, for all~$\eta \geq N^{-c}$
	\begin{align}\label{eq:error_exists}
		\bbP\big[\exists \sfC \in {\rm Meso}(\omega_t) \text{ and } A \in \{{\rm T,B,L,R}\}\,: \, \Delta_t^{\rm err} {A}(\sfC) \neq 0 \big] \leq C \eta.
	\end{align}
\end{proposition}

The rest of this section is dedicated to the proof of Proposition~\ref{prop:error_exists}. 
To start, we state a bound on the number of mesoscopic clusters of any configuration~$\omega_t$.

\begin{lemma}\label{prop:mesoscopic_number}
	There exists~$C > 0$ such that, for any~$N$, any $0<\eta < 1$, and any~$t \geq0$, 
	\begin{align}\label{eq:mesoscopic_number}
		\phi_{\bbL_t}\big[\omega \text{ has more than~$\tfrac{\lambda}{\eta^{2}}$~$\eta$-mesoscopic clusters}\big] 
		\leq C e^{-\lambda/C} \qquad \text{ for all~$\lambda\geq 1$}.
	\end{align}
\end{lemma}

This is a direct consequence of the RSW theory and will appear immediate to experts.

\begin{proof}
	Fix~$N$,~$\eta$ and~$t$. We may restrict our attention to $\eta N$ sufficiently large, otherwise the statement is trivial. 
	Write~${\bf N}_{\rm meso}$ for the number of~$\eta$-mesoscopic clusters. 
	Notice that any such cluster needs to cross an annulus~$\Lambda_{\eta N/4}(x) \setminus \Lambda_{\eta N/8}(x)$ for some~$x \in \frac{\eta N}8\bbZ^2 \cap [-N,N]\times[0,N]~$ from the inside to the outside. For any such~$x$, write~${\bf N}_x$ for the number of such crossing clusters. 
	Then, 
	\begin{align}\label{eq:N_macro_sum}
		{\bf N}_{\rm meso} \leq \sum_x {\bf N}_x,
	\end{align}
	where the sum is over all~$x$ as above. There are~$O(\eta^{-2})$ terms in the sum. 
	
	It is a standard consequence of~\eqref{eq:RSW_iso} that the variables~${\bf N}_x$ have exponential tails. 
	Moreover, since~\eqref{eq:RSW_iso} is uniform in boundary conditions, the variables~${\bf N}_x$ may even be stochastically dominated by variables~$\tilde {\bf N}_x$ with exponential tails and which are independent for points at mutual distances larger than~$\eta N$. 
	The conclusion follows by splitting the sum in~\eqref{eq:N_macro_sum} into a bounded number of sums, each of which is bounded by i.i.d.\ variables~$\tilde{\bf  N}_x$ with exponential tails. 
\end{proof}

We now turn to the proof of Proposition~\ref{prop:error_exists} which is also based on the RSW theory and on the IIC-mixing estimate of Lemma~\ref{lem:IIC}.
There are no major difficulties here, but the proof is tedious due to the multiple sources of possible errors. 

\begin{proof}[Proof of Proposition~\ref{prop:error_exists}]
	The fact that $\Delta_t^{\rm err} {A}(\sfC) \neq 0$ for some~$\sfC \in {\rm Meso}(\omega_t)$ and~$A \in \{{\rm T,B,L,R}\}$ implies that one of the following must occur: 
	\begin{itemize}
	\item[(1)] an overlap error occurs at step~$t$, or
	\item[(2)] a coupling error occurs at step~$t$, or 
	\item[(3)] we have~$\Delta_t {\rm T}(\sfC) \neq \Delta_t^{\rm IIC} {\rm T}(\sfC)$ for some~$\sfC \in {\rm Meso}(\omega_t)$, even though the environment in~$\Lambda_r({\rm Top}(\sfC))$ is given by the corresponding IIC configuration, or the same for~${\rm B}, {\rm L}$ or~${\rm R}$ instead of~${\rm T}$. 
	\end{itemize}
	We analyse each scenario separately.% in the case of $\Delta_t{\rm T}$ first. \smallskip 

	\noindent{\bf (1)}
	A relatively straightforward analysis based on~Proposition~\ref{prop:arms_bounds_iso} shows that 
	\begin{align}\label{eq:overlap_err}
		\bbP[\text{overlap error occurs at step~$t$}] \leq C \eta, 
	\end{align}
	for some universal constant~$C >0$. 
	Indeed, pave the observation window by overlapping boxes~$\{\Lambda_{4R}(z): z \in R\bbZ^2 \cap [-N,N] \times [0,N]\}$. 
	For an overlap error to occur, at least one of the following needs to occur in $\omega_t$:
	\begin{itemize}
	\item[(1.a)] there exists a box~$\Lambda_{4R}(z)$ in the mixed block 
	with three arms in a half-plane and at least one additional primal arm to a distance~$\eta N/2 - 8R$, or 
	\item[(1.b)]  there exists a box~$\Lambda_{4R}(z)$ that intersects two distinct blocks of the lattice and has three arms in a half-plane to a distance~$\eta N - 4R$.
	\end{itemize} 
	In the above, by ``three-arm in a half-plane'' we mean the events in~\eqref{eq:3hp_arm} or~\eqref{eq:3hp_arm2}.
	There are~$O(N/R)^2$ potential boxes in the first category and~$O(N/R)$ in the second. 
	By a union bound and the bounds on the 1-arm and 3-arm probabilities from Proposition~\ref{prop:arms_bounds_iso}, we may show 
	that~\eqref{eq:overlap_err} holds as long as $\eta >N^{-c}$  for appropriately chosen universal constants. We give details below.
	
	Start by bounding the probability of the events comprising three arms in the lower half-plane.
	For any given $z$, due to \eqref{eq:one_arm_bound_iso} and  \eqref{eq:3hp_arm},
	\begin{align}
	\phi_{\bbL_t} \big[\text{(1.a) occurs for $\Lambda_{4R}(z)$ with 3 arms in lower half-plane}\big]  \leq C_0\big(\tfrac{R}{\eta N}\big)^{2 + \onearmbound},
	\end{align}
	for some universal constant $C_0 >0$ (henceforth, $C_i: i\geq 0$ denote universal constants).
	Applying a union bound over the $O(N/R)^2$ potential boxes in the first category, we find 
	\begin{align}\label{eq:1a_sum_lhp}
	\phi_{\bbL_t} \big[\text{(1.a) occurs  with 3 arms in lower half-plane}\big]  \leq C_1 \big( \tfrac{R}{\eta N}\big)^{\onearmbound}\eta^{-2} \leq C_1 \eta,
	\end{align}
	where the last inequality is due to the choice~\eqref{eq:Rchoice}  of $R$.
	
	We turn to (1.b) with three arms in the lower half-plane. Here, there are $O(N/R)$ potential boxes, and \eqref{eq:3hp_arm} implies 
	\begin{align}
	\phi_{\bbL_t} \big[\text{(1.b) occurs  with 3 arms in lower half-plane}\big]  \leq C_2  \tfrac{R}{\eta^2 N} \leq C_2 \eta,
	\end{align}
	again due to the choice of $R$. \smallskip 
	
	The same arguments apply for when the three arms occur in the upper half-plane. 
	When three arms occur in the lateral half-planes, a more delicate analysis is needed due to the weaker bound \eqref{eq:3hp_arm3}. 
	The same argument as above applies for (1.b) also in the lateral case. 
	
	Finally, we turn to bounding the probability of (1.a) for arms in the left half-plane. 
	Fix $z \in R\bbZ^2 \cap [-N,N] \times [0,N]$ and  let $d$ be the distance between $z$ and the nearest of the two other blocks (if this distance is smaller than $R$, set $d = R$). 
	For the event in (1.a) to occur, 
	there must exist three-arms in the left half-plane plus an additional primal arm from $\Lambda_{4R}(z)$ to distance $d$, 
	and the same between distances $d$ and $\eta N/2-8R$. 
	The first event is entirely contained in a lattice of the type $\bbL_{\rm mix}$. 
	Using \eqref{eq:mixing},~\eqref{eq:3hp_arm2} and~\eqref{eq:3hp_arm3}, we conclude that 
	\begin{align*}
	\phi_{\bbL_t} \big[\text{(1.a) occurs for  $\Lambda_{4R}(z)$ with 3 arms in left half-plane}\big]  
	\leq C_3 \big(\tfrac{R}{d}\big)^{2+\xi_1} \times \big(\tfrac{d}{\eta N}\big)^{1+c}.	
	\end{align*}
	Summing the above over all possible values of $z$ yields the same upper bound as in \eqref{eq:1a_sum_lhp}. This concludes the proof of \eqref{eq:overlap_err}. The same works for arms in the right half-plane.
	
	\medskip 
	
	\noindent{\bf (2)} Lemma~\ref{lem:IIC} implies that
\begin{align}
	\bbP[{\rm err}_i = 1] = d_{\rm TV} \leq C\big(\tfrac{r}{R}\big)^{c_{\rm IIC}} \leq C \eta^{3},
\end{align}
where the last inequality is due to the choice of~$r$. 
Combining the above with Lemma~\ref{prop:mesoscopic_number}, we conclude that 
\begin{align}\label{eq:coupling_err}
	\bbP[\text{coupling error occurs at step~$t$}] \leq C \eta
\end{align}
for some universal constant~$C > 0$.\medskip

	\noindent{\bf (3)} We start by analysing the case of ${\rm Top}$. 
	This type of error may appear for two reasons: 
	\begin{itemize}
	\item[(3.a)] the configuration in~$\Lambda_{r/2}({\rm Top}(\sfC))$ after applying~$\bfS_{t+1}\circ \bfS_t$ is not identical to the one given by the IIC, or 
	\item[(3.b)] the top of the cluster~$\sfC$ or that of the corresponding IIC cluster shifts by more than~$r/2$ to the left or to the right when applying~$\bfS_{t+1}\circ \bfS_t$, thus exiting~$\Lambda_{r/2}({\rm Top}(\sfC))$.
	\end{itemize}
	
	Scenario (3.a) is very unlikely to occur.   Indeed, as explained in Section~\ref{sec:star-triangle}, the track exchanges are local, due to the occurrence of the events $F_k$ depicted in Figure~\ref{fig:track_exchange}. By the finite energy property, for any given cluster~$\sfC$, scenario (3.a) has a probability bounded above by~$C e^{-c r}$ for universal constants~$C,c$. Performing a union bound over the mesoscopic clusters of~$\omega_t$ and using Lemma~\ref{prop:mesoscopic_number}, we conclude that the probability of (3.a) is bounded above by~$C \eta^{-2}e^{-c r}$ for potentially modified universal constants~$C,c$.
	
	For scenario (3.b) to occur, an almost-top of some mesoscopic cluster~$\sfC$ needs to exist at a distance at least~$r/2$ from~${\rm Top}(\sfC)$ in~$\omega_t$, or the equivalent for the IIC clusters; see Figure~\ref{fig:controlling_err}. 
	By a union bound over the possible positions of~${\rm Top}(\sfC)$ and that of the almost-top, we may use~\eqref{eq:mixing} and~\eqref{eq:3hp_arm} to bound the probability of such an event. We find % that the probability of the scenario (3.b) in the case of top is bounded above by a multiple of
	\begin{align}\label{eq:(3.b)} 
	\phi_{\bbL_t} \big[\text{(3.b) occurs for  some ${\rm Top}(\sfC)$}\big] 
	\leq C_0\eta^{-2}  \sum_{k\geq r}k^{-2} = C_1 \eta^{-2} r^{-1},
	\end{align}
	for universal constants $C_0,C_1$. 
	Recall that~$r$ was fixed above \eqref{eq:iic_err_couple} to be equal to~$\eta^C N$ for some fixed constant~$C$.
	Thus, the probability above is bounded by~$\eta^{-2-C}/N$. % \leq \eta$, provided $\eta > N^{-\frac{1}{3+C}}$.  
	This concludes the bound for the error of type (3) for the top of clusters. The same holds for the bottom. 
	\smallskip 	
	
	We turn to bounding the probability of (3) for the ``lateral'' directions, and focus on ${\rm Right}$. 
	For (3.a) the same argument as for ${\rm Top}$ applies. 
	For (3.b) a more careful analysis relying on~\eqref{eq:3hp_arm2} and~\eqref{eq:3hp_arm3} instead of~\eqref{eq:3hp_arm} is needed. 
	It is similar to that for (1.a) and we only sketch it below.
For scenario (3.b) to occur for $\rm Right$, an {\em almost-right} of some mesoscopic cluster~$\sfC$ needs to exist at a distance at least~$r/2$ from~${\rm Right}(\sfC)$ in~$\omega_t$, or the equivalent for the IIC clusters. 
	For any two points $z,z'$, we express the probability that they are the right and  almost-right of some mesoscopic cluster in terms of the distance $k$ between $z$ and $z'$, but also the distances $d$ and $d'$ between the points and the other blocks of the lattice. We do so using~\eqref{eq:mixing},~\eqref{eq:3hp_arm2} and~\eqref{eq:3hp_arm3}, as in the corresponding bound in case (1.a). 
	Summing over all possible positions of $z$ and $z'$, we find that \eqref{eq:(3.b)}  also applies to ${\rm Right}$.

	\begin{figure}
	\begin{center}
	\includegraphics[width = .8\textwidth]{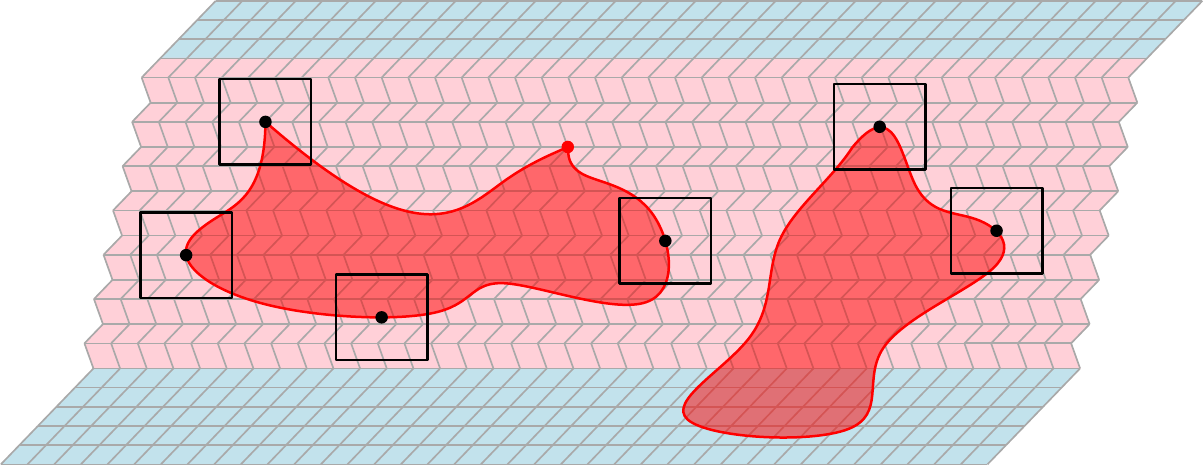}
	\caption{A step in the transformations~$\bbL_t$. The mixed block is pink, the two other ones are blue; 
	notice that the bottom block is of angle~$\alpha$, which indicates that~$t > K$.
	The two red clusters are mesoscopic; the extrema are marked, as are the extremum boxes. 
	All extrema of the left cluster are active, for the right one, only the top and right are active, the other are frozen. 
	\newline
	Notice the potential error for the move of the top of the left cluster: there exists a second point (marked in red) in the cluster
	that may move up and overtake the top.}
	\label{fig:controlling_err}
	\end{center}
	\end{figure}
	
	\medbreak
	We have now estimated the probabilities of case (1), (2) and (3) and are ready to  conclude.
	Summing the different potential sources of error gives	\begin{align*}
		\bbP\big[\exists \sfC \in {\rm Meso}(\omega_t) \text{ and } A \in \{{\rm T,B,L,R}\}\,: \, \Delta_t^{\rm err} {A}(\sfC) \neq 0 \big]&
		\leq C (\eta + C \eta^{-2} e^{-c r} + \eta^{-2-C}/N) \\
		&\leq 3 C \eta,
	\end{align*}
	provided that~$\eta \geq N^{-c_0}$ for some fixed (small) constant~$c_0 >0$. 
\end{proof}

\subsection{Stability of mesoscopic clusters}\label{sec:stability_meso}

Write
\begin{align}\label{eq:drift}
 	{\rm Drift}_{\rm vert} &= \tfrac12 \bbE[\Delta_t^{\rm IIC}{\rm T}] = \tfrac12 \bbE[\Delta_t^{\rm IIC}{\rm B}],\\
	{\rm Drift}_{\rm lat} &=  \tfrac12\bbE[\Delta_t^{\rm IIC}{\rm L}] =  \tfrac12\bbE[\Delta_t^{\rm IIC}{\rm R}].
\end{align} 
The prefactor~$1/2$ appears since~$\Delta_t^{\rm IIC}{\rm T}$ is an increment over two time steps of the process~$(\omega_t)_{t\geq 0}$.
The second equality on each line is a simple consequence of symmetry. 
Moreover, we have
\begin{align}-\sin\beta < {\rm Drift}_{\rm vert} < \sin \alpha.\end{align} Indeed, 
the non-strict inequalities are immediate: applying~$\bfS_{1} \circ \bfS_{0}$ to~$\phi^{\rm IIC,T}_{\bbL_{\rm mix}}$ the top of the IIC cluster may increase or decrease by at most one cell in the vertical direction. The height of a cell of~$\bbL_{\rm mix}$ is~$\sin\alpha + \sin\beta$ and the transformation~$\bfS_{1} \circ \bfS_{0}$ shifts all cells vertically by~$\sin\alpha - \sin\beta$. Thus~$-2\sin \beta \leq \Delta_t^{\rm IIC}{\rm T}\leq 2\sin \alpha$ almost surely. The finite energy property allows us to derive the strict inequalities. %The same holds for ${\rm Bottom}$.

Fix a mesoscopic cluster~$\sfC$ of~$\omega_0$.
For~$t >0$,  assuming that the cluster remains mesoscopic in each~$\omega_s$ with~$0\leq s \leq t$, write~$\sfC_t$ for the cluster in~$\omega_t$ and set
\begin{align}
	{A}_t(\sfC)  = {A}(\sfC_t) \qquad \text{ for~$A \in \{{\rm T,B,L,R}\}$}. 
\end{align}

Write~$\tau_{\rm start}^{\rm T}(\sfC)$ for the first time that~${\rm Top}(\sfC_t)$ is not frozen. 
Write~$\tau_{\rm meso}(\sfC)$ for the first time~$t$ that~$\sfC_t$ ceases to be mesoscopic. 
Finally, set~$\tau_{\rm end}^{\rm T}(\sfC)$ for the first time~$t > \tau_{\rm start}^{\rm T}(\sfC)$ such that, either~$\sfC_t$ is not mesoscopic or~${\rm Top}(\sfC_t)$ is frozen again. 
The evolution of~${\rm T}_t(\sfC)$ is then constant up to~$\tau_{\rm start}^{\rm T}(\sfC)$ and after~$\tau_{\rm end}^{\rm T}(\sfC)$ (the latter follows by convention if~$\tau_{\rm end}^{\rm T}(\sfC) = \tau_{\rm meso}(\sfC)$). 

Apply the same construction for~${\rm Bottom}(\sfC)$ to define~$\tau_{\rm start}^{\rm B}(\sfC)$ and~$\tau_{\rm end}^{\rm B}(\sfC)$. 
Note that $\tau_{\rm start}^{\rm T}(\sfC)$, $\tau_{\rm start}^{\rm B}(\sfC)$,
$\tau_{\rm end}^{\rm T}(\sfC)$, $\tau_{\rm end}^{\rm B}(\sfC)$ and~$\tau_{\rm meso}(\sfC)$ are determined by the process $(\calH(\omega_t))_t$, with the first two being even determined by its initial value:
\begin{align}
\tau_{\rm start}^{\rm T}(\sfC) = K -\frac{ {\rm T}_0(\sfC)}{\sin\beta} \quad\text{ and }\quad \tau_{\rm start}^{\rm B}(\sfC) = K - \frac{{\rm B}_0(\sfC)}{\sin\beta}.
\end{align}
Here, and throughout the rest of the proof, we ignore integer parts in the definition of time-steps.

Let us start by considering a simplified situation in which the incremental steps of~$({\rm T}_t(\sfC), {\rm B}_t(\sfC), {\rm L}_t(\sfC), {\rm R}_t(\sfC))$ are deterministic and exactly given by~${{\rm Drift}}_v$ and~${{\rm Drift}}_h$, respectively, and that~$\tau_{\rm meso}(\sfC)= \infty$. Assuming this, we would find that, for~$A \in \{{\rm T}, {\rm B}\}$,
\begin{align}\label{eq:tau_drift}
A_t(\sfC) &= A_0(\sfC) + \big((t \wedge \tau_{\rm end}^{{\rm Drift}, A} (\sfC) - \tau_{\rm start}^A(\sfC))\vee 0\big)\cdot {\rm Drift}_{\rm vert},\end{align}
where
\begin{align}
%{\rm T}_t(\sfC) &= T_0(\sfC) + \big((t \wedge \tau_{\rm end}^{{\rm Drift}, {\rm T}} (\sfC) - \tau_{\rm start}^{\rm T}(\sfC))\vee 0\big) {{\rm Drift}}_v \text{ and }\\
%{\rm B}_t(\sfC) &= B_0(\sfC) + \big((t \wedge \tau_{\rm end}^{{\rm Drift}, {\rm B}} (\sfC) - \tau_{\rm start}^{\rm B}(\sfC))\vee 0\big) {{\rm Drift}}_v 
	\tau_{\rm end}^{{\rm Drift}, A} (\sfC) &= {A}_0(\sfC) \frac{1 + {\rm Drift}_{\rm vert} /\sin \beta}{\sin\alpha - {\rm Drift}_{\rm vert}} + K % \qquad \text{ for~$A \in \{{\rm T}, {\rm B}\}$}
%	\text{ and } 
%	\tau_{\rm end}^{{\rm Drift}, {\rm B}} (\sfC) = {\rm B}_0(\sfC) \frac{1 + {\rm Drift}_{\rm vert} /\sin \beta}{\sin\alpha - {\rm Drift}_{\rm vert}} + K 
\end{align} 
is the equivalent~$\tau_{\rm end}^{A}$ under the assumption that~$A_t$ moves deterministically with an increment of~${\rm Drift}_{\rm vert}$ at every step when~$A$ is active. 
%In particular, the top of~$\sfC$ would be active between
%$\tau_{\rm start}^{\rm T}(\sfC)$
%and~$\tau_{\rm end}^{{\rm Drift}, {\rm T}}(\sfC)$,
%so for a total of~${\rm T}_0(\sfC) \frac{\sin\alpha + \sin\beta}{(\sin \alpha -{\rm Drift}_{\rm vert})\sin\beta}$  transformations. 

For~${\rm L}_t(\sfC)$ and~${\rm R}_t(\sfC)$, the dynamics, even under the assumption that their steps are deterministic, is not that straightforward. This is because it is not clear when exactly~${\rm Left}(\sfC_t)$ and~${\rm Right}(\sfC_t)$ are frozen or active. 
Nevertheless, we can expect that
\begin{align*}
\big|{\rm L}_t(\sfC) - {\rm L}_0(\sfC) -  \big((t \wedge \tau_{\rm end}^{{\rm Drift}, {\rm T}} (\sfC) - \tau_{\rm start}^{\rm T}(\sfC))\vee 0\big)\cdot {\rm Drift}_{\rm lat} \big| \leq 
C \big({\rm T}_0(\sfC) -  {\rm B}_0(\sfC)\big),
%C \big(\tau_{\rm start}^{\rm T}(\sfC) -   \tau_{\rm start}^{\rm B}(\sfC) \big),
\end{align*}
for some constant~$C$. This is because the ambiguity in whether~${\rm Left}(\sfC_t)$ is frozen or active only occurs when~${\rm Top}(\sfC_t)$ and~${\rm Bottom}(\sfC_t)$ are in different blocks of~$\bbL_t$, and this only occurs for~$O({\rm T}_0(\sfC) -  {\rm B}_0(\sfC))$ steps of the process. The same holds for~${\rm R}_t(\sfC)$. 

With this in mind, we introduce the following definition.

\begin{definition}For constants~$c,C> 0$ to be fixed below, we say that a cluster~$\sfC \in {\rm Meso}(\omega_0)$ is {\em$(c,C)$-stable} if $\tau_{\rm meso}(\sfC) \geq K+K'$ and, for all~$0 \leq t \leq K+K'$ and~$A \in \{{\rm T,B,L,R}\}$, 
\begin{align}\label{eq:stable_C}
	&\big|A_t(\sfC) - A_0(\sfC) -  \big((t \wedge \tau_{\rm end}^{{\rm Drift}, {\rm T}} (\sfC) - \tau_{\rm start}^{\rm T}(\sfC))\vee 0\big)\cdot {\rm Drift}_* \big| \leq C \eta^c N ,
\end{align}
where~${\rm Drift}_*$ is~${\rm Drift}_{\rm vert}$ for~$A = {\rm T,B}$ and~${\rm Drift}_{\rm lat}$ for~$A = {\rm L,R}$.
\end{definition}
 Note that we use~$\tau_{\rm end}^{{\rm Drift}, {\rm T}}(\sfC)$ for all directions, even the bottom one. This will prove irrelevant as the difference between~$\tau_{\rm end}^{{\rm Drift}, {\rm T}}(\sfC)$ and~$\tau_{\rm end}^{{\rm Drift}, {\rm B}}(\sfC)$ will be absorbed in the error term~$C \eta^c N$. 

\begin{remark}\label{rem:M_origin}
	The link between the definition~\eqref{eq:stable_C} of stable cluster and the map~$M_{\beta,\alpha}$ may appear mysterious. 
	To understand this, consider~$x =(x_1,x_2) \in \bbR \times \bbR_+$ such that 
	\begin{align}
	\langle x,e_{\rm vert}\rangle  =  {\rm T}_0(\sfC) \qquad\text{and}\qquad
	\langle x,e_{\rm lat}\rangle   =  {\rm R}_0(\sfC)
	\end{align} 
	for some cluster~$\sfC \in {\rm Meso}^+(\omega_0)$. 
	Assume in addition that~$\sfC$ is~$(c,C)$-stable and write~$y =(y_1,y_2)$ for the unique point such that 
   \begin{align}
   	\langle y,e_{\rm vert}\rangle  =  {\rm T}_{ K+K'}(\sfC)\qquad\text{and}\qquad
	\langle y,e_{\rm lat}\rangle   =  {\rm R}_{ K+K'}(\sfC). 
   \end{align}
	Then, a straightforward computation shows that 
	\begin{align}
	y_1 &= x_1 +    \big({\rm Drift}_{\rm lat} \tfrac{1}{\sin\beta}  +  {\rm Drift}_{\rm vert} \tfrac{\cos\beta}{\sin\beta} \big)\cdot x_2 \tfrac{\sin\alpha + \sin\beta}{(\sin \alpha -{\rm Drift}_{\rm vert})\sin\beta} + O(\eta^c N),  \\
	y_2 &= x_2 + {\rm Drift}_{\rm vert} \cdot x_2 \tfrac{\sin\alpha + \sin\beta}{(\sin \alpha -{\rm Drift}_{\rm vert})\sin\beta}+ O(\eta^c N)
	\end{align}
	which translates to~$y^{T} = M_{\beta,\alpha}x^{T} + O(\eta^c N)$. 
	
	Indeed, the term
	$$x_2 \frac{\sin\alpha + \sin\beta}{(\sin \alpha -{\rm Drift}_{\rm vert})\sin\beta} = \tau_{\rm end}^{{\rm Drift}, {\rm T}} (\sfC) - \tau_{\rm start}^{{\rm T}} (\sfC)$$ is the approximate number of transformations that affect the cluster~$\sfC$. 
	For the vertical coordinate, observe that, on average, each transformation pushes the cluster up by~${\rm Drift}_{\rm vert}$.
	Similarly, due to the relation between the canonical coordinates and the coordinates obtained from $({\rm T}_{ t}(\sfC),{\rm R}_{ t}(\sfC))$, 
	each transformation, on average, pushes the cluster to the right by~${\rm Drift}_{\rm lat} \tfrac{1}{\sin\beta}  +  {\rm Drift}_{\rm vert} \tfrac{\cos\beta}{\sin\beta}$. 
\end{remark}

We are now in a position to discuss a key result which states that the dynamics only affects the extrema of mesoscopic clusters by the linear map~$M_{\beta,\alpha}$. 
This will be later shown to imply that the actual shapes of macroscopic clusters are only affected through~$M_{\beta,\alpha}$. 
We cannot expect the statement below to apply to all mesoscopic clusters of~$\omega_0$, since clusters that are barely mesoscopic in~$\omega_0$ may cease being mesoscopic during the process. For that reason, we introduce the more restrictive notion of~${\rm Meso}^+(\omega_0)$.

Fix some large constant~$C_{M}> 2 + 2 \max\{|{\rm Drift}_{\rm vert}|, |{\rm Drift}_{\rm lat}|\}>0$. 
Let ${\rm Meso}^+(\omega_0)$ be the set of {\em mesoscopic} clusters~$\sfC$ of~$\omega_0$ such that
\begin{align}\label{eq:Meso^+}
	C_{M}\eta N\leq {\rm T}_0(\sfC)  - {\rm B}_0(\sfC)  &\leq \tfrac{\sqrt \eta}{C_{M}} N,\nonumber\\
	{\rm T}_0(\sfC)&\leq \tfrac{1}{C_{M}}N,\nonumber\\
	 - \tfrac{1}{C_{M}}N\le {\rm L}_0(\sfC)\le {\rm R}_0(\sfC)&\leq \tfrac{1}{C_{M}}N.\qquad
\end{align}
The only difference from standard mesoscopic clusters is the presence of the constant $C_M$, which makes the defining conditions slightly more restrictive while giving us some room in subsequent estimates.

The next proposition states that all mesoscopic clusters are stable with probability close to 1.
\begin{proposition}\label{prop:stability_meso}
	There exist constants~$c,C >0$ such that for~$N$ large enough and~$\eta > N^{-c}$,
	\begin{align}\label{eq:stable}
	\bbP[\text{all~$\sfC \in {\rm Meso}^+(\omega_0)$ are~$(c,C)$-stable}] \geq 1 - C \eta^{c}.
	\end{align}
\end{proposition}

% The use of the constant~$C_M$ in the restrictions of~\eqref{eq:Meso^+} is so that the clusters of~${\rm Meso}^+(\omega_0)$ do not exit the observation window during our process due to the drift. 

\begin{proof}
	For any cluster~$\sfC \in {\rm Meso}(\omega_0)$ and~$A \in \{{\rm T}, {\rm B}\}$, 
	as long as~$t \leq \tau_{\rm end}^A(\sfC)$,
	\begin{align}\label{eq:Ext_t}
		{A}_t(\sfC) =  	{A}_0(\sfC)  + 
		  \sum_{\tau_{\rm start}^A(\sfC) \leq s < t}\Delta_s^{\rm IIC} {A}(\sfC) + \sum_{\tau_{\rm start}^A(\sfC) \leq s < t} \Delta_s^{\rm err} {A}(\sfC).
	\end{align}
	(Formally, the above is only valid for~$t$ even, and the sum is over~$s$ even.) Set
	\begin{align}
		{A}^{\rm IIC}_t(\sfC) =  	{A}_0(\sfC)  +  \sum_{\tau_{\rm start}^A(\sfC) \leq s < t}\Delta_s^{\rm IIC} {A}(\sfC). 
	\end{align}
	Extend the notation to~$A  \in \{{\rm L}, {\rm R}\}$ by  setting 
	\begin{align}\label{eq:Ext_t}
		{A}^{\rm IIC}_t(\sfC) =  	{A}_0(\sfC)  +  \sum_{\tau_{\rm start}^{\rm T}(\sfC) \leq s < t}\Delta_s^{\rm IIC} {A}(\sfC).
	\end{align}
	Notice that in the latter definition, we use~$\tau_{\rm start}^{\rm T}(\sfC)$ as a starting time, even though~${\rm L}_t(\sfC)$ may be constant beyond this time. 
	Also, for technical reasons, we continue~$	{A}^{\rm IIC}_t(\sfC)$ beyond~$\tau_{\rm end}^{.}(\sfC)$ with IIC increments which are irrelevant for the true process~${A}_t(\sfC)$.
	
	Define 
	\begin{align}
	{\rm TotErr} := \sum_{t = 0}^{K+K'} \ind[\{\exists \sfC \in {\rm Meso}(\omega_t) \text{ and } A \in \{{\rm T,B,L,R}\}\,: \, \Delta_t^{\rm err} {A}(\sfC) \neq 0 \}]
	\end{align}
	for the total number of steps at which errors occur. 
	Recall that the variables~$\Delta_t^{\rm err} {A}$ for~$A \in \{{\rm T,B,L,R}\}$ are a.s.\ bounded. 
	Thus, for some universal constant~$C$, 
	\begin{align}\label{eq:TotErr}
		\big|{A}_t(\sfC) -  {A}_t^{\rm IIC}(\sfC)\big|\leq C \cdot {\rm TotErr} \qquad  \text{ for all~$0 \leq t \leq \tau_{\rm end}^A(\sfC)$ and
		$A \in \{{\rm T},{\rm B}\}$}.%~$A \in \{{\rm T},{\rm B},{\rm L},{\rm R}\}$}.
	\end{align}
	The case of lateral extrema is slightly more delicate and will be handled later, see \eqref{eq:TotErrLR}.
	
	Summing~\eqref{eq:error_exists} and applying the Markov inequality, we find 
	\begin{align}\label{eq:TotErr2}
		\bbP\big[ {\rm TotErr} \leq C\, \eta^c N \big] \geq 1- \eta^c
	\end{align}
	for some universal constants~$c,C > 0$, under the condition of Proposition~\ref{prop:error_exists} for $\eta$. 
	
	The coordinates~${A}_t^{\rm IIC}(\sfC)$ are independent random walks with constant drift. 
	Define the event
	\begin{align*}
		{\rm GoodRW} = 
		\left\{\begin{array}{c}\forall \sfC \in {\rm Meso}^+(\omega_0), \, A \in \{ {\rm T}, {\rm B}, {\rm L}, {\rm R}\}  \text{ and } \tau_{\rm start}^{*}(\sfC)\leq t\leq K+K',\\
		|{A}_t^{\rm IIC}(\sfC) -{A}_0^{\rm IIC}(\sfC) - (t - \tau_{\rm start}^{*}(\sfC)) \cdot {\rm Drift}_{A}| \leq N^{2/3} \end{array}\right\},
	\end{align*}
	where~$* = {\rm T}$ if~$A \in \{ {\rm T}, {\rm L}, {\rm R}\}$ and~$* = {\rm B}$ if~$A= {\rm B}$,
	and~${\rm Drift}_{A} =  {\rm Drift}_{\rm vert}$ when~$A \in \{{\rm T,B}\}$ and~$ {\rm Drift}_{A} =  {\rm Drift}_{\rm lat}$ when~$A \in \{{\rm L,R}\}$. Note that, as opposed to the notion of stability~\eqref{eq:stable_C}, we choose to work here with~$\tau_{\rm start}^{\rm T}$  for~${\rm T}^{\rm IIC}_t(\sfC)$  and~$\tau_{\rm start}^{\rm B}$ for~${\rm B}^{\rm IIC}_t(\sfC)$. This extra precision ensures that the vertical diameter of~$\sfC$ does not degenerate throughout the process.
		
	Standard random walk estimates %\im{Source} 
	combined with crude bounds on~$|{\rm Meso}^+(\omega_0)|$ provided by~\eqref{eq:mesoscopic_number} imply that 
	\begin{align}\label{eq:GoodRW}
		\bbP\big[ {\rm GoodRW}\big] \geq 1 - C \eta.
	\end{align}
	for any~$\eta > N^{-c}$ and~$c,C > 0$ some universal constants.

	Finally, define 
	\begin{align*}
		{\rm NoDieIIC} = \left\{\begin{array}{c} \forall \sfC \in {\rm Meso}^+(\omega_0) \text{ and } 0\leq t\leq K+K',\\
		 2 \eta  N \leq {\rm T}_t^{\rm IIC}(\sfC) - {\rm B}_t^{\rm IIC}(\sfC) \leq \tfrac12 \eta^{1/2} N, \\
		{\rm T}_t^{\rm IIC}(\sfC) \leq \tfrac1{2} N, \\
		 	-\tfrac1{2} N\le {\rm L}_t^{\rm IIC}(\sfC)\le{\rm R}_t^{\rm IIC}(\sfC) \leq \tfrac1{2} N
		 \end{array} \right\} 
	\end{align*}
	In other words,~${\rm NoDieIIC}$ is the event that the coordinates of the fictitious IIC processes
	$({\rm T}_t^{\rm IIC}, {\rm B}_t^{\rm IIC}, {\rm L}_t^{\rm IIC}, {\rm R}_t^{\rm IIC})$ for  $\sfC \in {\rm Meso}^+(\omega_0)$ satisfy the conditions of~\eqref{eq:meso}, 
	with constants twice as restrictive, for all~$0\leq t\leq K+K'$.
	Due to the conditions~\eqref{eq:Meso^+} for~$\sfC$ to be in~${\rm Meso}^+(\omega_0)$
	 (specifically the choice of~$C_{M}$ as sufficiently large) 
	we find that, whenever~${\rm GoodRW}$ occurs, so does~${\rm NoDieIIC}$.

	Assume henceforth that both~${\rm GoodRW}$ and~${\rm NoDieIIC}$ occur, that~${\rm TotErr} \leq C\, \eta^c N$ and that~$N$ is sufficiently large.
	We will prove that the event in~\eqref{eq:stable} occurs. 

	We start by proving that
	\begin{align}\label{eq:meso+survive}
		\text{$\tau_{\rm meso}(\sfC)> K+K'$\quad for all~$\sfC \in{\rm Meso}^+(\omega_0)$}.
	\end{align}
	Indeed, for any~$\sfC \in{\rm Meso}^+(\omega_0)$,~${\rm T}_t(\sfC)$ and~${\rm B}_t(\sfC)$ are close to their IIC counterparts due to~\eqref{eq:TotErr} and our assumption on~${\rm TotErr}$. Factoring in the occurrence of~${\rm NoDieIIC}$, we conclude that, for all~$0\leq t\leq K+K'$,
	\begin{align} \label{eq:meso+survive2}
		 \eta  N \leq {\rm T}_t(\sfC) - {\rm B}_t(\sfC) \leq \eta^{1/2} N  \quad\text{ and }\quad {\rm T}_t(\sfC) \leq  N.
	\end{align}
	The lateral coordinates~${\rm L}_t(\sfC)$ and~${\rm R}_t(\sfC)$ require slightly more attention, due to the uncertain times when they evolve. 
	However, we know they are stationary before~$\tau_{\rm start}^{\rm T}(\sfC)$ and after~$\tau_{\rm end}^{\rm T}(\sfC)$ and 
	we know they are evolving at all times between~$\tau_{\rm start}^{\rm B}(\sfC)$ and~$\tau_{\rm end}^{\rm B}(\sfC)$.
	The upper bound on~${\rm T}_t(\sfC) - {\rm B}_t(\sfC)$ in~\eqref{eq:meso+survive2} 
	ensures that there are at most~$O(\sqrt \eta N)$ steps when~${\rm L}_t^{\rm IIC}(\sfC)$ moves, but~${\rm L}_t(\sfC)$ does not, 
	and the same for~${\rm R}_t(\sfC)$. 
%	The upper bound on~${\rm T}_0(\sfC) - {\rm B}_0(\sfC)$ in~\eqref{eq:Meso^+} and the occurrence of~${\rm GoodRW}$ 	
%	ensure that there are at most~$O(\sqrt \eta N)$ steps when~${\rm L}_t^{\rm IIC}(\sfC)$ moves, but~${\rm L}_t(\sfC)$ does not, 
%	and the same for~${\rm R}_t(\sfC)$. 
	This adds an error term when comparing the lateral coordinates of the clusters~$\sfC \in {\rm Meso}^+(\omega_0)$ to their IIC equivalents, 
	but this error is sufficiently small that it allows us to deduce from~${\rm NoDieIIC}$ that
	\begin{align}
	{\rm L}_t(\sfC) \geq - N \text{ and } {\rm R}_t(\sfC) \leq  N \quad \text{ for all~$0\leq t\leq K+K'$,}
	\end{align}
	which, together with~\eqref{eq:meso+survive2}, implies~\eqref{eq:meso+survive}.\smallskip 
	
	We now turn to the stability of the clusters of~${\rm Meso}^+(\omega_0)$. 
	By the definition of~$\tau_{\rm end}^{\rm T} (\sfC)$, we have that 
	\begin{align}
		{\rm T}_{\tau_{\rm end}^{\rm T} (\sfC) }(\sfC) = (\tau_{\rm end}^{\rm T} (\sfC)  - K)\sin\alpha.
	\end{align}
	Comparing~${\rm T}_{\tau_{\rm end}^{\rm T} (\sfC) }(\sfC)$ to~${\rm T}^{\rm IIC}_{\tau_{\rm end}^{\rm T} (\sfC) }(\sfC)$ using \eqref{eq:TotErr}, 
	with the latter term being close to~${\rm T}_0(\sfC) +(\tau_{\rm end}^{\rm T} (\sfC)  - \tau_{\rm start}^{\rm T}(\sfC)) \cdot {\rm Drift}_{\rm vert}$ by ${\rm GoodRW}$, we conclude that 
	\begin{align*}
		{\rm T}_0(\sfC) + (\tau_{\rm end}^{\rm T} (\sfC) - \tau_{\rm start}^{\rm T}(\sfC)) \cdot {\rm Drift}_{\rm vert} - (\tau_{\rm end}^{\rm T} (\sfC)  - K)\sin\alpha = O(\eta^c N ).
	\end{align*}
	Recall that~$(K - \tau_{\rm start}^{\rm T}(\sfC)) \sin \beta = {\rm T}_0(\sfC)$, which implies that 
	\begin{align}\label{eq:tau_Drift_RW}
		\tau_{\rm end}^{\rm T} (\sfC) = {\rm T}_0(\sfC) \frac{1 + {\rm Drift}_{\rm vert}/ \sin \beta}{\sin\alpha - {\rm Drift}_{\rm vert}} + K + O(\eta^c N)
		= \tau_{\rm end}^{{\rm Drift, T}}(\sfC)+ O(\eta^c N).
	\end{align}
	Finally, since~${\rm T}_t(\sfC)$ is stationary after~$\tau_{\rm end}^{\rm T}(\sfC)$, for all~$0 \leq t\leq K+K'$,
	\begin{align*}
		\big|{\rm T}_t(\sfC) - {\rm T}_0(\sfC)&   -  \big((t \wedge \tau_{\rm end}^{{\rm Drift}, {\rm T}} (\sfC) - \tau_{\rm start}^{\rm T}(\sfC))\vee 0\big){\rm Drift}_{\rm vert}\big|\\
		\leq& 
		\big|{\rm T}_{t\wedge \tau_{\rm end}^{\rm T}} (\sfC)  -	{\rm T}_{t\wedge \tau_{\rm end}^{\rm T}}^{\rm IIC}(\sfC)\big|  \\
		&+\big|{\rm T}_{t\wedge \tau_{\rm end}^{\rm T}}^{\rm IIC}(\sfC) - {\rm T}_0(\sfC) - \big((t \wedge \tau_{\rm end}^{\rm T} (\sfC) - \tau_{\rm start}^{\rm T}(\sfC))\vee 0\big){\rm Drift}_{\rm vert}\big|  \\
		&+\big|(\tau_{\rm end}^{\rm T}(\sfC) -  \tau_{\rm end}^{{\rm Drift},{\rm T}}(\sfC)){\rm Drift}_{\rm vert}\big|.
	\end{align*}
	Each of the three terms in the right-hand side is~$O(\eta^c N)$ due to~${\rm TotErr} \leq C\,\eta^c N$, the occurrence of~${\rm GoodRW}$ and~\eqref{eq:tau_Drift_RW}, respectively. 

	The same computation applies to~${\rm B}_t(\sfC)$. For~${\rm L}_t(\sfC)$ and~${\rm R}_t(\sfC)$, 
	the only difference is the additional error due to the 
	uncertainty about whether the extremum is active or not (see \eqref{eq:Delta_err_def} and \eqref{eq:Ext_t}).
	Such errors only occur for time steps between~$\tau_{\rm start}^{\rm T} (\sfC)$ and~$\tau_{\rm start}^{\rm B} (\sfC)$ 
	and between~$\tau_{\rm end}^{\rm B} (\sfC)$ and~$\tau_{\rm end}^{\rm T} (\sfC)$. 
	As such, for $A \in \{{\rm L},{\rm R}\}$, \eqref{eq:TotErr} should be replaced by
	\begin{align}\label{eq:TotErrLR}
		\big|{A}_t(\sfC) -  {A}_t^{\rm IIC}(\sfC)\big|\leq C \Big(  {\rm TotErr}  + 
		\tau_{\rm start}^{\rm B} (\sfC)-\tau_{\rm start}^{\rm T} (\sfC) +  \tau_{\rm end}^{\rm T} (\sfC) - \tau_{\rm end}^{\rm B} (\sfC)
		\Big),
	\end{align}
	for~$0 \leq t \leq \tau_{\rm end}^A(\sfC)$.
	%However, the number of such transformations is bounded due to the upper bound on~${\rm T}_t (\sfC) -{\rm B}_t (\sfC)$ 
	The right-hand side of the above may be controlled by~\eqref{eq:meso+survive2} in addition to the assumed bound on ${\rm TotErr}$. 
	We conclude that 
	\begin{align*}
		\big|{A}_t(\sfC) - {A}_0(\sfC)    -  \big((t \wedge \tau_{\rm end}^{{\rm Drift}, {\rm T}} (\sfC) - \tau_{\rm start}^{\rm T}(\sfC))\vee 0\big)\cdot {\rm Drift}_*\big|
		=O(\eta^c N)	+ O(\sqrt \eta N)
	\end{align*}
	for all~$\sfC \in {\rm Meso}^+(\omega_0)$ and~$A \in \{{\rm T},{\rm B},{\rm L},{\rm R}\}$, where~${\rm Drift}_*$ refers to the drift in the direction associated to~$A$. 
	
	Thus, when~${\rm TotErr} \leq C\,\eta^cN$ and~${\rm GoodRW}$ occurs, all clusters~$\sfC \in{\rm Meso}^+(\omega_0)$ are stable in the sense of~\eqref{eq:stable_C}, provided that~$c,C >0$ are chosen appropriately. Finally,~\eqref{eq:TotErr2} and~\eqref{eq:GoodRW} ensure that the above events occur with probability at least~$1 - C \eta^c$, provided $\eta$ is chosen in accordance with Proposition~\ref{prop:error_exists}. This implies the desired conclusion.  
\end{proof}

\subsection{Homotopy classes}\label{sec:homotopy}

Fix some~$\delta > 0$ and~$N \geq 1$ and a lattice of the type~$\bbL(\pmb\alpha)$ with $\pmb\alpha \in \{\alpha,\beta\}^\bbZ$.
For a configuration~$\omega$,
a {\em nail} at scales~$(\delta,N)$ for a point~$x \in \delta^{1/2}N \bbZ^2$ is any cluster~$\sfC$ of~$\omega$ such that 
\begin{align}
	|{\rm T}(\sfC)-\langle x,e_{\rm vert}\rangle  | &\leq \delta N, \\
	|{\rm B}(\sfC)-\langle x,e_{\rm vert}\rangle  | &\leq \delta N,\\
	|{\rm L}(\sfC)-\langle x,e_{\rm lat}\rangle  | &\leq \delta N,\\
	|{\rm R}(\sfC)-\langle x,e_{\rm lat}\rangle  | &\leq \delta N.\label{eq:nails}
\end{align}
% We will identify the nails and the loops surrounding them in the loop representation of~$\omega$. 

\begin{figure}
\begin{center}
\includegraphics[width = 0.7\textwidth]{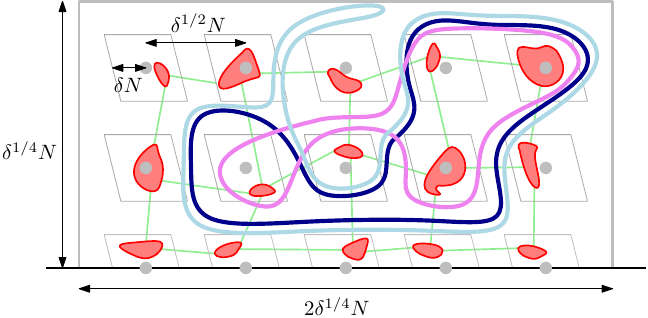}
\caption{A lattice of nails (in red) is formed of one nail for each point~$x$ in the grey lattice, within the rectangular window. The grey rhombi designate the regions in which nails should be contained, according to~\eqref{eq:nails}. The blue, light-blue and pink loops are all non-trivial. They all surround the same nails, but the homotopy class of the pink loop is different from that of the two blue loops: the homotopy class is given by the order in which the loops intersect the light green segments. 
The two blue loops are at a large~$d$-distance  of each other in the sense of~\eqref{eq:loop_dist} in spite of having the same homotopy class. This is because the light blue loop has a long arm wiggling between nails, without surrounding any of them; such a behaviour is unlikely due to~\eqref{eq:RSW_iso}.}
\label{fig:lattice_nails}
\end{center}
\end{figure}

A~$(\delta,N)$-lattice of nails is the choice of a nail~$\sfN(x)$ for each point~$x \in (\delta^{1/2}N )\cdot (\bbZ^2 \cap [-\delta^{-1/4}, \delta^{-1/4}] \times [0, \delta^{-1/4}] )$. Note that the lattice of nails roughly occupies the window $[-\delta^{1/4} N ,\delta^{1/4} N] \times [0,\delta^{1/4} N]$,  rather than the full observation window~$[-N,N] \times [0,N]$. This is for technical reasons that will be apparent below. 
A configuration may contain no lattice of nails or several such lattices. See Figure~\ref{fig:lattice_nails}.

For a fixed lattice of nails, a loop~$\ell$ in the loop representation of~$\omega$ is considered non-trivial if it surrounds at least two, but not all of the nails of the lattice. 
Consider the punctured plane~$\bbR^2 \setminus (\bigcup_x \sfN(x))$, where we mean that we remove the interior of the outer contour of each nail. 
Notice that any non-trivial loop is contained in~$\bbR^2 \setminus (\bigcup_x \sfN(x))$ and has a certain homotopy class. 

It will be important below to have a standard representation of homotopy classes in such a punctured plane, as we will compare homotopy classes with respect to different families of nails. We do this as follows. 
Write~$\vec E$ for the oriented edges of the  graph~$\bbZ^2 \cap [-\delta^{-1/4}, \delta^{-1/4}] \times [0, \delta^{-1/4}]~$ (each unoriented edge of the graph corresponds to two oriented edges). 
Fix a point~$\sfn_x$ in each nail~$\sfN(x)$ of the nail corresponding to~$x$. Identify~$\vec e = (xy)$ with the oriented segment between~$\sfn_x$ and~$\sfn_y$.

Let~$\calW$ be the set of finite words on the alphabet~$\vec E$ and denote the empty word by~$\emptyset$. Define the equivalence relation~$\sim$ on~$\calW$ generated by~$(u_i)_{1\le i\le p} \sim (v_j)_{1\le j\le q}$ if 
\begin{itemize}
	\item~$p=q$ and there exists~$k\in[1,p]$ such that~$u_1\ldots u_p = v_k \ldots v_p v_1\ldots v_{k-1}$ or 
	\item~$p=q+2$,~$u_1\ldots u_q = v_1 \ldots v_q$ and~$u_{q+2}$ is the same as the edge~$u_{q+1}$ but with the opposite orientation.
\end{itemize}
Define the set of {\em reduced words} as the quotient~${\calC\calW} := \calW / \sim$. 

Recall that the loops of~$\omega$ are oriented so as to have primal edges on their right. 
For a loop $\ell$, write~$w_0(\ell) \in \calW$ for the word formed of the sequence of edges~$\vec e$ crossed by~$\ell$ from left to right, 
and~$w(\ell) = w_\sfN(\ell)$ for the reduced word corresponding to~$w_0(\ell)$.

It is standard to check that this indeed encodes the homotopy class of every non-trivial loop and that~$w_\sfN(\ell)$ does not depend on the choice of points~$\sfn_x$. It may depend on the choice of the lattice of nails; however we will eventually prove that this only occurs with low probability (see Proposition~\ref{prop:homotopy_to_CN}).

\begin{definition}%[Homotopy distance]
We say that two percolation configurations~$\omega$ and~$\omega'$ are {\em homotopically similar} at scales~$(\delta,N)$ 
if they both contain lattices of nails~$\sfN$ and~$\sfN'$ at scales~$(\delta,N)$ such that
there exists a bijection~$\psi$ between the non-trivial loops of~$\omega$ 
with respect to~$\sfN$ and non-trivial loops of~$\omega'$ with respect to~$\sfN'$ with 
\begin{align}
	w_\sfN(\ell) = w_{\sfN'}(\psi(\ell)),
\end{align}
for all non-trivial loops~$\ell$ of~$\omega$. 
\end{definition}

The same also applies to piece-wise linear deformations of percolation configurations. 

\begin{proposition}\label{prop:homotopy_to_CN}
	%Let~$\bbL = \bbL(\beta)$ and~$\bbL' = \bbL(\alpha)$ be two isoradial rectangular lattices with constant angles 
	Let $\alpha, \beta \in (0,\pi)$ and~$M$ be an invertible linear transformation of~$\bbR^2$. 
	There exist constants~$C,c >0$ such that the following holds for all~$\delta > 0$ and~$N \geq \delta^{-C}$. 
	If~$\bbP$ is a coupling of~$\phi_{\bbL(\beta)}$ and~$\phi_{\bbL(\alpha)}$ such that 
	\begin{align*}
		\bbP\big[\text{$M(\omega)$ and~$\omega'$ are homotopically similar at scales~$(\delta,N)$}] \geq 1- \delta,
	\end{align*}
	then
	\begin{align}\label{eq:homotopy_to_CN}
		{d}_{\rm CN}\big[ \phi_{\delta^c\bbL(\beta)}\circ M^{-1}, \phi_{\delta^c\bbL(\alpha)}\big] \leq C\,\delta^c.
	\end{align}
\end{proposition}

The rest of the section is dedicated to the proof of the above. The idea is explained in Figures~\ref{fig:lattice_nails} and~\ref{fig:nails2}.

We start off with a helpful consequence of the~\eqref{eq:RSW_iso} property. Assume~$\delta$ and~$N$ fixed. 
For a point~$x \in \sqrt\delta N \bbZ^2$, let~$\sfD(x)$ be the rhombus-shaped region of~$\bbR^2$ defined as 
\begin{align}
	\sfD(x) = \big\{ y \in \bbR^2:\, 
	|\langle y - x,e_{\rm vert}\rangle  |\leq \delta N \text{ and }
	|\langle y - x,e_{\rm lat}\rangle  | \leq \delta N \big\}.
\end{align}
These are the regions containing the nails at scales~$(\delta, N)$.
A rectangle~$[0,4R] \times [0,R]$ is said to contain a {\em thin crossing} if it contains a primal or dual cluster~$\sfC$ that crosses the rectangle vertically, but surrounds none of the regions~$\sfD(x):x \in \sqrt\delta N \bbZ^2$. 
The same notion applies to crossings  of~$[0,R] \times [0,4R]$ in the horizontal direction. 

\begin{lemma}\label{lem:homotopy_to_CN}
	There exist constants~$c_0,C_0 > 0$ such that for every $\alpha\in(0,\pi)$, 
	any~$\delta > 0$,~$N\geq 1$, any~$R \geq \sqrt \delta N$, 
	\begin{align}
		\phi_{\bbL(\alpha)}\big[ \text{$[0,4R] \times [0,R]$ contains a thin crossing}\big]\leq C_0 (\sqrt \delta N / R)^{c_0}.
	\end{align}
\end{lemma}

\begin{proof}[Proof of Lemma~\ref{lem:homotopy_to_CN}]
	We only sketch this. Explore successive vertical crossings of~$[0,4R] \times [0,R]$ from left to right. 
	Whenever such a crossing~$\Gamma$ is revealed, the region to its right is unexplored. We can then identify~$R/\sqrt \delta N$ regions~$\sfD(x)$  
	that are within a distance~$2\sqrt \delta N$ of~$\Gamma$, but lie entirely in the unexplored region; see Figure~\ref{fig:surrounding_nail}.
	Then,~\eqref{eq:RSW_iso} states that each such region is surrounded by the cluster of~$\Gamma$ with probability at least a uniformly positive constant $c>0$. Furthermore, the events that  these regions are surrounded may be stochastically dominated by independent Bernoulli variables of uniformly positive parameter $c'=c'(c)>0$. 
	
	\begin{figure}
	\begin{center}
	\includegraphics[width = 0.65\textwidth]{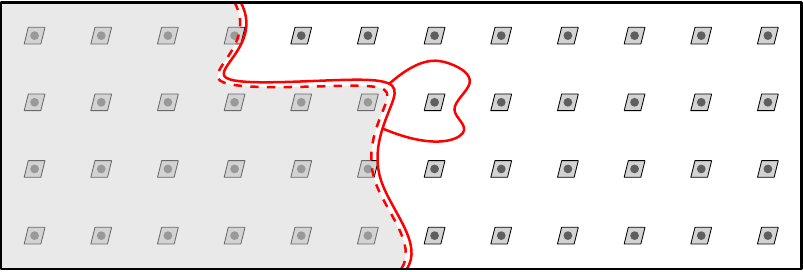}
	\caption{The left boundary~$\Gamma$ of a vertically crossing primal cluster; the conditioning only affects the grey area. 
	Under this conditioning, the cluster of~$\Gamma$ surrounds each one of the grey rhombi~$\sfD(x)$ closest to it with uniformly positive probability. }
	\label{fig:surrounding_nail}
	\end{center}
	\end{figure}
	
	We conclude that, for any such~$\Gamma$, the probability that the cluster of~$\Gamma$ forms a thin crossing is bounded above by~$(1-c')^{R/\sqrt \delta N}$. 
	Combining this with an exponential tail on the number of crossings provides the desired bound. 
\end{proof}

\begin{proof}[Proof of Proposition~\ref{prop:homotopy_to_CN}]
	%We only sketch this proof; a more detailed proof may be found in~\cite{DumKozKra20}. 
	For simplicity, we consider first the case where~$M = {\rm id}$. 
	For a configuration~$\omega$ and a lattice of nails~$\sfN$, we may define a reference configuration~$\omega_{\rm ref}$ 
	depending only on the homotopy information of~$\omega$ and which is homotopically similar to~$\omega$ at scales~$(\delta,N)$. The configuration $\omega_{\rm ref}$ is not actually a percolation configuration, but only a pair of families of loops which might not correspond to a percolation configuration.
	It is constructed as follows; see also Figure~\ref{fig:nails2}.
	
	The nails~$\sfN_{\rm ref}$ are single primal vertices, namely those closest to the points of~$(\delta^{1/2}N )\cdot (\bbZ^2 \cap [-\delta^{-1/4}, \delta^{-1/4}] \times [0, \delta^{-1/4}] )$.
	For each non-trivial loop~$\ell$ of~$\omega$, place a loop in~$\ell_{\rm ref} \in \omega_{\rm ref}$ of the same type as $\ell$, 
	formed of straight line segments between the midpoints of the edges of~$(\delta^{1/2}N )\cdot (\bbZ^2 \cap[-\delta^{-1/4}, \delta^{-1/4}] \times [0, \delta^{-1/4}])$,
	so that 
	\begin{align}
	w_\sfN(\ell) = w_{\sfN_{\rm ref}}(\ell_{\rm ref}).
	\end{align}
	
	We will prove that, with high probability, when~$\omega \sim \phi_{\bbL(\alpha)}$, 
	 \begin{align}\label{eq:omegaomega_ref}
		{d}_{\rm CN}\big[\tfrac{1}{\delta^{3/8} N} \omega, \tfrac{1}{\delta^{3/8} N} \omega_{\rm ref}\big] \leq \delta^{1/16}.
	\end{align}
	Applying\footnote{We actually bound the probability of~\eqref{eq:omegaomega_ref} for the configurations shifted vertically by~$-\frac12\delta^{1/4}N$, but this has no influence on our conclusion, as the measures are translationally invariant.} this to~$\omega$ and~$\omega'$ directly implies~\eqref{eq:homotopy_to_CN}, since~$\omega_{\rm ref} =\omega'_{\rm ref}$. 
	
	% In~$\tfrac{1}{\delta^{3/8} N} \omega$, nails are separated by a distance~$\delta^{1/8}$ and cover a~$2\delta^{-1/8}\times \delta^{-1/8}$ rectangular window. 
	Assume that~\eqref{eq:omegaomega_ref} fails.
	Then, there exists a loop~$\ell \in \omega$ contained in~$[-\delta^{1/4}N,\delta^{1/4}N] \times [0,\delta^{1/4}N]~$ 
	so that 
	\begin{align}\label{eq:long_arm}
		d\big(\tfrac{1}{\delta^{3/8} N}\ell, \tfrac{1}{\delta^{3/8} N}\ell_{\rm ref}\big) \geq \delta^{1/16}.
	\end{align}
	This requires~$\ell$ to contain a long ``arm'' of diameter~$\delta^{7/16} N$ that avoids all nails of~$\sfN$; 
	see Figure~\ref{fig:nails2}.
	
	In particular, if we take~$R =\frac12\delta^{7/16} N$ and if we pave the rectangle~$[-\delta^{1/4}N,\delta^{1/4}N]\times [0,\delta^{1/4}N]$ with the translates of~$[0,4R] \times [0,R]$ and~$[0,R] \times [0,4R]$ by points of~$R\bbZ^2$, 
	then~\eqref{eq:long_arm} implies that there exists at least one such rectangle that contains a thin crossing.
	Since there are~$O(\delta^{-3/8})$ such rectangles, Lemma~\ref{lem:homotopy_to_CN} and the union bound imply that 
	\begin{align}\label{eq:homotopy_to_CN0}
		\phi_{\bbL(\alpha)} \Big[{d}_{\rm CN}\big(\tfrac{1}{\delta^{3/8} N} \omega, \tfrac{1}{\delta^{3/8} N} \omega_{\rm ref}\big) > \delta^{1/16}\Big]
		\leq C_1  \delta^{-3/8}	 \exp(-c_0 \delta^{-1/16})
		\leq C_2\delta^{c_2},
	\end{align}
	where the constants~$C_1,C_2,c_2 > 0$ are universal. 
	This implies~\eqref{eq:homotopy_to_CN}, as explained above.
	
	\begin{figure}
	\begin{center}
	\includegraphics[width = .8\textwidth, page = 2]{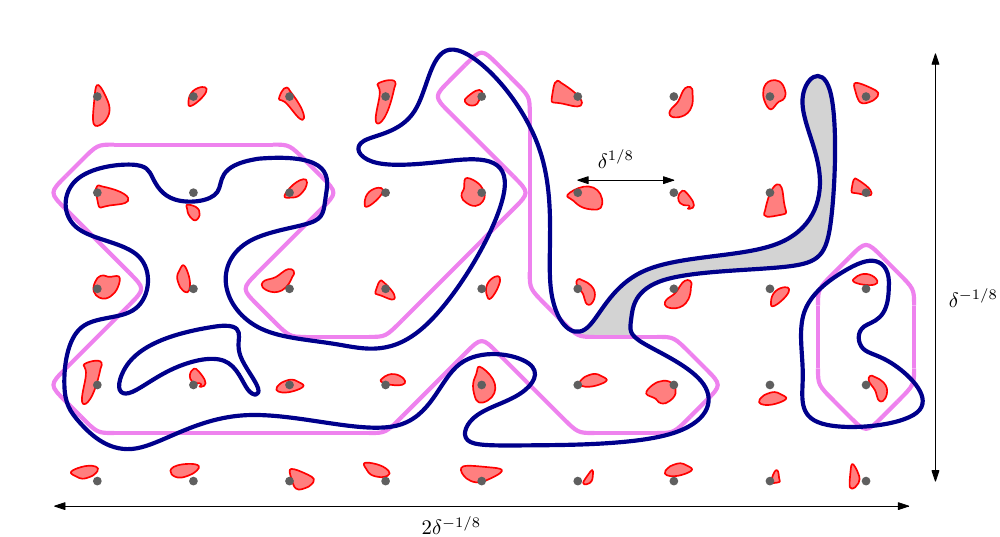}
	\caption{A configuration containing a lattice of nails~$\sfN$ in red, and several blue loops in a rescaled configuration~$\tfrac{1}{\delta^{3/8} N} \omega$. Of the blue loops, one is trivial, the other two have corresponding pink loops in~$\omega_{\rm ref}$. The large blue loop~$\ell$ is not close to~$\ell_{\rm ref}$ because of the long grey ``arm''. Here we assume that this arm has length at least~$\delta^{1/16}$. The probability of such an arm occurring may be bounded by~\eqref{eq:RSW_iso}.}
	\label{fig:nails2}
	\end{center}
	\end{figure}

	When~$M \neq {\rm id}$, we may still prove~\eqref{eq:homotopy_to_CN0} for both~$\omega$ and~$\omega'$. 
	Then, it suffices to observe that, since~$M$ is linear,~$M(\omega_{\rm ref}) = \omega_{\rm ref}'$.	
\end{proof}

\subsection{Proof of Theorem~\ref{thm:linear}}\label{sec:5linear_proof}

We are finally ready for the proof of Theorem~\ref{thm:linear}. 
Fix~$\alpha,\beta \in (0,\pi)$. 

In light of Proposition~\ref{prop:homotopy_to_CN}, our goal is to prove that for~$N$ large enough 
and~$\delta \geq N^{-c_\delta}$ for some small constant~$c_{\delta}>0$, 
we may couple~$\omega\sim \phi_{\bbL(\beta)}$ and~$\omega'\sim \phi_{\bbL(\alpha)}$ so that
$M_{\beta,\alpha}(\omega)$ and~$\omega'$ are homotopically similar at scales~$(\eta,N)$. 
The coupling will be that of Section~\ref{sec:5MarkovChain}, as explained below. 
Recall the quantity $\eta > N^{-c}$ appearing in Proposition~\ref{prop:stability_meso}
and in the definition of mesoscopic cluster. 
Fix $\delta > \sqrt \eta$. All definitions below depend on $\eta$ and $\delta$, which will be chosen below as small negative powers of $N$. 

We call a cluster~$\sfC$ of~$\omega$ a {\em nail$^+$} at~$x \in (\delta^{1/2}N )\cdot (\bbZ^2 \cap [-\delta^{-1/4}, \delta^{-1/4}] \times [0, \delta^{-1/4}])$
if 
\begin{align*}
	|{\rm T}(\sfC)-\langle x,e_{\rm vert}\rangle  | &\leq \tfrac{\delta}{2C_{M}} N, \\
	|{\rm B}(\sfC)-\langle x,e_{\rm vert}\rangle  | &\leq\tfrac{\delta}{2C_{M}} N,\\
	|{\rm L}(\sfC)-\langle x,e_{\rm lat}\rangle  | &\leq \tfrac{\delta}{2C_{M}} N,\\
	|{\rm R}(\sfC)-\langle x,e_{\rm lat}\rangle  | &\leq \tfrac{\delta}{2C_{M}} N.
\end{align*}
with~$C_M$ chosen as in the definition of~${\rm Meso}^+$.
We consider here the more restricted notion of nail$^+$ rather than just nail, 
for the same reason that we considered~${\rm Meso}^+(\omega_0)$ rather than~${\rm Meso}(\omega_0)$ in Proposition~\ref{prop:stability_meso}:
a buffer is needed to offset the small variations that may occur during the process~$(\omega_t)_{0 \leq t \leq K+K'}$.

\begin{lemma}\label{lem:lattice_nails_exists}
	There exist constants~$C,c >0$ such that, for all~$\delta,\eta >0$ with~$\delta$ small enough and $\eta <\delta^2$, and any~$N$ large enough,
	\begin{align}\label{eq:existence_nails}
	\phi_{\bbL(\beta)}[\text{there exists a lattice of nails$^+$ in~${\rm Meso}^+$ at scales~$(\delta,N)$}] \ge 1 - C\eta^{c}  \delta^{-1/2} .
	\end{align}
\end{lemma}

\begin{proof}
	This is a direct consequence of~\eqref{eq:RSW_iso}. 
	Observe that, for given~$\delta$, there are at most~$O(\delta^{-1/2})$ points in the lattice~$(\delta^{1/2}N )\cdot (\bbZ^2 \cap [-\delta^{-1/4}, \delta^{-1/4}] \times [0, \delta^{-1/4}])$. 
	Furthermore, for~$\delta$ small enough, any such point~$x$ satisfies
	\begin{align}
		\langle x,e_{\rm vert}\rangle \leq \tfrac{N}{C_{M}} \quad \text{ and }\quad 
		- \tfrac{N}{C_{M}}\leq	\langle x,e_{\rm lat}\rangle \leq \tfrac{N}{C_{M}}, 
	\end{align}
	which is to say that it is contained in the window covered by the definition of~${\rm Meso}^+$.
	Finally, for each such point~$x$, since~$\sqrt \eta \leq \delta$,~\eqref{eq:RSW_iso} implies that 
	\begin{align}
		\phi_{\bbL(\beta)}[ \text{there exists a nail$^+$~$\sfN(x) \in {\rm Meso}^+$ at~$x$}] \geq 1 - C_0\, \eta^{c_0}. 
	\end{align}
	for universal constants~$c_0,C_0$. Performing a union bound yields the result. 
\end{proof}

\begin{proof}[Proof of Theorem~\ref{thm:linear}]
%We are finally ready to choose~$\eta$ and~$\delta$. 
For~$N\geq 1$, write~$\eta = N^{-c_\eta}$ and~$\delta = N^{-c_\eta^2}$, for some small~$c_\eta >0$ to be determined below. 
In particular, we assume~$c_\eta < 1/2$ so that Lemma~\ref{lem:lattice_nails_exists} applies. 

Recall the definition of ${\rm Stability}(\eta,N)$ for the event in Proposition~\ref{prop:stability_meso}, 
with the constants~$c,C$ provided by the proposition. 
By choosing~$c_\eta > 0$ small enough, Lemma~\ref{lem:lattice_nails_exists} and Proposition~\ref{prop:stability_meso} ensure that 
\begin{align}
	\bbP[\text{$\omega_0$ contains a lattice of nails$^+$ in~${\rm Meso}^+(\omega_0)$ at scales~$(\delta,N)$}] &\ge 1 - N^{-c}\,\,\text{ and }
	% Here~$c = \frac14 c_0 c_\eta$, there~$c_0$ is the constant given in Lemma~\ref{lem:lattice_nails_exists}, and we assume that~$c_\eta \leq c_0/2$. 
	\nonumber\\
	\bbP[\text{Stability}(\eta,N)]& \ge 1 - N^{-c},
	\label{eq:stability+}
\end{align}
for all~$N$ large enough, where~$c$ is some constant which depends on $c_\eta$, but not on $N$.

Assume henceforth that both of the events above occur and that~$N$ is large. 
Fix a lattice of nails$^+$~$\sfN_0 \subset {\rm Meso}^+(\omega_0)$. 
By~$\text{Stability}(\eta,N)$, these clusters survive throughout the process~$(\omega_t)_{0\leq t\leq K+K'}$.
Denote by~$\sfN_t$ their collection in~$\omega_t$.
Furthermore, due to~$\text{Stability}(\eta,N)$ and the definition of  nail$^+$, 
$\sfN_t$ is a lattice of nails, up to a piecewise-linear transformation, for every~$t \leq K+K'$.
As such, it makes sense to encode the homotopy information of~$\omega_t$ with respect to~$\sfN_t$, in the same way as for~$\omega_0$. 

We will now argue that this homotopy information remains unchanged throughout the process. 
More precisely, any non-trivial loop~$\ell_0$ surrounding a primal or dual cluster has a corresponding loop~$\ell_t$ in~$\omega_t$ and 
$w_{\sfN_t}(\ell_t) = w_{\sfN_0}(\ell_0)$. 
Indeed, when passing from~$\omega_t$ to~$\omega_{t+1/2}$, as~$\ell_t$ has diameter larger than~$\eta N$, 
it either avoids all extremal boxes or~$\omega_t = \omega_{t+1/2}$. 
In both cases, we conclude that~$w_{\sfN_{t+1/2}}(\ell_{t+1/2}) = w_{\sfN_t}(\ell_t)$.
Furthermore, when passing from~$\omega_{t+1/2}$ to~$\omega_{t+1}$, the fact that track exchanges do not break loops ensures that the homotopy class of any non-trivial loop is preserved. As such~$w_{\sfN_{t+2}}(\ell_{t+2}) = w_{\sfN_{t+1}}(\ell_{t+1}) = w_{\sfN_{t}}(\ell_{t})$.

Finally, as explained in Remark~\ref{rem:M_origin}, the effect of the drift on the lattice~$\sfN_0$ is described by~$M_{\beta,\alpha}$.
More precisely, for all~$x\in (\delta^{1/2}N )\cdot (\bbZ^2 \cap [-\delta^{-1/4}, \delta^{-1/4}] \times [0, \delta^{-1/4}])$, 
\begin{align}
%	\big|\langle y, e_{\rm vert} \rangle - \langle M_{\beta,\alpha}(x), e_{\rm vert} \rangle \big|&\leq \tfrac{\delta}{2C_{M}} N + C N^{1-2c} \quad\text{ and } \\
%	\big|\langle y, e_{\rm lat} \rangle - \langle M_{\beta,\alpha}(x), e_{\rm lat} \rangle \big| &\leq  \tfrac{\delta}{2C_{M}} N + C N^{1-2c} \qquad 
	\big|\langle y, e_{*} \rangle - \langle M_{\beta,\alpha}(x), e_{*} \rangle \big| &\leq  \tfrac{\delta}{2C_{M}} N + C N^{1-c} \qquad 
	 \forall y \in \sfN_{K+K'}(x) \text{ and } * \in \{{\rm vert},{\rm lat}\}.
\end{align}
By again assuming that~$c_\eta$ is sufficiently small and~$N$ large, we conclude that 
$M_{\beta,\alpha}^{-1}(\sfN_{K+K'})$ is a~$(\delta,N)$-lattice of nails.
Then,~$\omega_0$ and~$M_{\beta,\alpha}^{-1}(\omega_{K+K'})$ are homotopically similar at scales~$(\delta,N)$. 
Applying Proposition~\ref{prop:same_law}, Proposition~\ref{prop:homotopy_to_CN}, and using~\eqref{eq:stability+}, we conclude that 
\begin{align}
	{d}_{\rm CN}\big[ \phi_{\delta\bbL(\alpha)}, \phi_{\delta\bbL(\beta)}\circ M_{\beta,\alpha} \big] \leq C\,\delta^c \leq C N^{-c'},
\end{align}		
for all~$N$ large enough, where~$C,c,c' >0$ are universal constants. The assumption on~$N$ may be removed by altering the constants. 
\end{proof}

\subsection{An equality of drifts}\label{sec:5equality_drifts}

As an addition to Theorem~\ref{thm:linear}, we make an apparently obvious but crucial observation about the value of the drift vectors. This observation will be useful in Section~\ref{sec:drift0}.

\begin{proposition}\label{prop:drift_RT}
	For any~$0  < \beta < \pi$,
	\begin{align}\label{eq:drift_RT}
		{\rm Drift}_{\rm lat}(\beta,\beta/2) = {\rm Drift}_{\rm vert}(\beta,\beta/2).
	\end{align}
\end{proposition}

The above is actually a surprisingly profound fact. It ultimately shows that the isoradial embedding of the lattices~$\bbL(\alpha)$ (and more generally of bi-periodic graphs) is the right embedding of the inhomogeneous random-cluster model to ensure universality in the sense of Theorem~\ref{thm:universalCNSS}. Indeed, except in this proposition, there is no use of the exact formula linking isoradial embeddings to the inhomogeneity of the random-cluster parameters. 
Behind~\eqref{eq:drift_RT} lies the fact that the line exchanges producing the drift are composed of star-triangle transformations, which do not depend on the direction of the tracks being exchanged. 

\begin{proof} 
	Fix~$0 < \alpha  < \beta < \pi$. 
	We will prove a more general statement, which will imply~\eqref{eq:drift_RT} when taking~$\alpha = \beta /2$. 
	Several scales will appear in the construction below: consider integers~$r \ll N \ll R$, with~$r = N^{1-c}$ for some small constant~$c$.
	
	Consider the finite isoradial graph~$\bbL$ containing a block of~$R$ horizontal tracks with transverse angle~$\beta$ crossed by~$2R$ vertical tracks of transverse angle~$0$. 
	We call this the~$\beta$-block. Position it so that~$0$ is the central point on the bottom part of this block. 
	Complete the block by~$r$ horizontal tracks of transverse angle~$\alpha$ that run along its top and right side. 
	See Figure~\ref{fig:vTR}.
	
	\begin{figure}
	\begin{center}
	\includegraphics[width = .98\textwidth]{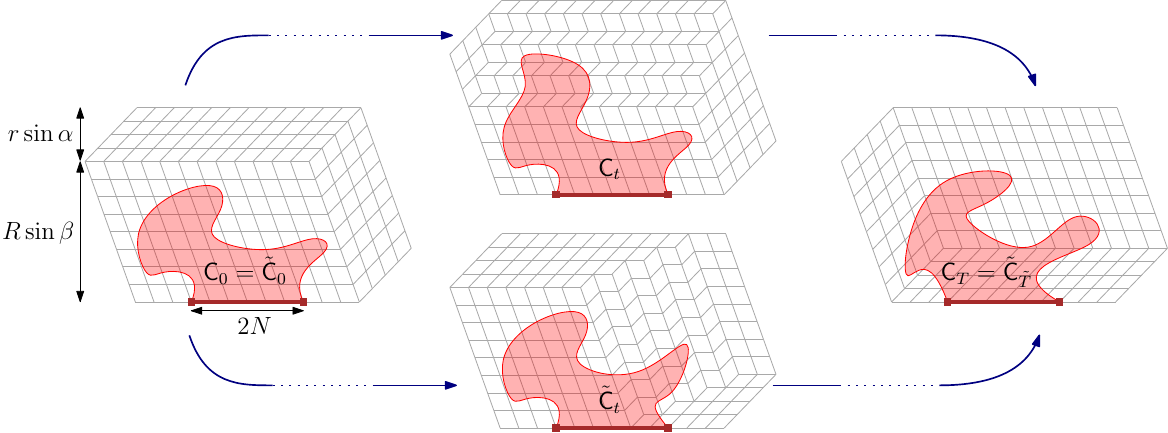}
	\caption{From~$\bbL$ on the left, there are two series of transformations leading to~$\bbL'$ on the right. 
	On the top, we proceed by exchanging horizontal tracks, bringing them down progressively. This produces a sequence of lattices~$\bbL_{t}$, and a sequence~$\sfC_t$ of transformations of the cluster of the bottom red segment. 
	On the bottom, we proceed by progressively exchanging vertical tracks. This leads to sequences~$\tilde \bbL_t$ and~$\tilde \sfC_t$.
	The results of the two processes are equal in law.}
	\label{fig:vTR}
	\end{center}
	\end{figure}

	Let~$\bbL'$ be the graph obtained by switching the position of the tracks of angle~$\alpha$ with the~$\beta$-block. In~$\bbL'$, the tracks of angle~$\alpha$ run on the left and bottom side of the~$\beta$-block. 
	
	There are two ways to pass from~$\bbL$ to~$\bbL'$ using compositions of star-triangle transformations. 
	As illustrated in Figure~\ref{fig:vTR}, we may define~$\bbL = \bbL_0,\dots, \bbL_{T} = \bbL'$, by essentially ``sliding'' the tracks of angle~$\alpha$ down through the~$\beta$-block. 
	To pass from~$\bbL_{t}$ to~$\bbL_{t+1}$ we perform a series of finite track exchanges, similarly to those in the transformation~$\bfS_t$ appearing in~\eqref{eq:5S_transformation}. 
	For~$r < t \leq  T-r$, the lattice~$\bbL_t$ contains a mixed block of~$2r$ tracks which slides down as~$t$ increases. 
	To pass from~$t$ to~$t+1$, exactly~$r$ track exchanges are performed, each consisting of~$2N$ sequential star-triangle transformations. 
	Let~$(\omega_t)_{0\leq t\leq T}$ be the Markov chain of configurations on the lattices~$\bbL_t$ defined as in Section~\ref{sec:5MarkovChain}, with~$\eta = N^{1-2c} \ll r$.

	A second way to  transform~$\bbL$ into~$\bbL'$ is to slide the tracks of angle~$\alpha$ from the right to the left; 
	write~$\tilde \bbL_0,\dots, \tilde \bbL_{\tilde T}$ for this set of transformations and~$(\tilde\omega_t)_{0 \leq t\leq \tilde T}$ for the corresponding chain of configurations. 
	
	The distortion of the contours of the large clusters in each configuration may be controlled in the same way as in Theorem~\ref{thm:linear}. 
	In particular, if we write~$\sfC_t$ and~$\tilde\sfC_{t}$ to be the cluster of the interval~$[-N,N] \times \{0\}$ in~$\omega_t$ and~$\tilde\omega_t$, respectively, 
	we find that 
	\begin{align}\label{eq:drift_vh}
	 	\lim_{R \to \infty} \bbE\big[{\rm R}(\sfC_{T}) - {\rm R}(\sfC_{0})\big] 
		&= %\frac{r\sin\beta}{\sin\beta + {\rm Drift}_{\rm vert}(\beta,\alpha)} 
		\frac{2r(\sin \alpha+ \sin\beta)}{ \sin \alpha+ \sin\beta + {\rm Drift}_{\rm vert}(\beta,\alpha)}
		\cdot {\rm Drift}_{\rm lat}(\beta,\alpha) + o(r)  \text{ and } \nonumber\\
	 	\lim_{R \to \infty} \bbE\big[{\rm R}(\tilde\sfC_{\tilde T}) - {\rm R}(\tilde\sfC_{0})\big]
%		& = % \frac{r\sin(\pi-\beta)}{\sin(\pi-\beta) + {\rm Drift}_{\rm vert}(\pi -\beta,\pi-\beta + \alpha)} 
%		\frac{2r(\sin (\pi-\beta+\alpha)+ \sin(\pi-\beta))}{ \sin (\pi-\beta+\alpha)+ \sin(\pi-\beta) + {\rm Drift}_{\rm vert}(\pi-\beta,\pi-\beta+\alpha)}
%		\cdot {\rm Drift}_{\rm vert}(\pi-\beta,\pi-\beta+\alpha)+ o(r) \nonumber\\
		& = %\frac{r\sin\beta}{\sin\beta + {\rm Drift}_{\rm vert}(\beta,\beta-\alpha)} 
		\frac{2r(\sin (\beta - \alpha)+ \sin\beta){\rm Drift}_{\rm vert}(\beta,\beta-\alpha)}{ \sin (\beta-\alpha)+ \sin\beta + {\rm Drift}_{\rm vert}(\beta,\beta-\alpha)}
		+ o(r).
	\end{align}
	In the above, the fractions are the approximate number of transformations that affect any mesoscopic cluster. 
	In the second line, the track-exchanges effectively occur between tracks of angle~$\pi-\beta$ and~$\alpha + \beta$ and the rightmost point of~$\tilde\sfC_t$ acts as a topmost point from the point of view of the track exchanges. We used the horizontal symmetry to state that~${\rm Drift}_{\rm vert}(\pi -\beta,\pi -\beta + \alpha) =  {\rm Drift}_{\rm vert}(\beta, \pi - \beta - \alpha)$; note however that~${\rm Drift}_{\rm lat}(\pi -\beta,\pi -\beta + \alpha) =- {\rm Drift}_{\rm lat}(\beta, \pi - \beta - \alpha)$, as the direction of the vector used for reference is reversed. 
	
	Finally, notice that~${\rm R}(\sfC_{0})$ and~${\rm R}(\tilde\sfC_{0})$ have the same law, as do~${\rm R}(\sfC_{T})$ and~${\rm R}(\tilde\sfC_{\tilde T})$. 
	Combining this with~\eqref{eq:drift_vh} and taking~$r$ and $N$ to infinity, we find
	\begin{align*}
		\frac{(\sin \alpha+ \sin\beta)\cdot {\rm Drift}_{\rm lat}(\beta,\alpha) }{ \sin \alpha+ \sin\beta + {\rm Drift}_{\rm vert}(\beta,\alpha)}	
	=
		\frac{(\sin (\beta - \alpha)+ \sin\beta)\cdot {\rm Drift}_{\rm vert}(\beta,\beta-\alpha)}{ \sin (\beta-\alpha)+ \sin\beta + {\rm Drift}_{\rm vert}(\beta,\beta-\alpha)}.
	\end{align*}
	Taking~$\alpha = \beta/2$ in the above, we conclude~\eqref{eq:drift_RT}.
\end{proof}

\subsection{Proof of Lemma~\ref{lem:M_cont}}\label{sec:lemM_cont}

%We now dive in the proof of Lemma~\ref{lem:M_cont}. 
We now focus on proving the continuity of $M_{\beta,\alpha}$ as a function of $\beta$ and $\alpha$.

\begin{proof}[Proof of Lemma~\ref{lem:M_cont}]
	We start with the proof of continuity. 
	We will focus here on the vertical drift; the proof is identical for the lateral drift. 
	Fix two angles~$\alpha, \beta$. 
	Let~$\bbL_{\rm mix}(\alpha,\beta)$ be the mixed lattice with angles~$\alpha,\beta$. 
	Also, define~$\bbL_{\rm mix}(\alpha,\beta) \cap \Lambda_R~$ for the restriction of this lattice to~$\Lambda_R$, 
	with the addition of the rhombi needed to perform~$\bfS_1 \circ \bfS_0$. 
	
	Fix~$\eps > 0$. 
	Due to Lemma~\ref{lem:IIC}, specifically to the fact that it is uniform in the angles~$\alpha$ and~$\beta$ 
	%(outside of a neighbourhood of~$0$ and~$\pi$) 
	and to the locality of the star-triangle transformation, 
	for~$R$ large enough
	we may couple the application of 
	$\bfS_1 \circ \bfS_0$ to~$\phi^{\rm IIC,T}_{\bbL_{\rm mix}(\alpha,\beta)}$ 
	and to~$\phi^1_{\bbL_{\rm mix}(\alpha,\beta) \cap \Lambda_R}[\cdot \,|\,{\rm Top}(\sfC_{x_n}) = 0 ]$ so that the increments 
	$\Delta_t^{\rm IIC} {\rm T}$ and~$\tilde \Delta_t^{\rm IIC} {\rm T}$ for the top of the IIC and of~$\sfC_{x_n}$, respectively,  satisfy
	\begin{align*}
		\bbP [\Delta_t^{\rm IIC} {\rm T} \neq \tilde \Delta_t^{\rm IIC} {\rm T}] \leq \eps. 
	\end{align*}
	Moreover, the choice of~$R$ may be uniform in a vicinity of~$(\alpha,\beta)$. 
	
	Since~$\Delta_t^{\rm IIC} {\rm T}$ is bounded, we conclude that 
	$(\alpha,\beta)\mapsto \phi^{\rm IIC,T}_{\bbL_{\rm mix}(\alpha,\beta)} [\Delta_t^{\rm IIC} {\rm T}]$ is 
	the locally uniform limit of 
	$(\alpha,\beta)\mapsto \phi^1_{\bbL_{\rm mix}(\alpha,\beta) \cap \Lambda_R}[\Delta_t^{\rm IIC} {\rm T}\,|\,{\rm Top}(\sfC_{x_n}) = 0 ]$.
	The latter is a continuous function, and therefore so is the former.

	Finally, we turn to the invertibility of~$M_{\beta,\alpha}$. By the explicit formula~\eqref{eq:Mbetaalpha}, 
	this is equivalent to 
	\begin{align}
		{\rm Drift}_{\rm vert}  \neq -\sin\beta.
	\end{align}
	As explained in Section~\ref{sec:stability_meso}, the finite energy property and a deterministic bound on the increments of the IIC yield~${\rm Drift}_{\rm vert}  > -\sin\beta$.
%	The increments~$\Delta_t^{\rm IIC} {\rm T}$ are deterministically bounded by~$-\sin \beta$,
%	and the finite energy property is enough to deduce that~${\rm Drift}_{\rm vert}  > -\sin\beta$.
\end{proof}

\subsection{Uniformity in angles}\label{sec:unif_angles}

Here, we give a version of Theorem~\ref{thm:linear} that is (partially) uniform in the choice of angles. 
We do not claim it to be optimal. 

\begin{proposition}\label{prop:linear_unif}
	Fix~$q\in[1,4]$,~$\beta \in (0,\pi)$ and assume $M_{\beta,\alpha} = {\rm id}$ for all $\alpha \in (0,\pi)$. 
	There exist constants~$c,C > 0$ such that, for any $\alpha \in (0,\pi)$, 
	\begin{align}\label{eq:linear_unif}
		{d}_{\rm CN}\big[ \phi_{\delta\bbL(\beta)}, \phi_{\delta\bbL(\alpha)}\big] \leq C\,\delta^c 
		\quad \text{ for all~$\delta > 0$}.
	\end{align}
\end{proposition}

We will not provide a full proof of the above, but simply explain which steps of the proof of Theorem~\ref{thm:linear} need to be adapted. 
The assumption on $M_{\beta,\alpha}$ being the identity is so as to remove any complications coming from the drift. 

It is obvious that  Theorem~\ref{thm:linear} is uniform in $\alpha$ and $\beta$ away from $0$ and $\pi$. We discuss below the uniformity when $\alpha$ approaches~$0$.

In the process of Theorem~\ref{thm:linear}, the number of steps diverges as $\alpha \to 0$ as $O(N/\alpha)$.
However, the probability of an error occurring at any given time step is of the order of $\epsilon(N) \alpha N$, where $\epsilon(N)$ is uniform in $\alpha$, and tends to zero as a power of $N$ when $N$ tends to infinity. Thus, when summing the probabilities of error, we conclude that \eqref{eq:TotErr2} remains valid uniformly in $\alpha$.

\section{Rotational invariance and universality: proofs of Theorems~\ref{thm:rotation_invariance_infinite_vol} and~\ref{thm:universalCNSS}}\label{sec:conclusion_rot_inv}

With Theorem~\ref{thm:linear} now proved, we turn to our two main results: 
the asymptotic rotational invariance of~$\phi_{\bbL(\pi/2)}$ (Theorem~\ref{thm:rotation_invariance_infinite_vol}) 
and the universality among lattices~$\bbL(\alpha)$ with~$\alpha \in (0,\pi)$  (Theorem~\ref{thm:universalCNSS}).
In light of Theorem~\ref{thm:linear} (or rather its uniform version Proposition~\ref{prop:linear_unif}), the latter is equivalent to~${\rm Drift}_{\rm lat}(\pi/2,\alpha)  = {\rm Drift}_{\rm vert}(\pi/2,\alpha) = 0$ for all~$\alpha \in (0,\pi)$. 

We start by proving in Section~\ref{sec:deducing_rot_inv} the asymptotic rotational invariance of~$\phi_{\bbL(\pi/2)}$ without using Theorem~\ref{thm:universalCNSS}. 
In Section~\ref{sec:drift0}, we use the asymptotic rotational invariance of~$\phi_{\bbL(\pi/2)}$ to deduce that~${\rm Drift}_{\rm lat}(\pi/2,\alpha)  = {\rm Drift}_{\rm vert}(\pi/2,\alpha) = 0$ for a dense set of angles~$\alpha \in (0,\pi)$, 
which in turn implies Theorem~\ref{thm:universalCNSS}.

%Finally, Section~\ref{sec:BA} contains a sketch of the direct computation of the drift using estimates of eigenvalues of the six-vertex models obtained via the Bethe-ansatz. This section is not necessary for the proof of any of the results here. it illustrates a different view on the interplay between inhomogeneous models and their isoradial embedding. 

\subsection{Rotational invariance: proof of Theorem~\ref{thm:rotation_invariance_infinite_vol}}\label{sec:deducing_rot_inv}

Write~$S_{\frac\alpha2}$ for the orthogonal reflection with respect to~$e^{{i}\alpha/2}$.
We start off with a result about the invariances of~$\phi_{\delta\bbL(\beta)}$ for general angles~$\beta$. 
Recall that~$\phi_{\bbL(\alpha)}$ is invariant under~$S_{\frac\alpha2}$, and since 
$\phi_{\bbL(\beta)}$ and~$\phi_{\bbL(\alpha)}$ are related via Theorem~\ref{thm:linear}, it follows that 
$\phi_{\bbL(\beta)}$ is asymptotically invariant under~$S_{\frac\alpha2}$ conjugated with~$M_{\beta,\alpha}$.

\begin{proposition}\label{prop:Tinv}
	For any~$\alpha, \beta \in (0,\pi)$, there exist constants~$c,C > 0$ such that
	\begin{align}\label{eq:autoT}
		{d}_{\rm CN}\big[\phi_{\delta\bbL(\beta)},\, \phi_{\delta\bbL(\beta)}\circ M_{\beta,\alpha}^{-1}\circ S_{\alpha/2}\circ M_{\beta,\alpha}  \big] 	 
		\leq C\, \delta^c \qquad \text{for all~$\delta > 0$}.
	\end{align}
\end{proposition}

The proposition is an immediate consequence of the argument described above. We will not give a detailed proof.

We now turn to the proof of Theorem~\ref{thm:rotation_invariance_infinite_vol}. For~$\alpha \in (0,\pi)$ write  
\begin{align*}
	T_\alpha = M_{\frac{\pi}2,\alpha}^{-1}\circ S_{\alpha/2}\circ M_{\frac{\pi}2,\alpha}.
\end{align*}
In light of~\eqref{eq:autoT}, we will say that~$\phi_{\bbL(\frac{\pi}2)}$ is {\em asymptotically invariant} under~$T_\alpha$.
Since this is the case for all~$\alpha \in (0,\pi)$, we conclude that~$\phi_{\bbL(\frac{\pi}2)}$ is asymptotically invariant 
with respect to the group generated by~$\{T_\alpha \,:\, \alpha \in (0,\pi)\}$.
To prove the asymptotic rotation invariance of~$\phi_{\bbL(\frac{\pi}2)}$, we will show that the group generated by~$\{T_\alpha \,:\, \alpha \in (0,\pi)\}$ contains all rotations.
To start, we list some properties of~$T_{\alpha}$.

\begin{proposition}\label{prop:T} We have the following properties:
	\begin{itemize}
		\item[(i)] For each~$\alpha \in (0,\pi)$,~$T_\alpha$ has eigenvalues~$1$ and~$-1$
		with eigenvectors~$M_{\frac{\pi}2,\alpha}^{-1}  e^{{i}\alpha/2}$ and~$M_{\frac{\pi}2,\alpha}^{-1}  e^{{i}(\alpha/2 + \pi/2)}$, respectively;
		\item[(ii)]~$\alpha \mapsto T_\alpha$ is continuous over~$(0,\pi)$;
		\item[(iii)]~$\alpha \mapsto T_\alpha$ is not constant.
	\end{itemize}
\end{proposition}

\begin{proof}
	\noindent(i) It suffices to observe that~$M_{\frac{\pi}2,\alpha}$ acts as a change of basis, and maps the vectors 
		$M_{\frac{\pi}2,\alpha}^{-1} e^{{i}\alpha/2}$ and~$M_{\frac{\pi}2,\alpha}^{-1}  e^{{i}(\alpha/2 + \pi/2)}$ onto the eigenvectors of~$S_{\alpha/2}$. \smallskip

		\noindent(ii) This is a direct consequence of the continuity of~$\alpha \mapsto M_{\frac{\pi}2,\alpha}$, which was proved in Lemma~\ref{lem:M_cont}. \smallskip

		\noindent(iii) 
		Let us proceed by contradiction and assume that~$T_\alpha = T_{\pi/2} = S_{\pi/4}$ for all~$\alpha$. 
		Under the assumption that~$T_\alpha = T_{\pi/2} = S_{\pi/4}$ and due to point (i),~$M_{\frac{\pi}2,\alpha}$ maps~$e^{{i}\pi/4}$  onto a multiple of~$e^{{i}\alpha/2}$ 
		and~$e^{{i}3\pi/4}$  onto a multiple of~$e^{{i}(\alpha/2 + \pi/2)}$.
		Moreover, recall from~\eqref{eq:Mbetaalpha} that~$M_{\frac{\pi}2,\alpha}$ acts as the identity on the horizontal axis. 
		This completely determines~$M_{\frac{\pi}2,\alpha}$: 
		\begin{align}\label{eq:degenerate_A}
			M_{\frac{\pi}2,\alpha} = \begin{pmatrix}1 & \cos \alpha \\0 & \sin \alpha \end{pmatrix}.
		\end{align}
		
		Write~${\rm Square}_\alpha$ for the square~$\{a e^{{i}\alpha/2} + b e^{{i}(\alpha/2 + \pi/2)}: 0\leq a,b\leq 1 \}$ 
		and~$\calC_h({\rm Square}_\alpha)$ for the event that the~${\rm Square}_\alpha$ contains a ``horizontal'' crossing, 
		that is a crossing between~$\{b e^{{i}(\alpha/2 + \pi/2)}: 0\leq b\leq 1 \}$ and~$\{e^{{i}\alpha/2} + b e^{{i}(\alpha/2 + \pi/2)}: 0\leq b\leq 1 \}$.
		The same notation applies to other rectangles. 
		Theorem~\ref{thm:linear} implies that 
		\begin{align}
			\lim_{\delta\rightarrow0}\big|\phi_{\delta \bbL(\alpha)}[\calC_h({\rm Square}_\alpha)] - \phi_{\delta \bbL(\pi/2)}[\calC_h(M_{\pi/2,\alpha}^{-1}{\rm Square}_\alpha)]\big| =0. 
		\end{align}
		Notice now that, due to~\eqref{eq:degenerate_A},~$M_{\pi/2,\alpha}^{-1}{\rm Square}_\alpha$ 
		is a rectangle of aspect ratio~$\tan(\alpha/2)$, which tends to~$0$ as~$\alpha$ tends to $0$.
		This implies that its crossing probability under~$\phi_{\delta \bbL(\pi/2)}$ may be made arbitrarily close to~$1$ (uniformly in $\delta$).
		Thus, we find that
		\begin{align*}
			\lim_{\alpha \to 0}\liminf_{\delta \to 0}\phi_{\delta \bbL(\alpha)}[\calC_h({\rm Square}_\alpha)]  = 1.
		\end{align*}
		This contradicts~\eqref{eq:RSW_iso}, particularly the fact that this property is uniform over~$\alpha$. 
%		Write~${\rm Square}_\alpha$ for the square~$\{a e^{{i}\alpha/2} + b e^{{i}(\alpha/2 + \pi/2)}: 0\leq a,b\leq 1 \}$, 
%		seen as a quad\footnote{Recall from Corollary~\ref{cor:symmetric_quad} the notion of quad and quad crossings.} with the four marked points in the corners of the square, starting with~$0$ and distributed in counter-clockwise order. 
%		Theorem~\ref{thm:linear} implies that 
%		\begin{align}
%			\big|\phi_{\delta \bbL(\alpha)}[\calC({\rm Square}_\alpha)] - \phi_{\delta \bbL(\pi/2)}[\calC(M_{\pi/2,\alpha}^{-1}{\rm Square}_\alpha)]\big| \to 0 \text{ as~$\delta \to 0$}. 
%		\end{align}
%		Notice now that, due to~\eqref{eq:degenerate_A},~$M_{\pi/2,\alpha}^{-1}{\rm Square}_\alpha$ 
%		is a rectangle of aspect ratio~$\tan(\alpha/2)$, which tends to~$0$ as~$\alpha \to 0$.
%		This implies that its crossing probability under~$\phi_{\delta \bbL(\pi/2)}$ may be made arbitrarily close to~$1$.
%		Thus, we find that
%		\begin{align*}
%			\lim_{\alpha \to 0}\liminf_{\delta \to 0}\phi_{\delta \bbL(\alpha)}[\calC({\rm Square}_\alpha)]  = 1.
%		\end{align*}
%		This contradicts the RSW property~\eqref{eq:RSW_iso}, particularly the fact that this property is uniform over~$\alpha$. 
\end{proof}

We are now in a position to prove the asymptotic rotational invariance of~$\phi_{\bbL(\frac{\pi}2)}$.

\begin{proof}[Proof of Theorem~\ref{thm:rotation_invariance_infinite_vol}]
	As stated in Proposition~\ref{prop:Tinv},~$\phi_{\bbL(\frac{\pi}2)}$ is asymptotically invariant under all~$T_\alpha$ with~$\alpha \in (0,\pi)$. 
	In addition, it is also invariant under the vertical reflection~$S_0$, or equivalently with respect to the rotation by~$\pi/2$. As opposed to~$S_{\pi/4}$,~$S_0$ is not part of~$\{T_\alpha\,:\,\alpha \in (0,\pi)\}$.
	As a consequence, it is asymptotically invariant under all transformation in the group generated by~$\{T_\alpha\,:\,\alpha \in (0,\pi)\}$ and~$S_0$. 
	We will prove that this group contains all rotations. 
		
	First, let us show that for all~$\alpha \in (0,\pi)$,~$T_\alpha$ is an {orthogonal} reflection. 
	Fix~$\alpha$ and let~$u,v \in \bbC$ be the eigenvectors of~$T_\alpha$ of Euclidean norm~$1$ and eigenvalues~$1$ and~$-1$, respectively, 
	contained in the upper half-plane.
	We will proceed by contradiction. 
	Suppose that the angle between~$u$ and~$v$ is different from~$\pi/2$. Without loss of generality we may assume that it is strictly below~$\pi/2$. 
	Define the ellipse
	$$
		Q = 	\{x\, u  +y \, v : x,y \in \bbR,\, x^2 + y^2 \le 1 \}.
%	\{r\, u \sin \theta +r \, v\cos \theta : \theta \in [0,2\pi),\, r\in [0,1) \}
	$$
	Since~$u,v,-u,-v$ lie on its boundary (see Figure~\ref{fig:ellipse}), denote by~$(u,v), \dots, (-v,u)$ the boundary arcs delimited by these points. 
	Notice that~$T_\alpha$ maps~$Q$ onto itself, with the points~$u$,~$v$,~$-u$ and~$-v$ mapped to~$u$,~$-v$,~$-u$ and~$v$, respectively. 

	As a consequence of the asymptotic invariance of~$\phi_{\bbL(\frac{\pi}2)}$ under~$T_\alpha$, 
	the probabilities of crossing~$Q$ from~$(u,v)$ to~$(-u,-v)$ and from~$(v,-u)$ to~$(-v,u)$ are asymptotically equal: 
	\begin{align}
				\big|\phi_{\delta\bbL(\frac{\pi}2)}\big[(u,v)\xlra{Q} (-u,-v)\big]-\phi_{\delta\bbL(\frac{\pi}2)}\big[(v,-u)\xlra{Q} (-v,u)\big]\big| 
		\xrightarrow[\delta \to 0]{} 0.
	\end{align}
	As~$\phi_{\bbL(\frac{\pi}2)}$ is also invariant under the rotation of angle~$\pi/2$ (which may be seen as~$S_{0} \circ S_{\pi/4}$), we conclude that 
	the probability of crossing~$Q$ from~$(u,v)$ to~$(-u,-v)$ is asymptotically close to that of its rotation by~$\pi/2$:
	\begin{align}
		\big|\phi_{\delta\bbL(\frac{\pi}2)}\big[(u,v)\xlra{Q} (-u,-v)\big] 
			- \phi_{\delta\bbL(\frac{\pi}2)}\big[(e^{{i}\frac{\pi}2}u,e^{{i}\frac{\pi}2}v)\xlra{e^{{i}\frac{\pi}2}Q} (-e^{{i}\frac{\pi}2}u,-e^{{i}\frac{\pi}2}v)\big]\big| 
			\xrightarrow[\delta \to 0]{} 0.
	\end{align}
	Combining the above, we conclude that 
	\begin{align}\label{eq:ellipse}
		\big|\phi_{\delta\bbL(\frac{\pi}2)}\big[(v,-u)\xlra{Q} (-v,u)\big]
			- \phi_{\delta\bbL(\frac{\pi}2)}\big[(e^{{i}\frac{\pi}2}u,e^{{i}\frac{\pi}2}v)\xlra{e^{{i}\frac{\pi}2}Q} (-e^{{i}\frac{\pi}2}u,-e^{{i}\frac{\pi}2}v)\big]\big| 
			\xrightarrow[\delta \to 0]{} 0.\quad
	\end{align}
	
	Observe now (see Figure~\ref{fig:ellipse}) that
	\begin{align}
	\big\{(e^{{i}\frac{\pi}2}u,e^{{i}\frac{\pi}2}v)\xlra{e^{{i}\frac{\pi}2}Q} (-e^{{i}\frac{\pi}2}u,-e^{{i}\frac{\pi}2}v)\big\} \subset 
		\big\{(v,-u)\xlra{Q} (-v,u) \big\}.
	\end{align}
	Furthermore,~\eqref{eq:RSW_iso} implies that 
	\begin{align}
		\liminf_{\delta \to 0}\phi_{\delta\bbL(\frac{\pi}2)}\big[(v,-u)\xlra{Q} (-v,u) \big] - \phi_{\delta\bbL(\frac{\pi}2)}\big[(e^{{i}\frac{\pi}2}u,e^{{i}\frac{\pi}2}v)\xlra{e^{{i}\frac{\pi}2}Q} (-e^{{i}\frac{\pi}2}u,-e^{{i}\frac{\pi}2}v)\big] >0.
	\end{align}	
	This contradicts~\eqref{eq:ellipse}, and therefore 
	invalidates our assumption that~$T_\alpha$ is not an orthogonal reflection.\smallskip
	
	\begin{figure}
    	\begin{center}
        	\includegraphics[width = 0.44\textwidth]{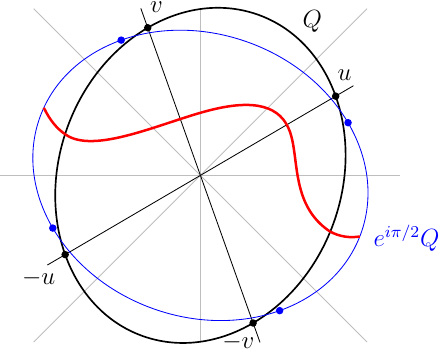}
        	\caption{If the angle between~$u$ and~$v$ is assumed to be strictly below~$\pi/2$, then~$Q$ is an ellipse (in black).
        	The rotation by~$\pi/2$ of~$(u,v)\xlra{Q} (-u,-v)$ is realised by the red path; it has strictly lower probability than~$(v,-u)\xlra{Q} (-v,u)$.}
        	\label{fig:ellipse}
    	\end{center}
	\end{figure}
	 
	We are now ready to conclude. 
	Write~$\theta(\alpha)$ for the angle between the real axis and the eigenvector~$M_{\frac{\pi}2,\alpha}^{-1} e^{{i}\alpha/2}$	of~$T_\alpha$. 
	Then, the composition of~$S_0$ with~$T_\alpha$ is the rotation of angle~$2\theta(\alpha)$, 
	and we conclude that~$\phi_{\bbL(\pi/2)}$ is asymptotically invariant under this rotation. 
	
	From  (ii) and (iii) of Proposition~\ref{prop:T}, we conclude that~$\{\theta(\alpha): \, \alpha \in (0,\pi)\}$ 
	has a non-empty interior, whence we deduce that~$\phi_{\bbL(\pi/2)}$ is asymptotically invariant under any rotation. 
	
This concludes the proof of Theorem~\ref{thm:rotation_invariance_infinite_vol} for the ${d}_{\rm CN}$ distance. The statement for the ${d}_{\rm SS}$ distance follows by~\cite[Sec. 2.3]{garban2013pivotal}. 
\end{proof}

\begin{remark}
	In the proof above, we used a special symmetry of the lattice~$\bbL(\pi/2)$, namely that with respect to the vertical reflection~$S_0$ 
	(or equivalently with respect to the rotation by~$\pi/2$). 
	The invariance of~$\bbL(\pi/2)$ with respect to~$S_{\pi/4}$, which corresponds to the generic invariance of~$\bbL(\alpha)$ with respect to~$S_{\alpha/2}$, 
	is not sufficient to conclude. 
	The best conclusion one could obtain without using the invariance under~$S_0$ is that the orbit of a given point under the group generated by~$\{T_\alpha: \, \alpha \in (0,\pi)\}$ 
	is an ellipse with axis~$e^{\pm \pi/4}\bbR$. 
	It is the vertical symmetry~(which is specific to~$\bbL(\pi/2)$ and has no correspondence for other  lattices~$\bbL(\alpha)$) 
	that allows us to conclude that this ellipse is actually a circle. 
\end{remark}

\subsection{Universality: proof of Theorem~\ref{thm:universalCNSS}}\label{sec:drift0}

We turn to the proof of Theorem~\ref{thm:universalCNSS}, or equivalently to the fact that~$M_{\frac{\pi}2,\alpha} = {\rm id}$ for all~$\alpha$. 
Recall that~$\phi_{\bbL(\beta)}$ is said to be asymptotically rotationally invariant if, for any~$\alpha \in [0,2\pi)$,
\begin{align*}
	{d}_{\rm CN}(\phi_{\delta \bbL(\beta)},\phi_{e^{i\alpha}\delta \bbL(\beta)}) \xrightarrow[\delta \to 0]{}  0. 
\end{align*}
Similarly, say that~$\phi_{\bbL(\alpha)}$ and~$\phi_{\bbL(\beta)}$ are {\em asymptotically similar} if 
\begin{align*}
	{d}_{\rm CN}(\phi_{\delta \bbL(\alpha)},\phi_{\delta \bbL(\beta)})\xrightarrow[\delta \to 0]{} 0. 
\end{align*}
Theorem~\ref{thm:linear} states that~$\phi_{\bbL(\alpha)}$ and~$\phi_{\bbL(\beta)}$ are asymptotically similar if~$M_{\beta,\alpha} = {\rm id}$.

The key to the proof of this section is the following lemma. 

\begin{lemma}\label{lem:rotation_inv_alpha}
	Fix~$\beta \in (0,\pi)$ and assume that~$\phi_{\bbL(\beta)}$ is asymptotically rotationally invariant. 
	Then~$M_{\beta,\beta/2}  = {\rm id}$. 
\end{lemma}

\begin{proof}
	Fix such a value of~$\beta$ and set~$\alpha = \beta/2$. 
	Then, due to~\eqref{eq:drift_RT} and to the special choice of~$\alpha$, 
	we have~${\rm Drift}_{\rm vert}(\beta,\beta/2) = {\rm Drift}_{\rm lat}(\beta,\beta/2)$. 
	When injected in~\eqref{eq:Mbetaalpha} and after basic computation, we find that~$M_{\beta,\alpha}$ has the special form 
	\begin{align*}
	M_{\beta,\beta/2} 
	= \begin{pmatrix}
	1 &  \frac{v}{\sin \beta/2}\\[8pt]
	0 & 1  + \frac{v}{\cos \beta/2}
	\end{pmatrix}
%\quad \text{ and }\quad 
%	M_{\beta,\alpha}^{-1}
%= \begin{pmatrix}
%1 & - \frac{v \cos \alpha}{\sin \alpha(v + \cos \alpha)}\\[8pt]
%0 &  \frac{\cos \alpha}{v + \cos \alpha}
%\end{pmatrix}
\qquad \text{where }v=  \frac{{\rm Drift}_{\rm lat}(\beta,\beta/2)(1 +2\cos\beta/2)}{2(\sin \beta/2 -{\rm Drift}_{\rm lat}(\beta,\beta/2))}.
\end{align*}
	In particular, we notice that the vectors~$1$ and~$e^{{i}\beta/2}$ are eigenvectors with eigenvalues~$1$ and~$\lambda:=1  + \frac{v}{\cos \beta/2} > 0$, respectively. 
	
	Recall from Proposition~\ref{prop:Tinv} that~$\phi_{\bbL(\beta)}$ is asymptotically invariant under~$T:= M_{\beta,\beta/2}^{-1} \circ S_{\beta/4} \circ M_{\beta,\beta/2}$.
	Since~$S_{\beta/4}$ exchanges the vectors~$1$ and~$e^{{i}\beta/2}$, 
	the transformation~$T$ maps the vector~$1$ to~$\lambda^{-1}e^{{i}\beta/2}$ and~$e^{{i}\beta/2}$ to~$\lambda$. 
	
	Consider the rhombic region~$R = \{x + y e^{{i}\beta/2}: \,x,y\in [-1,1]\}$ as a quad with four points~$a,b,c,d$ at its corners, 
	in counter-clockwise order, starting with the top left corner. By the asymptotic invariance of~$\phi_{\delta\bbL(\beta)}$ with respect to~$T$, 
	\begin{align*}
	\lim_{\delta\to0} \phi_{\delta\bbL(\beta)}\big[(ab) \xlra{R} (cd)\big] - \phi_{\delta\bbL(\beta)}\big[(T(a)T(b)) \xlra{T(R)} (T(c)T(d))\big]  = 0. 
	\end{align*}
	
	Now,~$\phi_{\delta\bbL(\beta)}$ is invariant under~$S_{\beta/2}$ and is assumed asymptotically rotationally invariant,
	which implies that it is asymptotically invariant under all reflections. In particular, it is asymptotically invariant under $S_{\beta/4}$. Since~$R$ is stable under~$S_{\beta/4}$, we conclude that
	\begin{align*}
	\lim_{\delta\to0} \phi_{\delta\bbL(\beta)}\big[(ab) \xlra{R} (cd)\big] - \phi_{\delta\bbL(\beta)}\big[(bc) \xlra{R} (da)\big]  = 0. 
	\end{align*}
	Combining the above, we find that 
	\begin{align*}
	\lim_{\delta\to0} \phi_{\delta\bbL(\beta)}\big[(bc) \xlra{R} (da)\big] - \phi_{\delta\bbL(\beta)}\big[(T(a)T(b)) \xlra{T(R)} (T(c)T(d))\big]  = 0. 
	\end{align*}
	
	\begin{figure}
	\begin{center}
	\includegraphics[width = 0.55\textwidth]{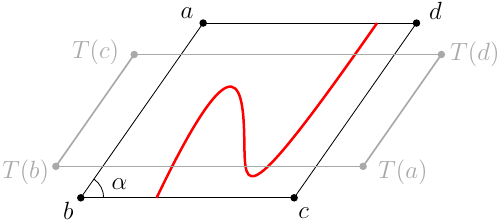}
	\caption{The rhombi~$R$ and~$T(R)$ under the assumption that~$\lambda> 1$. Any crossing in~$R$ between the arcs~$(bc)$  and~$(da)$ produces a crossing in~$T(R)$ between~$(T(a)T(b))$ and~$(T(c)T(d))$. However, the latter crossing is strictly easier to achieve, uniformly in the scale~$\delta$ of the lattice.}
	\label{fig:two_rhombi}
	\end{center}
	\end{figure}
This is only possible if~$\lambda = 1$.
	Indeed, the quad~$T(R)$ may easily be identified as $\{\lambda y + x \lambda^{-1} e^{{i}\alpha}: \,x,y\in [-1,1]\}$.
	If~$\lambda > 1$, then~$T(R)$ is wider and shorter than~$R$ (see Figure~\ref{fig:two_rhombi}) and~\eqref{eq:RSW_iso} implies that 
\begin{align*}
	\limsup_{\delta\to0} \phi_{\delta\bbL(\beta)}\big[(bc) \xlra{R} (da)\big] - \phi_{\delta\bbL(\beta)}\big[(T(a)T(b)) \xlra{T(R)} (T(c)T(d))\big]  < 0. 
\end{align*}
	Conversely, if~$\lambda < 1$,~\eqref{eq:RSW_iso} shows that
\begin{align*}
	\liminf_{\delta\to0} \phi_{\delta\bbL(\beta)}\big[(bc) \xlra{R} (da)\big] - \phi_{\delta\bbL(\beta)}\big[(T(a)T(b)) \xlra{T(R)} (T(c)T(d))\big]  > 0. 
\end{align*}	
	We conclude that~$\lambda = 1$, which translates into $v=0$ and hence~$M_{\beta,\beta/2}  = {\rm id}$. 	
\end{proof}

\begin{proof}[Proof of Theorem~\ref{thm:universalCNSS}]
	Define the set 
	$$\calR = \{\alpha \in (0,\pi) : \,\phi_{\bbL(\alpha)} \text{ asymptotically similar to~$\phi_{\bbL(\pi/2)}$}\}.$$

	Theorem~\ref{thm:rotation_invariance_infinite_vol} implies that $\phi_{\bbL(\pi/2)}$ is asymptotically rotationally invariant. As a consequence of the definition of $\calR$, so is $\phi_{\bbL(\alpha)}$ for any~$\alpha \in \calR$.
	In light of Lemma~\ref{lem:rotation_inv_alpha},~$\calR$ is stable by~$\alpha \mapsto \alpha/2$. 
	Moreover, by horizontal symmetry,~$\calR$ is stable by~$\alpha \mapsto \pi- \alpha$. 
	Repeatedly applying  these transformations shows that~$\calR$ is dense in~$(0,\pi)$. 
	
	Fix~$\alpha \in \calR$. 
	That~$\phi_{\bbL(\pi/2)}$ and~$\phi_{\bbL(\alpha)}$ are asymptotically similar yields~$M_{\frac{\pi}2, \alpha} = {\rm id}$. 
	Due to the continuity of~$\alpha \mapsto M_{\frac{\pi}2,\alpha}$ (see Lemma~\ref{lem:M_cont})
	and to~$\calR$ being dense in~$(0,\pi)$, we conclude that~$M_{\frac{\pi}2, \alpha} = {\rm id}$ for all~$\alpha \in (0,\pi)$. 
	
	Apply Proposition~\ref{prop:linear_unif} (that is, the uniform version of Theorem~\ref{thm:linear}) to conclude Theorem~\ref{thm:universalCNSS} for the ${d}_{\rm CN}$ distance. 
	The statement for the ${d}_{\rm SS}$ distance follows by~\cite[Sec. 2.3]{garban2013pivotal}. 
\end{proof}

\section{Consequences of Theorem~\ref{thm:rotation_invariance_infinite_vol}}\label{sec:consequences}

The proofs in this section are consequences of Theorem~\ref{thm:rotation_invariance_infinite_vol} and only use random-cluster model techniques on $\bbZ^2$. 
They are independent of the arguments appearing in the proof of Theorems~\ref{thm:rotation_invariance_infinite_vol},~\ref{thm:universalCNSS} and~\ref{thm:linear}. 

We start with the proof of Corollary~\ref{cor:rotation_invariance_domain}. We distinguish the case when the domain is bounded from that of unbounded domains. We only sketch the case of bounded domains; it is proved in greater detail in \cite{He26}.

\begin{proof}[Proof sketch of Corollary~\ref{cor:rotation_invariance_domain} when $\Omega$ is bounded] 
To obtain the result in a finite bounded domain, we use the domain Markov property and the fact that one may approximate~$\phi_{\Omega_\delta}^0$ by asking that there exists a loop~$\Gamma$ within distance~$\eta$ of~$\partial\Omega$ in the infinite-volume measure. More precisely, let~$A(\Omega,\eta)$ be the event that there exists a loop~$\pmb\Gamma\in \calF_0(\omega_\delta)$ which is included in~$\Omega$ and such that~$d(\pmb\Gamma,\partial\Omega)\le \eta$ ($d$ is the distance between loops defined in the introduction). Note that whether~$A(\Omega,\eta)$ occurs or not can be measured in the Schramm-Smirnov topology (we leave this as an exercise). 

Now, fix~$\ep_0>0$. We use the characterization of the Schramm-Smirnov distance provided in \cite[Proposition~3.9]{GPS}. There exists a family of non-degenerate quads~$Q_1,\dots,Q_n$ in~$\Omega$ such that if the sets of quads in~$Q_1,\dots,Q_n$ that are crossed are the same in~$\omega_\delta$ and~$\omega'_\delta$, then~$d_{\rm SS}(\omega_\delta,\omega'_\delta)\le \ep_0$. In particular, we deduce that if 
$H_{\vec Q}(I)$ denotes the event that~$Q_i$ is crossed if and only if~$i\in I$, then there exists a coupling~$\mathbf P$ between~$\omega_\delta\sim \phi_{\Omega_\delta}^0$ and~$\omega'_\delta\sim \phi_{e^{{\rm i}\alpha}\Omega_\delta}^0$ such that 
 \[
 \mathbf P[d_{\rm SS}(\omega_\delta,\omega'_\delta)\ge \ep_0]\le \ep_0
 \]
 if for every~$I\subset \{1,\dots,n\}$,
 \begin{align}\label{eq:ah1}
\big|\phi_{\Omega_\delta}^0[H_{\vec Q}(I)]-\phi_{e^{{\rm i}\alpha}\Omega_\delta}^0[H_{e^{{\rm i}\alpha}\vec Q}(I)]\big|\le \ep_0/2^n=:\ep.
\end{align}

We will now endeavour to prove \eqref{eq:ah1} for $\delta$ sufficiently small. 
Theorem~\ref{thm:rotation_invariance_infinite_vol} implies that for every~$\delta<\delta_0(\Omega,\eta,\ep)$,
\begin{align}
\big|\phi_{\delta\bbZ^2}[H_{\vec Q}(I)|A(\Omega,\eta)]-\phi_{\delta\bbZ^2}[H_{e^{{\rm i}\alpha}\vec Q}(I)|A(e^{{\rm i}\alpha}\Omega,\eta)]\big|\le\tfrac13\ep.
\end{align}
We therefore wish to prove that by choosing $\eta$ small enough, one gets
\begin{align}\label{eq:ahahah1}
\big|\phi_{\Omega_\delta}^0[H_{\vec Q}(I)]-\phi_{\delta\bbZ^2}[H_{\vec Q}(I)|A(\Omega,\eta)]\big|\le \tfrac13\ep
\end{align}
for all $\delta$  sufficiently small. 
The same can be done for the rotated version, so that the previous displayed equations imply~\eqref{eq:ah1} and conclude the proof.

To show~\eqref{eq:ahahah1}, let~$\pmb\Omega_\delta$ be the interior of the outer-most loop in~$\calF_0(\omega)$ satisfying the conditions of~$A(\Omega,\eta)$. Using the spatial Markov property, it suffices to show that 
\begin{align}
\big|\phi_{\Omega_\delta}^0[H_{\vec Q}(I)]-\phi_{\pmb\Omega_\delta}^0[H_{\vec Q}(I)]\big|\le \tfrac13\ep.
\end{align} 
Note that there is an increasing coupling between~$\pmb\omega_\delta\sim\phi_{\pmb\Omega_\delta}^0$ and~$\omega_\delta\sim\phi_{\Omega_\delta}^0$ ($\pmb\omega_\delta\le \omega_\delta$ because of~$\pmb\Omega_\delta\subset\Omega_\delta$), so that for~$\omega_\delta$ to belong to~$H_{\vec Q}(I)$ but not~$\pmb\omega_\delta$ or vice versa, it must be that one of the quads~$Q_i$ must be crossed in one but not in the other. We deduce that it suffices to show that for every possible realization of~$\pmb\Omega_\delta$,
\begin{align}
\phi_{\Omega_\delta}^0[\calC(Q_i)]-\phi_{\pmb\Omega_\delta}^0[\calC(Q_i)]\le \tfrac{1}{3n}\ep.
\end{align} 
Therefore, the result boils down to the following claim, which we prove separately. 
\end{proof}

\begin{claim}\label{claim:omom'}
For every~$\epsilon>0$, every bounded simply connected domain~$\Omega$ with~$C^1$-smooth boundary, and every quad~$Q$ inside~$\Omega$,  there exists~$\eta=\eta(\Omega,Q,\epsilon)>0$ such that for every~$\Omega'\subset\Omega$ with~$d(\partial\Omega',\partial\Omega)\le \eta$,
\[
	\phi_{\Omega_\delta}^0[\calC(Q)]\le \phi_{\Omega'_\delta}^0[\calC(Q)]+ \epsilon
\]
for~$\delta$ small enough.
\end{claim}

\begin{figure}
	\begin{center}
	\includegraphics[width=.9\textwidth]{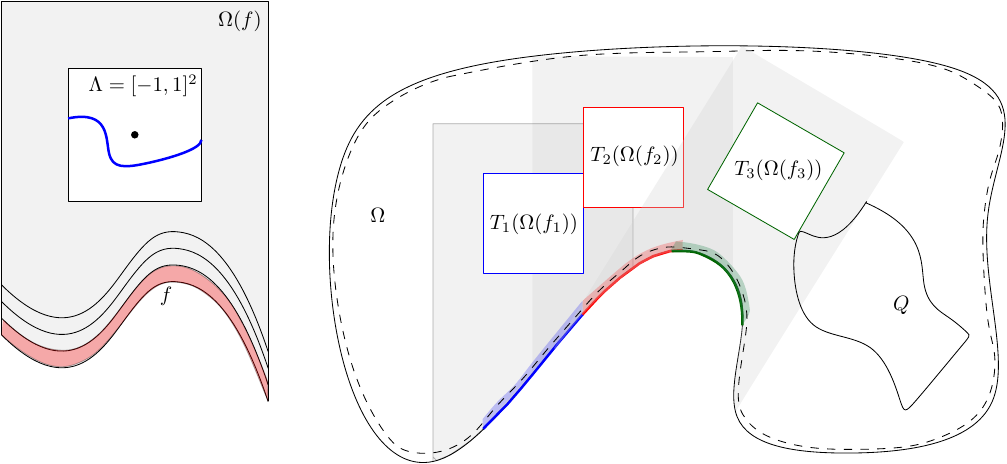}
	\caption{{\em Left:} an example of a domain~$\Omega(f)$. Note that the sets~$\Omega(f+k\alpha)$ have a nested structure (the red part denotes~$\Omega(f)\setminus\Omega(f+\alpha)$). {\em Right:} the impact of changing the boundary is compared with the impact of changing the boundary in a family of subdomains which are images under simple transformations of domains~$\Omega(f)$ (with potentially different functions~$f$). The existence of such a decomposition is made possible by the fact that the boundary of~$\Omega$ is~$C^1$.}
	\label{fig:finitedomainproof}
	\end{center}
\end{figure}

\begin{proof}[Proof of Claim~\ref{claim:omom'}]
    We only sketch the proof. Consider first the ``epigraph`` domains indexed by continuous functions~$f$ from~$[-2,2]$ to~$\bbR$ given by
     \[
     \Omega(f):=\{x=(x_1,x_2)\in \bbR^2:x_1\in(-2,2), f(x_1)<x_2<2\}
     \]
    (see Figure~\ref{fig:finitedomainproof}). Define~$\Lambda:=[-1,1]^2$. 
    For~$\alpha>0$, a quite lengthy application of the techniques developed in \cite[Lemma 5.3]{DumMan20} (see \cite{He26} for details) implies that for every~$f\le -2$ and~$1\le k\le \tfrac12\lfloor 1/\alpha\rfloor=:K$, 
    \[
    \phi^0_{\Omega(f)_\delta}[\calC(\Lambda)]-\phi^0_{\Omega(f+\alpha)_\delta}[\calC(\Lambda)]\le C\big(\phi^0_{\Omega(f+k\alpha)_\delta}[\calC(\Lambda)]-\phi^0_{\Omega(f+(k+1)\alpha)_\delta}[\calC(\Lambda)]\big).
    \]
    Summing over~$1\le k\le K$, we deduce that 
    \begin{align}\label{eq:ahahah5}
    \phi^0_{\Omega(f)_\delta}[\calC(\Lambda)]-\phi^0_{\Omega(f+\alpha)_\delta}[\calC(\Lambda)]
    \le \tfrac{C}{K}\big(\phi^0_{\Omega(f)_\delta}[\calC(\Lambda)]-\phi^0_{\Omega(f+K\alpha)_\delta}[\calC(\Lambda)]\big)
    \le \tfrac{C}{K}\le 4C\alpha. \quad
    \end{align}
    Note that a similar argument works for any rotation, translate, or rescaling of the domains above.
    
    We now use our assumption that~$\partial\Omega$ is~$C^1$-smooth. Since~$\partial\Omega$ is given by a curve~$\gamma$ which is~$C^1$ and has non-vanishing differential, one may find (see Figure~\ref{fig:finitedomainproof}) constants~$\kappa=\kappa(\Omega)>0$ and~$C=C(\Omega)>0$,  functions~$f_s:[-2,2]\rightarrow (-\infty,-2]$ and~$T_s:\bbR^2\rightarrow\bbR^2$ for~$1\le s\le S$, where~$S$ depends on~$\Omega$ (through the modulus of continuity of the derivative for the function parametrizing~$\partial\Omega$) but not on~$\eta$, satisfying the following properties:
    \begin{itemize}[noitemsep]
    \item~$T_s$ is the composition of a rotation, a translation, and the multiplication by~$\kappa$;
    \item~$T_s(\Omega(f_s))$ is included in~$\Omega$ for every~$s$;
    \item for all~$\eta$ small enough,~$\{x\in\Omega:d(x,\Omega^c)\le \eta\}$ is included in the union of the sets
    \[ A_s:=T_s(\{x=(x_1,x_2):x_1\in [-1,1], f(x_1)<x_2<f(x_1)+C\eta\}).\]
    \end{itemize}
    Introducing the domains
    $ \Omega_s:=\Omega\setminus \bigcup_{t=1}^s A_t$, and using again \cite{DumMan20} for the first and second inequalities, one can prove the existence of~$C_i=C_i(\Omega,Q,\kappa)>0$ such that 
     \begin{align}
     \phi_{\Omega_{s-1}}^0[\calC(Q)]- \phi_{\Omega_{s}}^0[\calC(Q)]&\le C_1\big( \phi_{\Omega_{s-1}}^0[\calC(T_s(\Lambda))]-\phi_{\Omega_s}^0[\calC(T_s(\Lambda))]\big)\nonumber\\
    % &\le C_2\big(\phi_{T_s(\Omega(f_s))}^0[\calC(T_s(\Lambda))]-\phi_{T_s(\Omega(f_s+C\eta))}^0[\calC(T_s(\Lambda))]\big)\nonumber\\
     &\le C_2\big(\phi_{\Omega(f_s)}^0[\calC(\Lambda)]-\phi_{\Omega(f_s+C\eta)}^0[\calC(\Lambda)]\big)\label{eq:ahahah4}\\
     &\le C_3\eta,\nonumber
     \end{align}
    where the last line is due to~\eqref{eq:ahahah5} applied to~$\alpha=C\eta$.
    
    Choose~$\eta=\eta(\Omega,\epsilon,S)>0$ small enough. Summing~\eqref{eq:ahahah4} over~$s$ gives
    \[
    \phi_{\Omega_\delta}^0[\calC(Q)]-\phi^0_{\Omega_\delta'}[\calC(Q)]\le \sum_{s=1}^{S}\phi_{\Omega_{s-1}}^0[\calC(Q)]- \phi_{\Omega_{s}}^0[\calC(Q)]\le \epsilon, 
    \]
    as claimed. 
 \end{proof}

\begin{proof}[Proof of Corollary~\ref{cor:rotation_invariance_domain} when $\Omega$ is unbounded] 
	For every~$\ep>0$, to determine the Schramm-Smirnov distance up to a precision of~$\ep>0$, only quads in~$B(0,1/\ep)$ need to be considered. Consider a bounded domain~$\Omega^{(\ep)}$ that coincides with~$\Omega$ on~$B(0,1/\ep^C)$. By the mixing property, one has that for every~$\delta>0$ and every event~$E$ depending on edges in~$\delta\bbZ^2\cap B(0,1/\ep)$ only, 
\[
|\phi_{\Omega^{(\ep)}_\delta}[E]-\phi_{\Omega_\delta}[E]|\le C_{\rm  mix}\ep^{\rm c_{\rm mix}(C-1)}\phi_{\Omega_\delta}[E].
\]
Using the invariance by rotation in~$\Omega^{(\ep)}_\delta$ and taking~$\delta$ to 0 then~$\ep$ to 0 concludes the proof.%Picking~$C>0$ large enough and using the result for~$\Omega^{(\ep)}$, we obtain this case as well.
\end{proof}

\begin{proof}[Proof of Corollary~\ref{cor:crossing}]    
	 When one considers a quad~$Q$ that remains at a distance at least~$\ep$ from the boundary of~$\Omega$, the result follows directly from   Corollary~\ref{cor:rotation_invariance_domain} %Theorem~\ref{thm:rotation_invariance_infinite_vol} 
	 and the measurability of~$\mathcal C(Q)$ in the Schramm-Smirnov topology.% (note that the event gets rewritten as~$Q\in \omega$ when~$\omega$ is seen as an element of~$\calH$). 
    
    Now, when~$1\le q<4$, to get the result without any assumption on the distance to the boundary, note that for a quad~$Q$, there exists a quad~$Q'$ whose distance to~$\partial\Omega$ is at least~$\ep$, and which is at Hausdorff distance at most~$2\ep$ from~$Q$. Using the strong version of crossing estimates from \cite{DumManTas20}, we easily obtain  (this type of reasoning is now classical, see for instance \cite[Lemma 3.12]{DumMan20}) that 
    \[
 \big|\phi_{\Omega_\delta}[\mathcal C(Q)]-\phi_{\Omega_\delta}[\mathcal C(Q')]  \big|\le C\ep^c
    \]
    for universal constants~$c,C>0$. The result follows readily by first choosing~$\ep$ small enough and then letting~$\delta$ tend to zero and using the rotational invariance result for~$Q'$.\end{proof}
    
    \begin{proof}[Proof of Corollary~\ref{cor:connectivity_correlations}] 
    We use a conditional mixing argument due to Garban, Pete, and Schramm \cite[Section 3]{garban2013pivotal} in the case of Bernoulli percolation and that can be extended to the random-cluster model using RSW theory. Consider the {\em Euclidean} ball~$B_n$ of radius~$n$, and its boundary~$\partial B_n$. Introduce the quantities
    \[
    \epsilon(n,N):=\phi_{\mathbb Z^2}^0[0\longleftrightarrow B_N^c|B_n\longleftrightarrow B_N^c]\quad\text{and}\quad\epsilon(n):=\lim_{N\rightarrow\infty}\epsilon(n,N).
    \]
 (The existence of $\epsilon(n)$ follows from the mixing between scales in the measure conditioned on the one arm event, which in turn follows from the RSW theory.)   The statement of conditional mixing from \cite{garban2013pivotal} implies the following claim (in \cite{garban2013pivotal} it is stated for the four-arm event, but a similar -- in fact simpler -- argument can be performed for the one-arm event, see e.g.~Proposition 5.3 of the same paper). For every~$\beta,\ep>0$, there exists~$\eta=\eta(\beta,\ep)>0$ such that for every~$\Omega$ and every~$x_1,\dots,x_n$ at a distance~$\ep$ from each other and from the boundary, and every partition~$P$ of~$(x_1,\dots,x_n)$,
    \[
    \Big|\phi_{\Omega_\delta}^0[\calE(P,x_1,\dots,x_n)]-\epsilon({\tfrac\eta\delta})^n\phi_{\Omega_\delta}^0[\calE(P,B_{\eta/\delta}(x_1),\dots,B_{\eta/\delta}(x_n))]\Big|\le \beta \phi_{\Omega_\delta}^0[\calE(P,x_1,\dots,x_n)],
    \]
    where~$\calE(P,B_{\eta/\delta}(x_1),\dots,B_{\eta/\delta}(x_n))$ is the event that the balls~$B_{\eta/\delta}(x_i)$ are connected to each other if and only if they belong to the same element of the partition~$P$. The same formula applies in the rotated measure. 
    
    We conclude, by observing that~$e^{{\rm i}\alpha} B_{\eta/\delta}(x_i)$ and~$B_{\eta/\delta}(e^{{\rm i}\alpha}x_i)$ are equal, and that the event~$\calE(P,B_{\eta/\delta}(x_1),\dots,B_{\eta/\delta}(x_n))$ is measurable in the Schramm-Smirnov topology, so that its probability and the probability of its rotation by an angle $\alpha$ are close to each other by Corollary~\ref{cor:rotation_invariance_domain}.% Theorem~\ref{thm:rotation_invariance_infinite_vol}.
\end{proof}

\begin{proof}[Proof of Corollary~\ref{cor:Potts_correlations}]
Fix~$\tau_1,\dots,\tau_n\in \mathbb T_q$. For each spin $\sigma$, let~$I_\sigma\subset\{x_1,\dots,x_n\}$ be the set of~$x_j$ such that~$\tau_{j}=\sigma$. Call a partition~$P=(P_1,\dots,P_k)$ of~$\{x_1,\dots,x_n\}$ {\em compatible with~$\tau$} if each~$P_j$ is included in one of the~$I_\sigma$. Also, let~$|P|=k$ be the number of elements in the partition. The Edwards-Sokal coupling implies that 
\[
\mu_{\Omega_\delta}[\sigma_{x_i}=\tau_i,1\le i\le n]=\sum_{\text{compatible }P}q^{-|P|}\phi^0_{\Omega_\delta}[\mathcal E(P,x_1,\dots,x_n)].
\]
We deduce the corollary by using Corollary~\ref{cor:connectivity_correlations}. 
\end{proof}

\paragraph{Acknowledgments.} We thank E.~Averous, L.~Gassmann, U.T.~Hansen, T.~He, P.~Lammers and M.~Mohanarangan for their comments on the paper. The first author was funded by the ERC CriBLaM. The third and fourth authors were funded by the Swiss FNS. The first, third, fourth and fifth authors were partially funded by the NCCR SwissMap and the Swiss FNS.

\bibliographystyle{alpha}
\newcommand{\etalchar}[1]{$^{#1}$}

\end{document}